\documentclass[journal ]{new-aiaa}
\usepackage[utf8]{inputenc}
\usepackage{textcomp}

\usepackage{graphicx}
\usepackage{caption}
\usepackage{amsmath}
\usepackage{amssymb}
\usepackage[version=4]{mhchem}
\usepackage{siunitx}
\usepackage{longtable,tabularx}
\usepackage{latexsym}
\usepackage{titlesec}
\usepackage{mathrsfs}
\usepackage{epstopdf}
\usepackage{subfigure}
\usepackage{algorithm}
\usepackage{algorithmic}
\usepackage{setspace}
\usepackage{bm}
\usepackage{amsthm}
\usepackage{multicol}
\usepackage{multirow}
\usepackage{booktabs}
\usepackage{color,xcolor}
\usepackage{listings}
\definecolor{codegreen}{rgb}{0,0.6,0}
\title{Investigation of discontinuous Galerkin methods in adjoint gradient-based aerodynamic shape optimization}

\author{Yiwei Feng \footnote{Postdoctoral Researcher. fengyw@buaa.edu.cn (Corresponding Author).} }
\affil{China Academy of Aerospace Aerodynamics, Beijing, 100074, P.R. China.}
\author{Lili Lv and Tiegang Liu \footnote{Professor. liutg@buaa.edu.cn (Corresponding Author).} and Kun Wang}
\affil{LMIB and School of Mathematical Sciences, Beihang University, Beijing, 100191, P.R. China.}
\author{Bangcheng Ai}
\affil{China Academy of Aerospace Aerodynamics, Beijing, 100074, P.R. China.}

\begin{document}

\maketitle

\begin{abstract}
This work develops a robust and efficient framework of the adjoint gradient-based aerodynamic shape optimization (ASO) using high-order discontinuous Galerkin methods (DGMs) as the CFD solver. The adjoint-enabled gradients based on different CFD solvers or solution representations are derived in detail, and the potential advantage of DG representations is discovered that the adjoint gradient computed by the DGMs contains a modification term which implies information of higher-order moments of the solution as compared with finite volume methods (FVMs). A number of numerical cases are tested for investigating the impact of different CFD solvers (including DGMs and FVMs) on the evaluation of the adjoint-enabled gradients. The numerical results demonstrate that the DGMs can provide more precise adjoint gradients even on a coarse mesh as compared with the FVMs under coequal computational costs, and extend the capability to explore the design space, further leading to acquiring the aerodynamic shapes with more superior aerodynamic performance.
\end{abstract}

\section*{Nomenclature}

\noindent(Nomenclature entries should have the units identified)

{\renewcommand\arraystretch{1.0}
\noindent\begin{longtable*}{@{}l @{\quad=\quad} l@{}}
$A$                  &  area of airfoil \\
AOA                  &  angle of attack \\ 
$\bm{c}_E$           &  equality constraints' vector \\
$\bm{c}_I$           &  inequality constraints' vector \\
$C_d$                &  drag coefficient \\
$C_{d0}$             &  initial drag coefficient \\
$C_{d1}$             &  optimized drag coefficient \\
$C_l$                &  lift coefficient \\
$C_p$                &  pressure coefficient \\
$\bm{D}$             &  design variables' vector \\
$\mathcal{D}$        &  design space \\
$E$                  &  total energy of the fluid \\
$\bm{\mathcal{F}}_I$ &  inviscid fluxes of the Euler equation \\ 
$J$                  &  aerodynamic objective \\
Ma                   &  Mach number  \\
$N_b$                &  number of vertexes at boundary \\
$N_d$                &  dimension of design space \\ 
$N_e$                &  number of computational elements \\
$\bm{R}$             &  (discretized) left hand side of the governing PDE constraint \\
$p$                  &  pressure of the fluid \\
$\bm{u}$             &  conservative variable vector of the Euler equation \\
$\bm{U}$             &  numerical solution variables' vector \\
$\bm{U}_h$           &  global flow variables' vector \\
$V$                  &  volume of wing \\ 
$\bm{v}$             &  velocity vector of the fluid \\
$\bm{X}$             &  computational grid vertexes' vector \\
$\bm{X}_{vol}$       &  inner grid vertexes' vector \\
$\bm{X}_{surf}$      &  surface grid vertexes' vector \\
$\rho$               &  density of the fluid \\
$\gamma$             &  ratio of the specific heats \\
$\bm{\lambda}$       &  adjoint solution variables' vector \\ 
\end{longtable*}}

\section{Introduction}
\lettrine{B}enefiting from the increasing power of modern computers, aerodynamic shape optimization (ASO) based on Computational Fluid Dynamics (CFD) has become an indispensable component for effective and robust modern aerodynamic designs \cite{jameson1998optimum, carpentieri2007adjoint, chernukhin2013multimodality, martins2013multidisciplinary, skinner2018state}. CFD solver is a particularly crucial chain in ASO, on the one hand, CFD solver needs to be high-fidelity to ensure the high-accuracy evaluation of aerodynamic performance, on the other hand, CFD solver needs to be high-efficiency to avoid massive time cost during the optimization phase.

Over the past decades, high-order numerical methods, especially the high-order discontinuous Galerkin methods (DGMs) \cite{bassi2002numerical, fidkowski2005p, hartmann2010discontinuous, wang2013high}, have emerged as a competitive alternative in solving various CFD problems, and therefore have become apparent as CFD solvers in ASO. For example, Wang et al. \cite{wang2011adjoint, wang2012shape} applied DGMs in adjoint-based aerodynamic shape optimization; Zahr et al. \cite{zahr2016high} applied DG solver of the Navier-Stokes equations for time-dependent aerodynamic optimization; Wang et al. \cite{wang2021local} used DGMs in simultaneous perturbation stochastic approximation (SPSA) algorithm for aerodynamic optimization. In order to further acquire higher accuracy of the aerodynamic objective and constraints, Li and Hartamann \cite{li2015adjoint} firstly developed higher-order adaptive DG methods with control of the discretization error, and employed them into discrete adjoint optimization algorithm for minimizing drag of an inviscid transonic flow around the RAE2822 airfoil. Chen and Fidkowski \cite{chen2019discretization} developed an adaptive DGM for controlling errors of both objective and constraint outputs as well, and demonstrate its superiority in 2D drag minimization with a target lift coefficient. Later, Wang et al. \cite{wang2019adjoint} and Pezzano et al. \cite{pezzano2022geometrically} developed isogeometric adaptive DGMs, and applied them into adjoint-based optimization and Bayesian (BO) optimization which is also denoted surrogate model-based optimization, respectively. Recently, Coppeans et al. \cite{coppeans2023comparison} shown the strengths of DGMs in improving solution accuracy even on a coarse mesh as compared with FVMs, and provide an attractive alternative of the CFD solver in ASO.

In addition to high-accuracy objective value, the gradient-based optimization algorithms also have high requirements for gradients/derivatives of the objective and constrains. A high-confidence gradient is crucial for complex nonlinear optimization with massive design variables, and might extend the capability to explore the design space. For ASO, the evaluation of adjoint-enabled gradients is closely related to many intermediate variables and strategies \cite{kaland2015adaptive, ugolotti2019adjoint}. For example, the works of Li, Hartmann \cite{li2015adjoint} and Chen, Fidkowski \cite{chen2019discretization} have shown that the change of computational grid (through adaptive mesh refinement) results in different evaluated gradients; Luca et al. \cite{abergo2023aerodynamic} verified that even if using the same initial computational grid, different mesh deformation strategies (ELA and RBF based methods) might still lead to different gradients. Up to now, few authors analyzed whether high-order representations or numerical methods such as DGMs have strengths on adjoint-enabled gradient evaluation under coequal degrees of freedom (DoFs) and computational costs. In addition, since the DGMs are sometimes recognized lack of robustness and expensive in terms of computational costs for simulating high-speed flows, few authors applied DGMs as the CFD solver to the 3D ASO.

The major objective of the effort presented in this work is threefold: (a) based on our recently developed platform HODG version-1.0 \cite{he2023hodg}, a robust and efficient discrete adjoint-based optimization algorithm is developed for 3D aerodynamic shape design; (b) the discrete adjoint-enabled gradients based on different numerical frameworks are derived in detail, and the potential advantages of DG representation compared with finite volume representation are discovered and analyzed; (c) a number of 2D and 3D numerical experiments are designed to support the theoretical analysis, and the results can show superiority of using DGMs in terms of acquiring more precise adjoint gradient and further exploring better final aerodynamic shapes under similar computational costs in the adjoint gradient-based optimization algorithm.

The remainder of this paper is organized as follows. The problem description of aerodynamically orientated optimization is shown in Section II. The general workflow of the adjoint-based optimization and its detailed implementation based on the HODG platform are presented in Section III. Sensitivity analysis based on adjoint gradient is made in Section IV, and the impact and advantages of high-order DG representations on gradient evaluation are analyzed and emphasized as well. In Section V, numerical experiments are designed to verify the theoretical remarks in Section IV. Concluding remarks and perspectives are made in Section VI.

\section{Problem description of aerodynamic shape optimization}

An aerodynamically orientated optimization problem can be mathematically described as a partial differential equation constrained (PDE-constrained) optimization of an aerodynamically relevant objective $J$,
\begin{equation}
\begin{aligned}
\label{eq:opt}
 & \min_{\bm{D}\in\mathcal{D}}\;\;J(\bm{U},\bm{D}) \\
 s.t. \quad & \left\{
\begin{array}{lr}
      \bm{R}(\bm{U},\bm{D})=\bm{0}\\
      \bm{c}_E(\bm{U},\bm{D})=\bm{0} \\
      \bm{c}_I(\bm{U},\bm{D})\leq \bm{0}.
\end{array} \right.
\end{aligned}
\end{equation}

Here, the objective $J=J(\bm{U},\bm{D})$ is usually an aerodynamic coefficient (such as drag) which is associated with $\bm{D}$ and $\bm{U}$. $\bm{D}\in \mathcal{D}\subseteq \mathbb{R}^{N_d}$ is called the design variables' vector, which in general determines the aerodynamic shape, and $\bm{U}$ is called the solution or state variables' vector, which represents the state of fluid affecting the objective.

In (\ref{eq:opt}), $\bm{R}(\bm{U},\bm{D})=\bm{0}$ are PDE constraints which denote a fluid governing equation system (such as the compressible Euler or Navier-Stokes equations), thus the flow variables' vector $\bm{U}$ is the solution vector to the PDE determined by the design variable $\bm{D}$ and free stream, which can be implicitly expressed as $\bm{U}=\bm{U}\left(\bm{D}\right)$. $\bm{c}_E(\bm{U},\bm{D})=\bm{0}$ and $\bm{c}_I(\bm{U},\bm{D})\leq \bm{0}$ are algebraic equality and inequality constraints, respectively, such as the requirement of conservation of lift or non-decrease of aero-thickness.

For that it is quite difficult to obtain the explicit formulation of $\bm{U}=\bm{U}(\bm{D})$ governed by $\bm{R}(\bm{U},\bm{D})=\bm{0}$, grid variables' vector $\bm{X}$ is introduced as the intermediate variable, the problem (\ref{eq:opt}) now is rewritten as the following form,
\begin{equation}
\begin{aligned}
\label{eq:opt2}
 & \min_{\bm{D}\in\mathcal{D}}\;\;J\big(\bm{U},\bm{X}(\bm{D})\big) \\
 s.t. \quad & \left\{
\begin{array}{lr}
      \bm{R}(\bm{U},\bm{X}(\bm{D}))=\bm{0}\\
      \bm{c}_E(\bm{U},\bm{X}(\bm{D}))=\bm{0} \\
      \bm{c}_I(\bm{U},\bm{X}(\bm{D}))\leq \bm{0}.
\end{array} \right.
\end{aligned}
\end{equation}

Here, $\bm{X}=\bm{X}(\bm{D})$ contains geometric information of meshes (e.g. coordinates of gird vertices), and is generated or deformed according to the shape $\bm{D}$. $\bm{U}=\bm{U}(\bm{X})$ is the numerical solution vector of $\bm{R}(\bm{U},\bm{X})=\bm{0}$ based on the grid $\bm{X}$.

$\bm{X}$ is usually based on the initial grid and deformed according to the morph of the $\bm{D}$, thus $\bm{X}=\bm{X}(\bm{D})$ can be expressed by
\begin{equation}
\bm{X}(\bm{D})=\bm{X}_{vol}\left( \bm{X}_{surf}(\bm{D})\right),
\label{eq:XD}
\end{equation}
where $\bm{X}_{vol}$ denotes the inner volume grid, the form of $\bm{X}_{vol}\left( \bm{X}_{surf}\right)$ is provided by a specific mesh deformation method, $\bm{X}_{surf}$ denotes the surface grid, the form of $\bm{X}_{surf}(\bm{D})$ is described by a specific shape parameterization method.

$\bm{U}$ is implicitly represented through solving the discretized governing equations,
\begin{equation}
 \bm{R}(\bm{U},\bm{X}(\bm{D}))=\bm{0}\quad\Longrightarrow\quad \bm{U} = \bm{U}\left( \bm{X}(\bm{D}) \right).
\label{eq:UX}
\end{equation}

In general, the objective $J$ is directly associated with the global flow variable $\bm{U}_h$, which is defined by the numerical solution vector $\bm{U}$ and the grid $\bm{X}$,
\begin{equation}
\bm{U}_h = \bm{U}_h(\bm{U},\bm{X}).
\label{eq:flow}
\end{equation}
Here, $\bm{U}_h$ is determined by the framework of the selected CFD solver. For example, under the finite volume framework, $\bm{U}_h$ is defined by
\begin{equation*}
 \bm{U}_h = \sum_{e\leq N_e} \overline{\bm{U}}_e\cdot \chi_{\Omega_e}(\bm{x}),
\end{equation*}
where $\overline{\bm{U}}_e$ is the cell-average of the element $\Omega_e$, $\bm{U}=\left\{\overline{\bm{U}}_e\right\}_e$ is the solution vector, $\chi_{\Omega_e}(\bm{x})$ is the characteristic function of set $\Omega_e$ and related with $\bm{X}$.

After presenting the mathematical description of an aerodynamic optimization problem (\ref{eq:opt2}) and clarifying the dependency relations between variables (\ref{eq:XD}-\ref{eq:flow}), the entire workflow for the optimization problem is detailed in the next section.

\section{Workflow of adjoint gradient based ASO}

\begin{figure}[htbp]
  \centering
  \includegraphics[width=9cm]{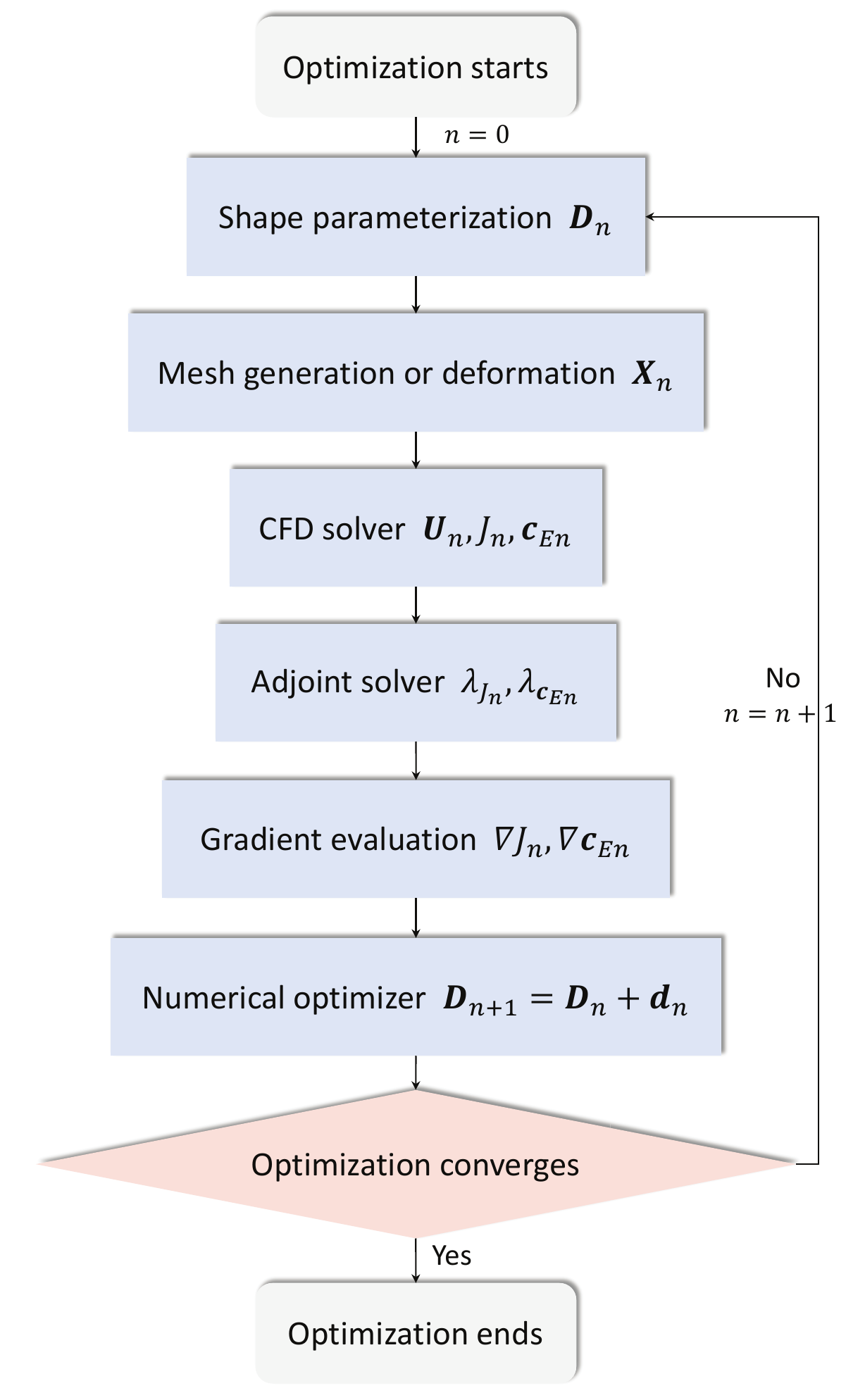}
  \caption{The workflow for the adjoint-enabled gradient-based ASO.}
  \label{fig:opt}
\end{figure}

To solve the problem (\ref{eq:opt2}), a popular choice is the CFD-based ASO through an adjoint-enabled gradient-based optimization algorithm. The entire workflow for the gradient-based optimization is presented in Fig. \ref{fig:opt}, including modules of the geometry parameterization, the mesh deformation, the CFD solver, the adjoint-based gradient evaluation, and the numerical optimizer.

There are various effective open-source frameworks (such as SU2 \cite{economon2016su2}, MACH-Aero \cite{martins2022aerodynamic}, etc.\cite{poirier2016efficient, he2020dafoam}) based on the above workflow. In this work, the discrete adjoint gradient-based ASO algorithm is developed and implemented based on our recently open sourced platform HODG \cite{he2023hodg}, which is a component-based framework of high-order discontinuous Galerkin methods (DGMs) for solving fluid governing equations.

The implementation of each module in Fig. \ref{fig:opt} is detailed as follows.
\subsection{Geometry parameterization}
\qquad The geometry parameterization (GP) module is responsible for mathematically parametric modeling and morph of the aerodynamic shape by $\bm{D}$. Hicks-Henne (HH) \cite{hicks1978wing} and free-form deformation (FFD) \cite{sederberg1986free, hsu1992direct} parameterizations are described and implemented in HODG.
\subsubsection{HH parameterization}
\qquad HH parameterization is a perturbation-based function-fitting method. The upper and lower surface curves of a 2D airfoil are parameterized as follow,
\begin{equation}
\begin{aligned}
 & y(x) = y_{\text{base}}(x)+\sum_{i=1}^N d_i b_i(x),\\
 & b_i(x) = \sin^t\left(\pi x^{\beta_i}\right),\quad i=1,2,..,N.
\end{aligned}
\label{eq:hh}
\end{equation}
Here, $y_{\text{base}}(x)$ is the baseline shape, $d_i$ are the design variables, $b_i(x)$ are basis functions. $N$ denotes the number of sample points on upper and lower surfaces, $t$ is a parameter controlling the compact range of the basis function, $\beta_i=\frac{\log 0.5}{\log p_i}$ represents the horizontal coordinate of the highest point of the basis function $b_i(x)$, that is, the chord position of the peak value of the shape function. During the implementation, $t$ is set to $3$ and $p_i$ satisfies
$$
\begin{aligned}
 &p_i = 0.5(1-\cos \theta_i),\\
 &\theta_i = \frac{i}{N+1}\pi, \quad i=1,..,N.
\end{aligned}
$$

HH parameterization does not need to fit the baseline shape, $d_i$ are set to $0$ at the beginning, and the shape is morphed by
\begin{equation}
 y(x)+\Delta y(x) = y_{\text{base}}(x) + \sum_{i=1}^N \Delta d_i b_i(x)
 \label{eq:hh2}
\end{equation}
\subsubsection{FFD parameterization}
\qquad FFD directly parameterized the locations of the body control points and does not need abstract the geometry. The idea is to embed the geometry body together with its surface control points into a box of flexible plastic which is known as the FFD box, when the FFD box is deformed, the surface control points are consistently shifted.

In detail, the global/physical coordinates of a surface control point $\bm{X}_{surf}=(x,y,z)\in\mathbb{R}^3$ are presented by
\begin{equation}
 \bm{X}_{surf}(u,v,w)=\sum_{i=0}^l \sum_{j=0}^m \sum_{k=0}^n B_i^l(u) B_j^m(v) B_k^n(w) \bm{P}_{i,j,k}.
\label{eq:ffd}
\end{equation}
Here, $(u,v,w)\in [0,1]^3$ are the local/parametric coordinates of the surface control point inside the control box, $\bm{P}_{i,j,k}$ are the design variables which denote the physical coordinates of the FFD box control points, $B_i^l(u),B_j^m(v),B_k^n(w)\in [0,1]$ are Bernstein polynomials of degree $l,m,n$ respectively, the form is expressed by
$$
 B_i^l(u) = \frac{l!}{i! (l-i)!}u^i(1-u)^{l-i}.
$$

In practice, the local parametric coordinates $(u,v,w)$ of each surface control point $\bm{X}_{surf}$ need to be calculated at the beginning, and remain unchanged during the morph of geometry. The movement of surface control points $\Delta\bm{X}_{surf}$ is presented by the change of the design variables $\Delta \bm{P}_{i,j,k}$,
\begin{equation}
 \bm{X}_{surf}(u,v,w)+\Delta\bm{X}_{surf}=\sum_{i=0}^l \sum_{j=0}^m \sum_{k=0}^n B_i^l(u) B_j^m(v) B_k^n(w) \left(\bm{P}_{i,j,k}+\Delta \bm{P}_{i,j,k}\right).
 \label{eq:ffd2}
\end{equation}

\subsection{Mesh Deformation}
\qquad The mesh deformation (MD) module is responsible for providing high-quality grid $\bm{X}$ for the CFD solver each time as provided geometry $\bm{D}$. The idea is to deform the inner volume grid $\bm{X}_{vol}$ according to the movement of the surface grid $\bm{X}_{surf}$ controlled by the design variable $\bm{D}$. The mesh deformation method based on explicit radial basis function (RBF) interpolation \cite{de2007mesh} is described and used in HODG.

In RBF mesh deformation method, relation between $\Delta \bm{X}_{vol}$ and $\Delta \bm{X}_{surf}$ is interpolated by the RBF, the nodes of inner volume grid that are closer to the boundary surface nodes will have a relatively larger displacement.

Denoting $\Delta \bm{X}_{vol}=\big(\Delta \bm{x}_v \big)$ as the displacement of $\bm{X}_{vol}$, where $\Delta \bm{x}_v$ is the displacement vector of an inner volume node $\bm{x}_v$, with the form of $\Delta \bm{x}_v=\big(\Delta x_{v}(\bm{x}_v),\Delta y_{v}(\bm{x}_v),\Delta z_{v}(\bm{x}_v)\big)$, $\Delta x_{v},\Delta y_{v},\Delta z_{v}$ are movement of $\bm{x}_v$ in the $x,y,z$ direction, respectively,
\begin{equation}
\left\{
\begin{aligned}
 & \Delta x_{v}(\bm{x}_v)=\sum_{b=1}^{N_b}\alpha_{x,b} \phi\left( \| \bm{x}_v-\bm{x}_b \| \right),\\
 & \Delta y_{v}(\bm{x}_v)=\sum_{b=1}^{N_b}\alpha_{y,b} \phi\left( \| \bm{x}_v-\bm{x}_b \| \right),\\
 & \Delta z_{v}(\bm{x}_v)=\sum_{b=1}^{N_b}\alpha_{z,b} \phi\left( \| \bm{x}_v-\bm{x}_b \| \right).
\end{aligned}
\right.
\label{eq:RBFs}
\end{equation}
Here, $\bm{x}_b$ denotes one of all $N_b$ boundary nodes, $\| \bm{x}_v-\bm{x}_b \|$ is the distance between the inner node $\bm{x}_v$ and boundary node $\bm{x}_b$, $\bm{\alpha}=(\bm{\alpha}_x,\bm{\alpha}_y,\bm{\alpha}_z)=\left(\alpha_{x,b},\alpha_{y,b},\alpha_{z,b}\right)$ are undetermined weights, $\phi(\|\bm{x}\|)$ is a RBF with the form of
\begin{equation}
 \phi\left(\|\bm{x}\|\right) = \left\{
\begin{aligned}
 &\left( 1-\frac{\|\bm{x}\|}{r} \right)^2  \quad & 0\leq \|\bm{x}\|\leq r \\
 &\quad 0 &\quad \|\bm{x}\|>r
\end{aligned}
\right.
\label{eq:RBF-function}
\end{equation}
where $r>0$ is the compact support radius.

To determine the undetermined weights $\bm{\alpha}$ in (\ref{eq:RBFs}), the following linear system with the size of $dim\times N_b$ needs to be solved
\begin{equation}
\Delta \bm{X}_{surf} = \bm{\Phi}_{b,b}\; \bm{\alpha}.
\label{eq:RBFs-matrix-vector}
\end{equation}

In (\ref{eq:RBFs-matrix-vector}), $\bm{\Phi}_{b,b}$ is the universal basis matrix, whose element in position $(b_i, b_j)$ is expressed as:
\begin{equation*}
 \Phi_{b_i,b_j} = \phi\left( \| \bm{x}_{b_i}-\bm{x}_{b_j} \| \right),
\end{equation*}
where $\bm{x}_{b_i},\bm{x}_{b_j}$ are two boundary nodes.

Based on (\ref{eq:RBFs}), the movement of inner volume gird $\Delta \bm{X}_{vol}$ is determined by
\begin{equation}
\Delta \bm{X}_{vol} = \bm{\Phi}_{v,b}\; \bm{\alpha}.
\label{eq:RBFs-deformation}
\end{equation}
$\bm{\Phi}_{v,b}$ with the size of $N_v\times N_b$ is the universal basis matrix whose element in position $(v_i, b_j)$ is expressed as:
\begin{equation*}
 \Phi_{v_i,b_j} = \phi\left( \| \bm{x}_{v_i}-\bm{x}_{b_j} \| \right),
\end{equation*}
where $\bm{x}_{v_i},\bm{x}_{b_j}$ denote a inner volume node and a boundary node, respectively.

As a result, relation between $\Delta\bm{X}_{vol}$ and $\Delta\bm{X}_{surf}$ can be summarized as follows,
\begin{equation}
\Delta \bm{X}_{vol} = \bm{\Phi}_{v,b}(\bm{X}_{vol},\bm{X}_{surf})\; \bm{\Phi}_{b,b}^{-1}(\bm{X}_{surf})\;\Delta\bm{X}_{surf}.
\label{eq:RBFs-deformation2}
\end{equation}

\subsection{High-order CFD solver}
\qquad The CFD module is responsible for solving the solution vector $\bm{U}$ and represent the flow variable $\bm{U}_h$ based on grid $\bm{X}$. In HODG, a high-order DGM is applied to solve the compressible Euler equations for analyzing inviscid steady flows.

The 3D steady compressible Euler equations can be expressed as
\begin{equation}
\label{eq:3deuler}
\left\{
\begin{aligned}
 &\nabla \cdot \bm{\mathcal{F}}_I(\bm{u})= \bm{0}  &&\bm{x}\in\Omega \\
 &\mathcal{B}\bm{u}(\bm{x})=\bm{0}      &&\bm{x}\in \partial \Omega=:\Gamma
\end{aligned}
\right.
\end{equation}

The conservative variable vector $\bm{u}$ and the inviscid fluxes $\bm{\mathcal{F}}_I(\bm{u})=\left[\bm{F}_{I,1}(\bm{u}),\bm{F}_{I,2}(\bm{u}),\bm{F}_{I,3}(\bm{u})\right]$ are defined by
\begin{equation}
\bm{U}=\begin{pmatrix}
\rho \\
\rho v_1 \\
\rho v_2 \\
\rho v_3 \\
\rho E
\end{pmatrix},\quad
\bm{F}_{I,1}= \begin{pmatrix}  \rho v_1   \\ \rho v_1^2+p \\ \rho v_1v_2\\ \rho v_1v_3  \\ v_1(E+p)  \end{pmatrix},\quad
\bm{F}_{I,2}= \begin{pmatrix}  \rho v_2   \\ \rho v_1v_2 \\ \rho v_2^2+p\\ \rho v_2v_3  \\ v_2(E+p)  \end{pmatrix},\quad
\bm{F}_{I,3}= \begin{pmatrix}  \rho v_3   \\ \rho v_1v_3 \\ \rho v_2v_3\\ \rho v_3^2+p  \\ v_3(E+p)  \end{pmatrix}.
\label{eq:eulerFlux}
\end{equation}
Here, $\rho,\;p$ and $E$ denote the density, pressure and total energy of the fluid, respectively. $\bm{v}=(v_1,v_2,v_3)$ is the velocity vector of the flow. The pressure variable $p$ can be computed from the equation of state
\begin{equation*}
 p=(\gamma-1)(E-\frac{1}{2}\rho \bm{v}^2),
\end{equation*}
where $\gamma$ is the ratio of the specific heats.

DGMs are used to discretize the weak form of the governing equations (\ref{eq:3deuler}) over the domain $\Omega$
 \begin{equation}
 \label{eq:weak}
  \int_{\Gamma}\bm{\mathcal{F}}_I(\bm{u})\cdot \bm{n} \; \phi\;d\Gamma -
  \int_{\Omega}\bm{\mathcal{F}}_I(\bm{u})\cdot \nabla \phi \; d\Omega =\bm{0},
 \end{equation}
where $\Gamma(=\partial \Omega)$ denotes the boundary of $\Omega$ and $\bm{n}$ is the unit outward normal vector to the boundary.

DGMs apply the following piecewise polynomials $\bm{U}_h$ to represent the global flow variable,
\begin{equation}
\label{eq:FE}
\bm{U}_h(\bm{x})=\sum_{e=1}^{N_e}\bm{U}_{h,e}(\bm{x})=\sum_{e=1}^{N_e}\left(\sum_{i=0}^{N_k} \bm{U}_e^{(i)}\cdot \phi_e^i(\bm{x}) \right)\chi_{\Omega_e}(\bm{x}).
\end{equation}
Here, $N_e,N_k$ denote the number of elements and basis functions, respectively, $\phi_e^i(\bm{x})$ are basis functions on element $\Omega_e$. Taylor basis functions \cite{luo2008discontinuous} are used in this work, whose p1 form can be expressed as
\begin{equation}
\label{eq:taylor-basis}
 \phi_e^0(\bm{x})=1, \quad \phi_e^1(\bm{x})=x-x_c,\quad \phi_e^2(\bm{x})=y-y_c, \quad \phi_e^3(\bm{x})=z-z_c,
\end{equation}
where $\Delta x=0.5(x_{\max}-x_{\min}),\;\Delta y=0.5(y_{\max}-y_{\min}),\;\Delta z=0.5(z_{\max}-z_{\min})$, $(x_c,y_c,z_c)$ is the geometric centroid of the element $\Omega_e$, and they are all related with grid $\bm{X}$.

Substituting test function $\phi$ in (\ref{eq:weak}) by elemental basis function $\phi_e^i(\bm{x})$, we can obtain the spatial DG discretization of (\ref{eq:weak}),
\begin{equation}
\label{eq:semi-DG}
 \oint_{\Gamma_e}\widehat{\bm{F}}_{I}(\bm{U}_h^-,\;\bm{U}_h^+,\;\bm{n}_e)\;\phi_e^id\Gamma - \int_{\Omega_e}\bm{\mathcal{F}}_{I}(\bm{U}_h)\cdot \nabla\phi_e^i\;d\Omega=\bm{0}, \qquad e = 1, 2,..,N_e.
\end{equation}
Here, $\widehat{\bm{F}}_I(\bm{U}_h^-,\;\bm{U}_h^+,\;\bm{n}_e)$ is an approximate Riemann solver for the inviscid numerical flux $\bm{\mathcal{F}}_I(\bm{u})\cdot \bm{n}_e$ along $\bm{n}_e$ direction, and $\bm{n}_e=(n_x,n_y,n_z)$ denotes the unit outward normal vector at interfaces $\Gamma_e=\partial \Omega_e$.

In practice, appropriate shock-processing techniques are needed to avoid high-order solution suffering from spurious oscillations. HODG employs a characteristic-compression embedded shock indicator~\cite{feng2021characteristic2} + artificial viscosity (AV)~\cite{2010Adaptive} + bound and positivity-preserving limiter~\cite{cheng2014positivity} to handle shock waves and near-vacuum regions in high-speed flows, greatly enhancing the robustness of DGMs in simulating high-speed flows.

The AV-modified spatial DG discretization is expressed by
\begin{equation}
\label{eq:semi-DG-artificial}
 \oint_{\Gamma_e}\widehat{\bm{F}}_{I}(\bm{U}_h^-,\;\bm{U}_h^+,\;\bm{n}_e)\;\phi_e^id\Gamma - \int_{\Omega_e}\bm{\mathcal{F}}_{I}(\bm{U}_h)\cdot \nabla\phi_e^i\;d\Omega + \epsilon_e\;\int_{\Omega_e}\nabla\bm{U}_h\cdot \nabla\phi_e^i\;d\Omega=\bm{0}.
\end{equation}

$\epsilon_e$ is namely the artificial viscosity coefficient, HODG adopts a Hartmann-type $\epsilon_e$ with strong residual term, whose form is presented below,
\begin{equation*}
\begin{aligned}
\epsilon_e=&C_{\epsilon}h_e^{2-\beta}\frac{\int_{\Omega_e}\left|\nabla \cdot \bm{F}(\bm{u}_h)\right|d\Omega}{|\Omega_e|},\quad or \quad
\epsilon_e=&C_{\epsilon}h_e^{2-\beta}\frac{\oint_{\Gamma_e}\left|\bm{F}(\bm{u}_h)\right|d\Gamma}{|\Gamma_e|},
\end{aligned}
\end{equation*}
where $C_{\epsilon}$ and $\beta$ above are hyper-parameters to be tuned, they are usually set to $0.001\sim0.1$, $0.01\sim0.5$, respectively.

Denoting $\bm{U} = \big[ \bm{U}_e^{(i)} \big]$ in (\ref{eq:FE}) as the global solution vector, and define the discretized residual vector as
\begin{equation}
 \bm{R}(\bm{U}) = \left[ \oint_{\Gamma_e}\widehat{\bm{F}}_{I}(\bm{U}_h^-,\;\bm{U}_h^+,\;\bm{n}_e)\;\phi_e^id\Gamma - \int_{\Omega_e}\bm{\mathcal{F}}_{I}(\bm{U}_h)\cdot \nabla\phi_e^i\;d\Omega + \epsilon_e\;\int_{\Omega_e}\nabla\bm{U}_h\cdot \nabla\phi_e^i\;d\Omega \right].
 \label{eq:residual}
\end{equation}

A pseudo-time term is introduced to accelerate the process of solving $\bm{R}(\bm{U})=\bm{0}$, (\ref{eq:semi-DG-artificial}) is then rewritten as the following matrix-vector form,
\begin{equation}
\label{eq:semi-DG-matrix}
\bm{M}\frac{d\bm{U}}{dt}+\bm{R}(\bm{U})=\bm{0},
\end{equation}

There are two techniques used in HODG to accelerate convergence process of the steady solution in ASO.

\textbf{(a) Solution remapping technique}. Solution vector $\bm{U}$ is initialized through remapping from the convergent solution on the previous shape and grid \cite{wang2021local}, instead of through free stream or copying degrees of freedom (DoFs) of the previous solution.

Denoting $\widetilde{\Omega}_e$ is deformed from $\Omega_e$, and $\bm{c}(\bm{x})=( c_x(\bm{x}), c_y(\bm{x}),c_z(\bm{x}) )$ as the displacement mapping from $\Omega_e$ to $\widetilde{\Omega}_e$ . Thus, the relation of DoFs between initial solution $\widetilde{\bm{U}}_h$ on $\widetilde{\Omega}_e$ and convergent solution $\bm{U}_h$ on $\Omega_e$ is presented by,
\begin{equation}
    \int_{\widetilde{\Omega}_e}  \widetilde{\bm{U}}_h \; \widetilde{\phi}_e^i \; d\widetilde{\Omega} = \int_{\Omega_e}  \bm{U}_h \; \phi_e^i \; d\Omega - \oint_{\partial \Omega_e} \widehat{H}(\bm{U}_h^-,\bm{U}_h^+,\bm{n}_e)\;\phi_e^i \;d\Gamma + \int_{\Omega_e} \left(\bm{c}(\bm{x})\cdot \nabla \phi_e^i\right)\;\bm{U}_h\;d\Omega,
    \label{eq:remapping}
\end{equation}
where
$$\widehat{H}(\bm{U}_h^-,\bm{U}_h^+,\bm{n}_e):=\frac{c_{\bm{n}}}{2}\left(\bm{U}_h^-+\bm{U}_h^+\right) - \frac{|c_{\bm{n}}|}{2}\left(\bm{U}_h^+-\bm{U}_h^-\right),$$
with $c_{\bm{n}}(\bm{x}):=\bm{c}(\bm{x})\cdot \bm{n}_e$.

The bound and positivity-preserving slope limiters \cite{cheng2014positivity} are used to treat the DoFs of the remapping solution, in order to avoid the loss of positivity of density or pressure variables in case of large mesh deformation.

\textbf{(b) Implicit temporal discretization}. Through the $1^{\text{st}}$ order backward Euler implicit time-integration, equations (\ref{eq:semi-DG-matrix}) are discretized as
\begin{equation}
    \label{eq:time}
    \bm{M}\frac{\Delta\bm{U}^n}{\Delta t}+\bm{R}(\bm{U}^{n+1})=\bm{0},
\end{equation}
using $1^{\text{st}}$ order approximation to treat $\bm{R}(\bm{U}^{n+1})=\bm{R}(\bm{U}^{n})+(\frac{\partial \bm{R}}{\partial \bm{U}})^n\Delta \bm{U}^n$, nonlinear system (\ref{eq:time}) is now linearized as the following form,
\begin{equation}
    \label{eq:time2}
    \bm{A}\Delta \bm{U}^n=-\bm{R}(\bm{U}^n),\;\text{  with } \;\bm{A}=\frac{\bm{M}}{\Delta t}\bm{I}+(\frac{\partial \bm{R}}{\partial \bm{U}})^n,
\end{equation}
where $\frac{\partial \bm{R}}{\partial \bm{U}}$ is the Jacobian matrix, and $\Delta \bm{U}^n=\bm{U}^{n+1}-\bm{U}^{n}$ is the solution difference between time level $n$ and $n+1$.

The key to speed up convergence in the implicit scheme is the use of local time step and the growth strategy of CFL, which are presented as follows,
\begin{equation}
\begin{aligned}
&\theta=\max\left\{ \min\left\{  \frac{\|\bm{R}(\bm{U}^{n})\|_2}{\|\bm{R}(\bm{U}^{n+1})\|_2} ,\;\theta_{\max} \right\},\; \theta_{\min} \right\},\\
&\text{CFL}^{n+1}=\theta\cdot\text{CFL}^{n}, \qquad \Delta t^n_e=\text{CFL}^n\frac{h_e}{\lambda_{\max}^n}.
\end{aligned}
\label{eq:cfl}
\end{equation}
Here, $\theta_{\min},\theta_{\max}$ are lower and upper limit of growth factor, respectively, $\theta_{\min}=0.8,\theta_{\max}=2.0$ are usually used; $h_e$ denotes the characteristic length of element $\Omega_e$, $\lambda_{\max}^n$ is maximum wave speed (max absolute value of eigenvalue).

This system of linear equations (\ref{eq:time2}) at each time step is solved by GMRES method with the LU-SGS preconditioner \cite{Hong1998A}.

\subsection{Adjoint-based gradient evaluation}
\qquad This module is responsible for evaluating the objective gradient $\nabla_{\bm{D}} J$. For the objective $J=J\Big(\bm{U}(\bm{X}(\bm{D})), \bm{X}(\bm{D})\Big)$, The derivative of objective with respect to the design variable $\frac{d J}{d \bm{D}}$ is evaluated as follows,
\begin{equation}
\frac{dJ}{d\bm{D}}=\frac{\partial J}{\partial \bm{U}}\frac{d \bm{U}}{d \bm{X}}\frac{d \bm{X}}{d \bm{D}}+\frac{\partial J}{\partial \bm{X}}\frac{d \bm{X}}{d \bm{D}}.
\label{eq:grad0}
\end{equation}

In practice, $\bm{U}(\bm{X})$ is implicitly represented by a discretized equation system $\bm{R}(\bm{U},\bm{X})=\bm{0}$, it becomes intractable complex to navigate to evaluate $\frac{d \bm{U}}{d \bm{X}}$ using difference approximation for high-dimensional design space.

The adjoint-based gradient evaluation is to avoid the calculation of $\frac{d \bm{U}}{d \bm{X}}$ through introducing Lagrange multipliers, the modified objective $\mathcal{L}$ is defined by
\begin{equation}
 \mathcal{L}(\bm{U},\bm{X},\bm{\lambda})=J(\bm{U},\bm{X})-\bm{\lambda}^T \bm{R}(\bm{U},\bm{X}).
\label{eq:lagrange}
\end{equation}
Here, the dimension of Lagrange multiplier $\bm{\lambda}\in \mathbb{R}^{N_e\times N_m\times N_k}$ is the same as that of the solution vector $\bm{U}$ and residual vector $\bm{R}$, where $N_e, N_m, N_k$ denote the number of element, the number of governing equation, the number of basis function, respectively.

For that $\bm{R}(\bm{U},\bm{X})=\bm{0}$, the newly defined objective does not change the original objective $\mathcal{L}(\bm{U},\bm{X},\bm{\lambda})=J(\bm{U},\bm{X})$, and its derivative satisfies
\begin{equation}
 \frac{d\mathcal{L}}{d\bm{D}}=\left(\frac{\partial J}{\partial \bm{U}}-\bm{\lambda}^T\frac{\partial \bm{R}}{\partial \bm{U}}\right)\frac{d \bm{U}}{d \bm{X}}\frac{d \bm{X}}{d \bm{D}}+\left(\frac{\partial J}{\partial \bm{X}}-\bm{\lambda}^T\frac{\partial \bm{R}}{\partial \bm{X}}\right)\frac{d \bm{X}}{d \bm{D}}.
\label{eq:adj-grad}
\end{equation}

Select Lagrange multiplier $\bm{\lambda}$ in (\ref{eq:adj-grad}) such that
\begin{equation}
\bm{\lambda}^T\frac{\partial \bm{R}}{\partial \bm{U}}=\frac{\partial J}{\partial \bm{U}}\quad \text{or} \quad (\frac{\partial \bm{R}}{\partial \bm{U}})^T\bm{\lambda}=(\frac{\partial J}{\partial \bm{U}})^T,
\label{eq:adjoint}
\end{equation}
(\ref{eq:adj-grad}) can be simplified into
\begin{equation}
\label{eq:adj-grad2}
\frac{d\mathcal{L}}{d\bm{D}}=\left(\frac{\partial J}{\partial \bm{X}}-\bm{\lambda}^T(\frac{\partial \bm{R}}{\partial \bm{X}})\right)\frac{d \bm{X}}{d \bm{D}}.
\end{equation}

Note that $\frac{d \bm{U}}{d \bm{X}}$ is now eliminated in (\ref{eq:adj-grad}), and the solution $\bm{\lambda}$ to system (\ref{eq:adjoint}) is called the discrete adjoint solution with respect to $J$.

The remaining terms $\frac{\partial J}{\partial \bm{X}},\frac{\partial \bm{R}}{\partial \bm{X}}$ in (\ref{eq:adj-grad2}) are determined by the selected CFD solver, $\frac{d \bm{X}}{d \bm{D}}=\frac{d \bm{X}_{vol}}{d \bm{X}_{surf}}\cdot\frac{d \bm{X}_{surf}}{d \bm{D}}$ is determined by the methodology of mesh deformation and geometry parameterization.

As a result, the adjoint-based derivative is evaluated by (\ref{eq:adj-grad2}), where the Jacobian matrix is approximated by (\ref{eq:grad-fd}), and the Lagrange multiplier (adjoint solution vector) is obtained through solving the discrete adjoint equations (\ref{eq:adjoint}). Evaluation of derivatives of other constraints $\frac{d \bm{c}_E}{d \bm{D}},\;\frac{d \bm{c}_I}{d \bm{D}}$ is nearly the same.

\subsection{Numerical optimizer}
\qquad The optimizer module is responsible for updating the design variable $\bm{D}$ according to the objective $J$, constraints $\bm{c}_E,\bm{c}_I$, and their derivatives $\frac{dJ}{d\bm{D}},\;\frac{d \bm{c}_E}{d \bm{D}},\;\frac{d \bm{c}_I}{d \bm{D}}$. The sequential least squares programming (SLSQP) algorithm \cite{johansen2004constrained} is adopted in HODG.

To solve the nonlinear KKT system of problem (\ref{eq:opt2}) by Newton iteration, SLSQP obtains the update direction at each Newton iteration step through transforming the linear system about the update direction into a more stable strictly convex quadratic programming sub-problem.

For an equality constrained optimization problem,
\begin{equation}
\begin{aligned}
\label{eq:opt3}
 & \min_{\bm{D}\in \mathcal{D}} \;J(\bm{D}) \\
 s.t. \quad & \bm{c}_E(\bm{D})=\bm{0},
\end{aligned}
\end{equation}
each SLSQP iteration step includes the following 3 sub-steps,
\begin{itemize}
 \item KKT system Assembly. The Lagrange function of (\ref{eq:opt3}) is defined by
  \begin{equation} \mathcal{L}(\bm{D},\bm{\mu})=J(\bm{D})-\bm{\mu}^T \bm{c}_E(\bm{D}). \label{eq:lagrange-opt} \end{equation}
  The KKT system is presented by
  \begin{equation}
    \left\{
    \begin{aligned}
    & \big(\frac{\partial \mathcal{L}}{\partial \bm{D}}\big)^T=\nabla J(\bm{D})-\nabla \bm{c}_E(\bm{D})\bm{\mu}=\bm{0}, \\
     & \big(\frac{\partial \mathcal{L}}{\partial \bm{\mu}}\big)^T = -\bm{c}_E(\bm{D})=\bm{0}.
    \end{aligned}
    \;\;\Longleftrightarrow \;\; \bm{F}(\bm{D},\bm{\mu})=
    \begin{bmatrix}
    \nabla J(\bm{D})-\nabla \bm{c}_E(\bm{D})\bm{\mu}\\
    \bm{c}_E(\bm{D})
    \end{bmatrix}=\bm{0}.
    \right.
    \label{eq:KKT}
    \end{equation}
 \item Newton iteration of the KKT system. Denoting $(\bm{D}_k,\bm{\mu}_k)$ as the variables at current iteration step, the update direction $(\bm{d}_{\bm{D}},\bm{d}_{\bm{\mu}})^T$ of the Newton iteration is the solution to the following linear system,
    \begin{equation}
    \nabla\bm{F}(\bm{D},\bm{\mu})^T
    \begin{bmatrix}
    \bm{d}_{\bm{D}}\\
    \bm{d}_{\bm{\mu}}
    \end{bmatrix}=-\bm{F}(\bm{D},\bm{\mu})\;\Longleftrightarrow\;
    \begin{bmatrix}
    \nabla^2\mathcal{L}_k  &-\nabla \bm{c}_{E,k} \\
    \nabla \bm{c}_{E,k}^T &\bm{O}
    \end{bmatrix}
    \begin{bmatrix}
    \bm{d}_{\bm{D}}\\
    \bm{\mu}_k+\bm{d}_{\bm{\mu}}
    \end{bmatrix}=
    -\begin{bmatrix}
    \nabla J_k \\
    \bm{c}_{E,k}
    \end{bmatrix}
    \label{eq:Newton}.
    \end{equation}
  \item Quadratic programming sub-problem. Transform linear system (\ref{eq:Newton}) into the following strictly convex quadratic programming sub-problem,
      \begin{equation}
        \begin{aligned}
        \min_{\bm{d}_k} \; &\nabla J_k^{T} \bm{d}_k +\frac{1}{2} \bm{d}_k^{T} \bm{B}_{k} \bm{d}_k, & \\
        s.t. \quad & \bm{c}_{E,k}+\nabla \bm{c}_{E,k}^{T} \bm{d}_k = \bm{0}.
        \end{aligned}
        \label{eq:SLSQP}
      \end{equation}
  Here, $\bm{B}_k$ is a positive definite approximation to $\nabla^2\mathcal{L}(\bm{D}_k,\bm{\mu}_k)$.
\end{itemize}

As a result, the variables at next iteration step are updated by
\begin{equation}
\begin{bmatrix}
\bm{D}_{k+1}\\
\bm{\mu}_{k+1}
\end{bmatrix}=
\begin{bmatrix}
\bm{D}_{k}\\
\bm{\mu}_{k}
\end{bmatrix}
+ \begin{bmatrix}
\bm{d}_{\bm{D}}\\
\bm{d}_{\bm{\mu}}
\end{bmatrix},
\label{eq:update-dv}
\end{equation}
where $\bm{d}_k=(\bm{d}_{\bm{D}},\bm{d}_{\bm{\mu}})$ is computed through (\ref{eq:lagrange-opt})-(\ref{eq:SLSQP}). More details about the SLSQP algorithm for inequality constraint optimization can be found in \cite{johansen2004constrained}.

\section{Sensitivity analysis based on high-order solvers}
\qquad The accuracy of gradient evaluation is essential for a gradient-based optimization algorithm. For high-dimensional design space, accurate gradients might be capable of accelerating convergence and exploring superior optimal solutions. In this section, the effects of different CFD solvers are investigated and analyzed from perspective of evaluation of both the objective and its adjoint-enabled gradient.

According the last section, the adjoint-based gradient of objective is presented by,
\begin{equation}
 \frac{dJ}{d\bm{D}} = \underbrace{\left( \frac{\partial J}{\partial \bm{X}} - \bm{\lambda}^T\frac{\partial \bm{R}}{\partial\bm{X}} \right)}_{\text{CFD-related}} \cdot \underbrace{ \left(\; \frac{d\bm{X}}{d\bm{D}} \;\right) }_{\text{GP,MD-related}}
 \label{eq:gradient}
\end{equation}

As we can observe the gradient (\ref{eq:gradient}) is the product of two terms, the former is related to the CFD solver, and the latter is related to both geometry parameterization (GP) and mesh deformation (MD).

In theory, the gradient evaluation is independent of the grid variable $\bm{X}$. While in practice, on the one hand, the selection of GP and MD methods directly affects the second term $\frac{d\bm{X}}{d\bm{D}}$ in (\ref{eq:gradient}), and might introduce pseudo-extremums where $\frac{d\bm{X}}{d\bm{D}}=\bm{0}$; On the other hand, the quality of $\bm{X}$ indirectly affects the accuracy of flow variable $\bm{U}$, and further affects the accuracy of the objective and gradients. In \cite{abergo2023aerodynamic}, the ELA and RBF based ASO are compared to show that different MD methods might result in different optimal solutions.

We now discuss the calculation of $\left(\frac{\partial J}{\partial \bm{X}} - \bm{\lambda}^T\frac{\partial \bm{R}}{\partial\bm{X}}\right)$ in (\ref{eq:gradient}) using different CFD solvers.

\subsection{Theoretical derivation}
\qquad Denote the current shape is $\bm{D}$, the current grid is $\bm{X}$, the current solution vector is $\bm{U}$, the current flow variable is $\bm{U}_h$, the dependency among them is presented in (\ref{eq:XD})-(\ref{eq:flow}).

In general, $J$ is a boundary integral of the geometry, $\bm{R}$ is a discrete form of the governing equation that includes both volume and surface integrals, their expressions are presented as follows,
\begin{equation}
\begin{aligned}
 & J(\bm{U},\bm{X})=\oint_{\Gamma_W} \bm{j}\left(\bm{U}_h(\bm{x})\right)\cdot d\bm{\Gamma}(\bm{x}) = \sum_{e}\int_{\Gamma_e(\bm{X})} \bm{j}\left(\bm{U}_h\right)\cdot\bm{n}_e d\Gamma,\\
 & \bm{R}(\bm{U},\bm{X})=\sum_{e}\left( \int_{\Gamma_e}\bm{r}(\bm{U}_h)\cdot d\bm{\Gamma} - \int_{\Omega_e}R(\bm{U}_h)\;d\Omega \right).
\end{aligned}
\label{eq:JandR}
\end{equation}

Considering the objective $J$, its partial derivative with respect to $\bm{X}$ is calculated by
\begin{equation}
 \frac{\partial J}{\partial \bm{X}} = \sum_{e}\left( \bm{j}\big(\bm{U}_{h}(\bm{x})\big)\cdot \bm{n}_e \Big|_{\partial \Gamma_e} \;\frac{d\Gamma_e}{d\bm{X}} + \int_{\Gamma_e} \bm{j}\left(\bm{U}_h\right)\cdot \frac{d\bm{n}_e}{d\bm{X}}\;d\Gamma + \int_{\Gamma_e} \frac{d\bm{j}}{d\bm{U}_h}\frac{\partial \bm{U}_h}{\partial \bm{X}}\cdot \bm{n}_e\;d\Gamma \right).
 \label{eq:grad-uh}
\end{equation}

The first two terms in (\ref{eq:grad-uh}) reflect the change of domain of the integration due to $\bm{X}$, and the last term reflects the change of the integrand due to $\bm{X}$, it can be observed clearer if using a numerical integration to approximate the integral in (\ref{eq:grad-uh}),
\begin{equation}
\left\{
\begin{aligned}
 &J(\bm{U},\bm{X}) = \sum_{e} \sum_{q} \Big( \bm{j}\big(\bm{U}_h(\bm{x}_q)\big)\cdot \bm{n}_e(\bm{x}_q) \; w_q \Big)\\
 &\frac{\partial J}{\partial \bm{X}} = \sum_{e}\sum_{q}\left( \bm{j}\big(\bm{U}_h(\bm{x}_q)\big)\cdot\bm{n}_e \;\frac{dw_q}{d\bm{X}} + \bm{j}\big(\bm{U}_h(\bm{x}_q)\big)\cdot \frac{d\bm{n}_e(\bm{x}_q)}{d\bm{X}} \; w_q + \frac{d\bm{j}}{d\bm{U}_h}\frac{\partial \bm{U}_h(\bm{x}_q)}{\partial \bm{X}}\cdot \bm{n}_e \; w_q\right)
\end{aligned}
\right.
\label{eq:grad-uh2}
\end{equation}
Here, $\bm{x}_q, w_q$ denotes quadrature points and weights, respectively.

For different numerical schemes, such as the FV and DG solutions,
\begin{equation}
 \bm{U}_h^{\text{FV}}\Big|_{\Omega_e} = \overline{\bm{U}}_e \equiv \text{const},\qquad \bm{U}_h^{\text{DG}}\Big|_{\Omega_e} = \sum_{i}\bm{U}_e^{(i)}\cdot \phi_e^i(\bm{x}).
 \label{eq:uhdg}
\end{equation}

The corresponding solution-based derivatives are calculated by
\begin{equation}
 \left(\frac{\partial J}{\partial \bm{X}}\right)_{\text{FV}} = \sum_{e}\left( \bm{j}\big(\bm{U}_{h}^{\text{FV}}\big)\cdot \bm{n}_e \Big|_{\partial \Gamma_e}\; \frac{d\Gamma_e}{d\bm{X}} + \int_{\Gamma_e} \bm{j}\left(\bm{U}_h^{\text{FV}}\right)\cdot \frac{d\bm{n}_e}{d\bm{X}}\;d\Gamma \right),
 \label{eq:grad-fv}
\end{equation}

\begin{equation}
 \left(\frac{\partial J}{\partial \bm{X}}\right)_{\text{DG}} = \sum_{e}\bigg( \underbrace{\bm{j}\big(\bm{U}_{h}^{\text{DG}}\big)\cdot \bm{n}_e\Big|_{\partial \Gamma_e}}_{\text{ho-obj}}\; \frac{d\Gamma_e}{d\bm{X}} + \int_{\Gamma_e} \underbrace{\bm{j}\left(\bm{U}_h^{\text{DG}}\right)}_{\text{ho-obj}}\cdot \frac{d\bm{n}_e}{d\bm{X}}\;d\Gamma + \underbrace{\int_{\Gamma_e} \frac{d\bm{j}}{d\bm{U}_h} \frac{\partial \bm{U}_h^{\text{DG}}}{\partial \bm{X}}\cdot\bm{n}_e\;d\Gamma}_{\text{ho-mod}} \bigg).
 \label{eq:grad-dg}
\end{equation}

The $\frac{\partial \bm{U}_h^{\text{DG}}}{\partial \bm{X}}$ depends on the form of basis functions and implies information of higher-order moments (such as the gradient) of the numerical solution $\bm{U}_h^{\text{DG}}$.

It can be conclude from (\ref{eq:grad-dg}) and the analysis above that the advantage of using high-order DG solution to evaluate $\frac{\partial J}{\partial \bm{X}}$ seems twofold: (a) DGMs are able to directly acquire high-order objective in the first two terms as compared with FVM; (b) DGMs can resolve the high-order modification which disappears even using higher-order FVMs. Similar conclusions can be made as treating $\frac{\partial \bm{R}}{\partial\bm{X}}$.

\subsection{Practical calculation}
\qquad During the calculation of partial derivatives, the representation relations are as follows,
\begin{equation}
 \bm{D} \stackrel{1}{\longrightarrow} \bm{X}_{surf} \stackrel{2}{\longrightarrow} \bm{X}_{vol} \stackrel{3}{\longrightarrow} \bm{U}_h \stackrel{4}{\longrightarrow} J, \bm{R}
 \label{eq:J-D}
\end{equation}

1. $\bm{X}_{surf}$ is a smooth function with respect to $\bm{D}$, regardless of HH (\ref{eq:hh2}) or FFD (\ref{eq:ffd2}).

2. Relation between $\bm{X}_{vol}$ and $\bm{X}_{surf}$ is implicitly and smoothly provided by the RBF method (\ref{eq:RBFs-deformation2}), while the relation is hard to represented explicitly, it is determined by the solution of a linear system.

3. $\bm{U}_{h}$ is the DG representation (\ref{eq:FE}), which is smoothly expressed by the $\bm{U}$ and $\bm{X}$.

4. $J,\bm{R}$ are integrals of $\bm{U}_h$ (\ref{eq:JandR}), which are smoothly defined as well.

For that the relations between $J,\bm{R}$ and $\bm{D}$ are smoothly dependent but implicitly expressed during the calculation of partial derivatives, in order to avoid the huge computational storage and cost brought by automatic differentiation (AD), the second-order difference is used to approximate the remaining three terms in the total derivative (\ref{eq:gradient}).

Define
\begin{equation}
\Delta \bm{D}_i = \Big(0,0,..,\underbrace{\Delta_i}_{\text{i-th}},..,0 \Big),
\label{eq:deltaD}
\end{equation}
where $\Delta_i\sim N(0,\sigma_i^2)$ is a random, sufficiently small quantity.

Thus, the partial derivatives in (\ref{eq:gradient}) are approximated by
\begin{equation}
\left\{
\begin{aligned}
 & \left( \frac{\partial J}{\partial \bm{X}}\frac{d \bm{X}}{d \bm{D}} \right)_i \approx \frac{J\Big(\bm{U}, \bm{X}(\bm{D}+\Delta\bm{D}_i)\Big)-J\Big(\bm{U}, \bm{X}(\bm{D}-\Delta\bm{D}_i)\Big)}{2\Delta \bm{D}_i} \\
 & \left( \frac{\partial \bm{R}}{\partial \bm{X}}\frac{d \bm{X}}{d \bm{D}} \right)_i \approx \frac{\bm{R}\Big(\bm{U}, \bm{X}(\bm{D}+\Delta\bm{D}_i)\Big)-\bm{R}\Big(\bm{U}, \bm{X}(\bm{D}-\Delta\bm{D}_i)\Big)}{2\Delta \bm{D}_i}
\end{aligned}
\right.
\label{eq:grad-fd}
\end{equation}

%Specifically, during the large amount of perturbations of grid $\bm{X}$, almost all of the computation time comes from solving the linear system (\ref{eq:RBFs-matrix-vector}) with dimension of $N_b$, while $\bm{\Phi}_{b,b}^{-1}(\bm{X}_{surf})$ remains unchanged under the current shape $\bm{D}$ and surface nodes $\bm{X}_{surf}$, only the LU decomposition of $\bm{\Phi}_{b,b}(\bm{X}_{surf})$ is therefore calculated and stored.
Note that the calculation of (\ref{eq:grad-fd}) does not require updating the solution vector $\bm{U}$ after the morph of shape $\bm{D}$, only needs to update the grid $\bm{X}$, the quadrature rule and the basis functions inside each element. Therefore, the large amount of perturbations of grid $\bm{X}$ still keep nearly independent of the dimension of the design space.

The core C++ code for the adjoint gradient-based optimization based on the HODG-1.0 framework \cite{he2023hodg} is presented as follows.
\begin{lstlisting}[language={C++}]
/* n:      optimization iteration step,
   Dn:     design variable
   Xsurf:  surface grid point,
   Xvol:   volume grid point
   Xn:     grid variable
   Un0:    initial solution DoFs on Xn
   Uns:    steady solution DoFs on Xn
   Jn/Cn:  current objective/constraints related to Uns
   Zn/Mun: adjoint solution with respect to Jn/Cn
   dJdD:   adjoint-enabled gradient of objective
   dCdD:   adjoint-enabled gradients of constraints  */
// Initialize Dn through FFD parameterization
n = 0;  parameterization->Shape2Para(Dn);
while ( !optimizer->reachCondition() ) {
    // Update surface grid according to design variables: Dn->Delta_Xsurf
    parameterization->Para2Shape(Dn);
    // Update volume grid based on surface grid: Delta_Xsurf->Delta_Xvol
    deformation->run();
    // Provide grid for CFD solver: Xn,Delta_Xsurf, Delta_Xvol->Xn
    deformation->afterRun();
    // Initialize DoFs on current grid: U(n-1)s,Xn->Un0
    solver->initDoFs();
    // Solve steady flows started from Un0 on Xn: Un0,Xn->Uns
    solver->run();
    // Compute objective, constrains and their adjoint solutions: Uns->Jn,Cn,Zn,Mun
    adjointSolver->run();
    // Evaluate adjoint gradient based on grid perturbation: Uns,Zn,Mun->dJdD,dCdD
    gradEval->run();
    // Update design variables through optimizer: Dn,Jn,Cn,dJdD,dCdD->D(n+1)
    optimizer->optimize(Dn,Jn,Cn,dJdD,dCdD);
    n += 1;
}
\end{lstlisting}

\section{Numerical experiments}
\qquad In this section, we test several classical cases based on the newly developed platform to investigate the performance of different CFD solvers (FVMs and DGMs) on evaluation of the adjoint-enabled gradient, and further the entire optimization process. In order to verify the high-order modification term in (\ref{eq:grad-dg}) resolved by DGMs can indeed improve the accuracy of adjoint gradient, two strategies, (a) mesh refinement and (b) order enrichment, are adopted in FVMs to compensate for the gradient errors of the first part in (\ref{eq:grad-fv}) caused by low accuracy of the objective.

Except for the different CFD solvers, the methods and techniques used in the rest of the ASO keep the same. For each test case, the FFD method is selected for shape parameterization, the RBF-based interpolation is used for mesh deformation, the discrete adjoint solver is employed for gradient evaluation, and SLSQP is used as the numerical optimizer. Each test case is performed on a ThinkStation-P620, with the computer dominant frequency of 2.0GHz and a total RAM of 128 Gb.

\subsection{NACA-0012 airfoil drag minimization}
\qquad The first test-case we consider is the NACA-0012 drag minimization, the target is to reduce the drag $C_d$ by modifying the shape $\bm{D}$ without decreasing the lift $C_l$ and the area $A$ of the airfoil. The free stream is set to $\text{Ma}=0.8, \text{AOA}=1.25^o$. The FFD method with 16 design variables are selected to parameterize and morph the airfoil, and in this case the design variables include only the y-coordinates of the FFD control points and are limited to a range of variation of no more than $25\%$.

\begin{figure}[htbp]
 \centering
 \subfigure[mesh-1]{
 \includegraphics[width=4.8cm]{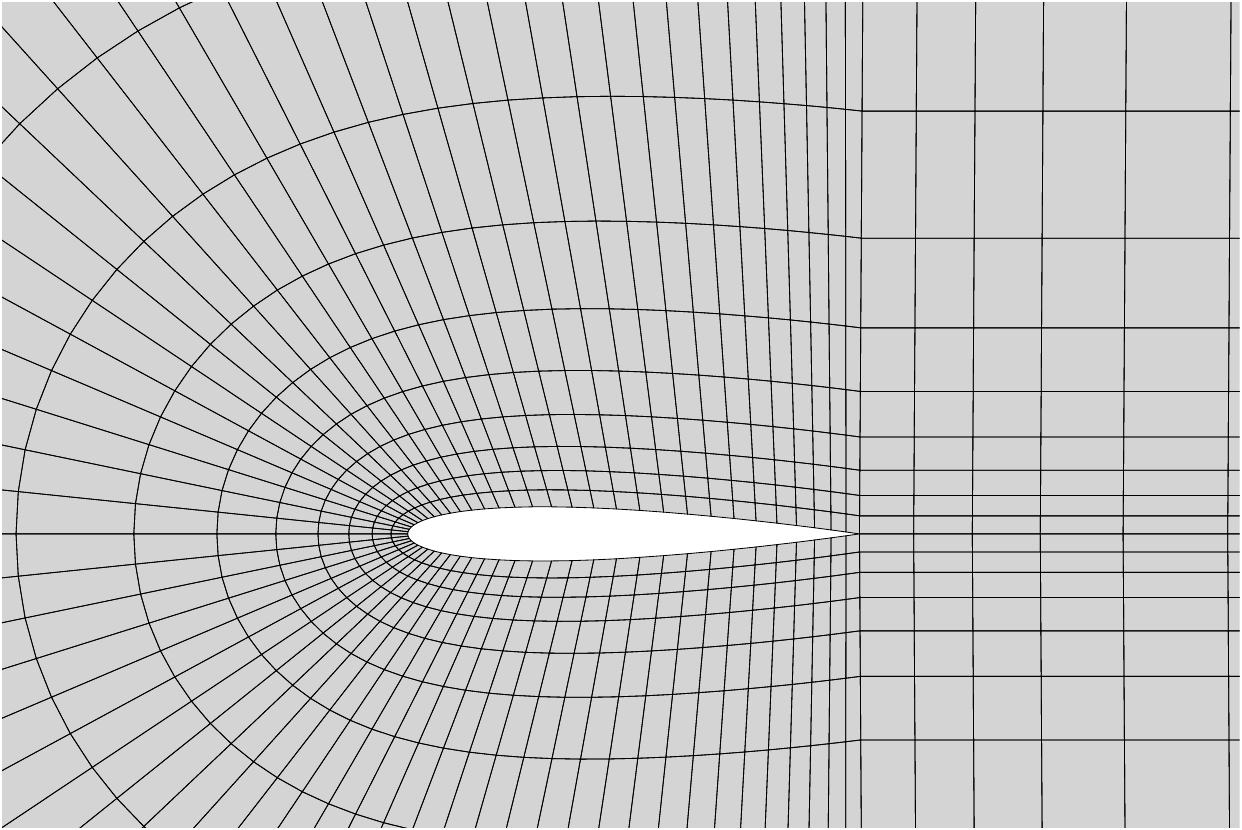}}
 \subfigure[mesh-2]{
 \includegraphics[width=4.8cm]{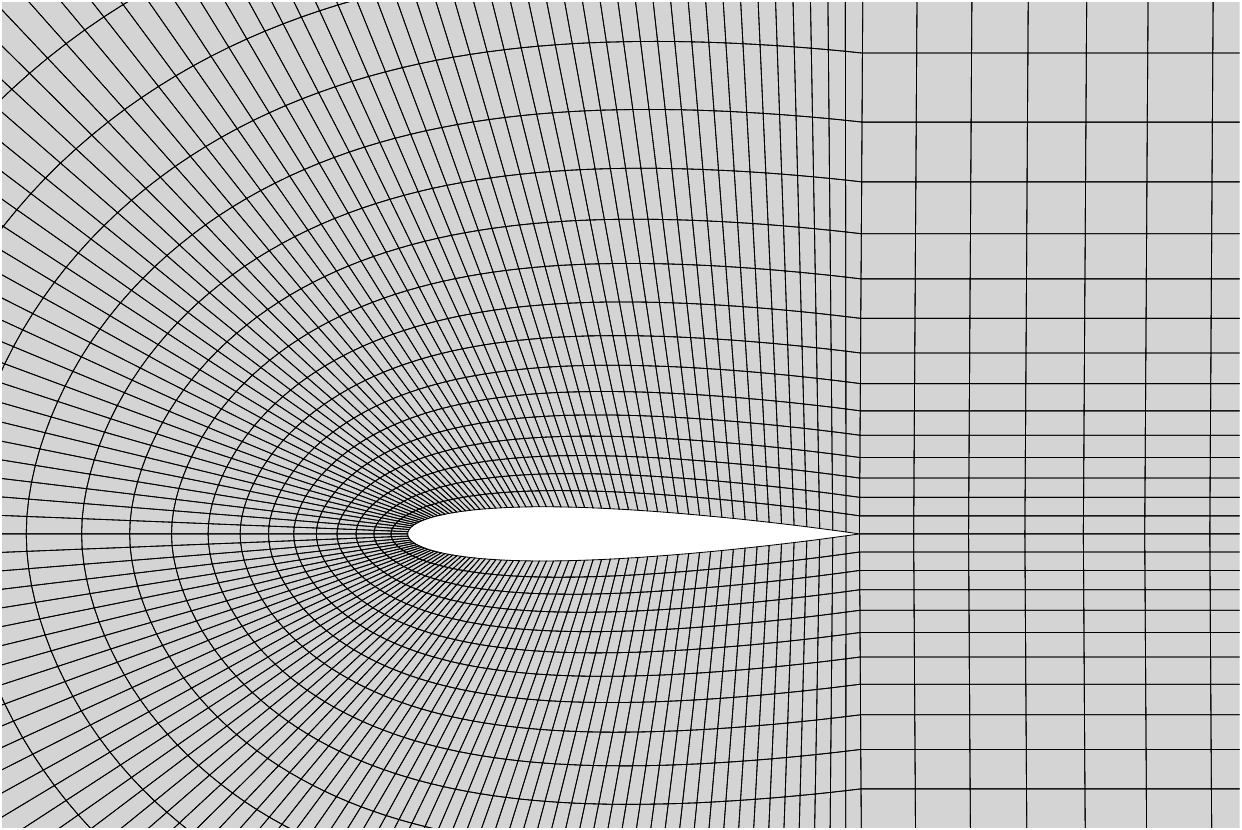}}
 \subfigure[mesh-3]{
 \includegraphics[width=4.8cm]{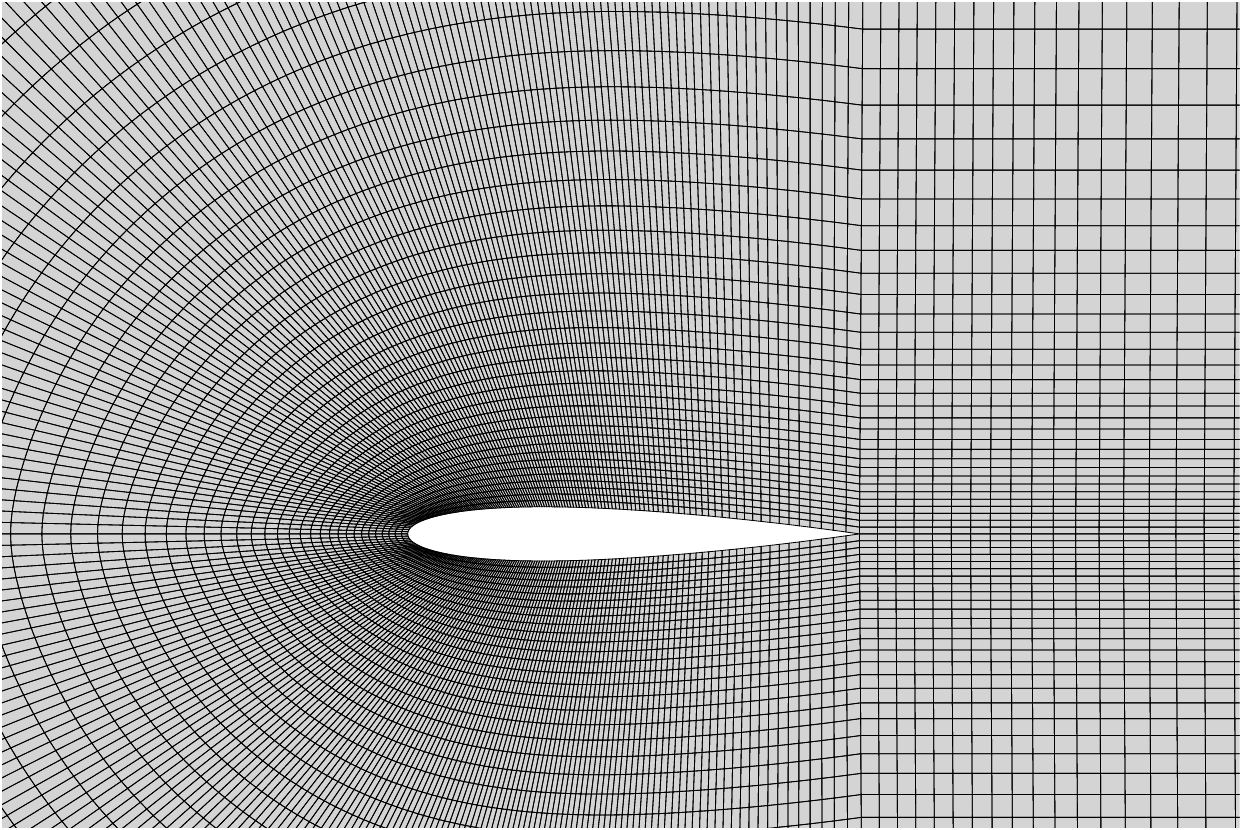}}
 \caption{Computational grids with different levels of the NACA-0012 airfoil test case.}
 \label{fig:naca0012mesh}
\end{figure}

\begin{figure}[htbp]
 \centering
 \subfigure[1stFV on mesh-3]{
 \includegraphics[width=4.8cm]{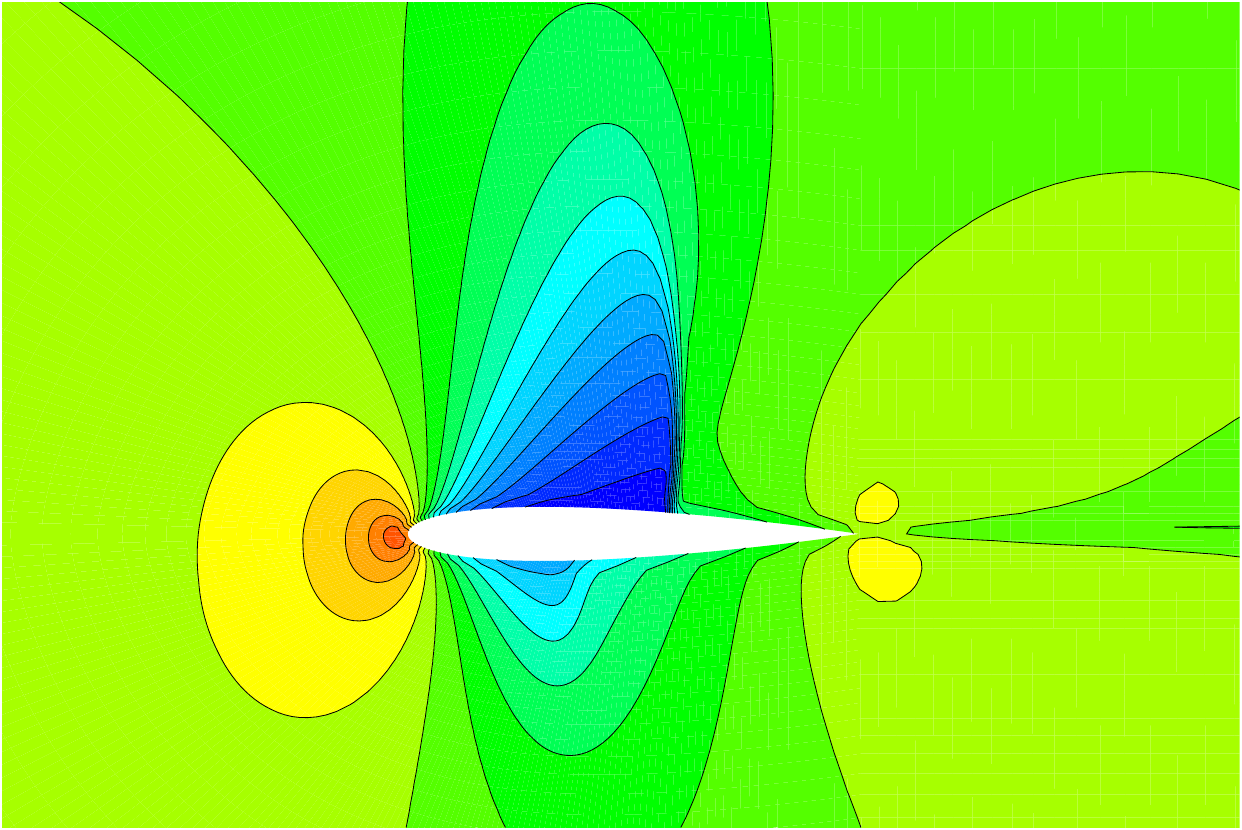}}
 \subfigure[2ndFV on mesh-1]{
 \includegraphics[width=4.8cm]{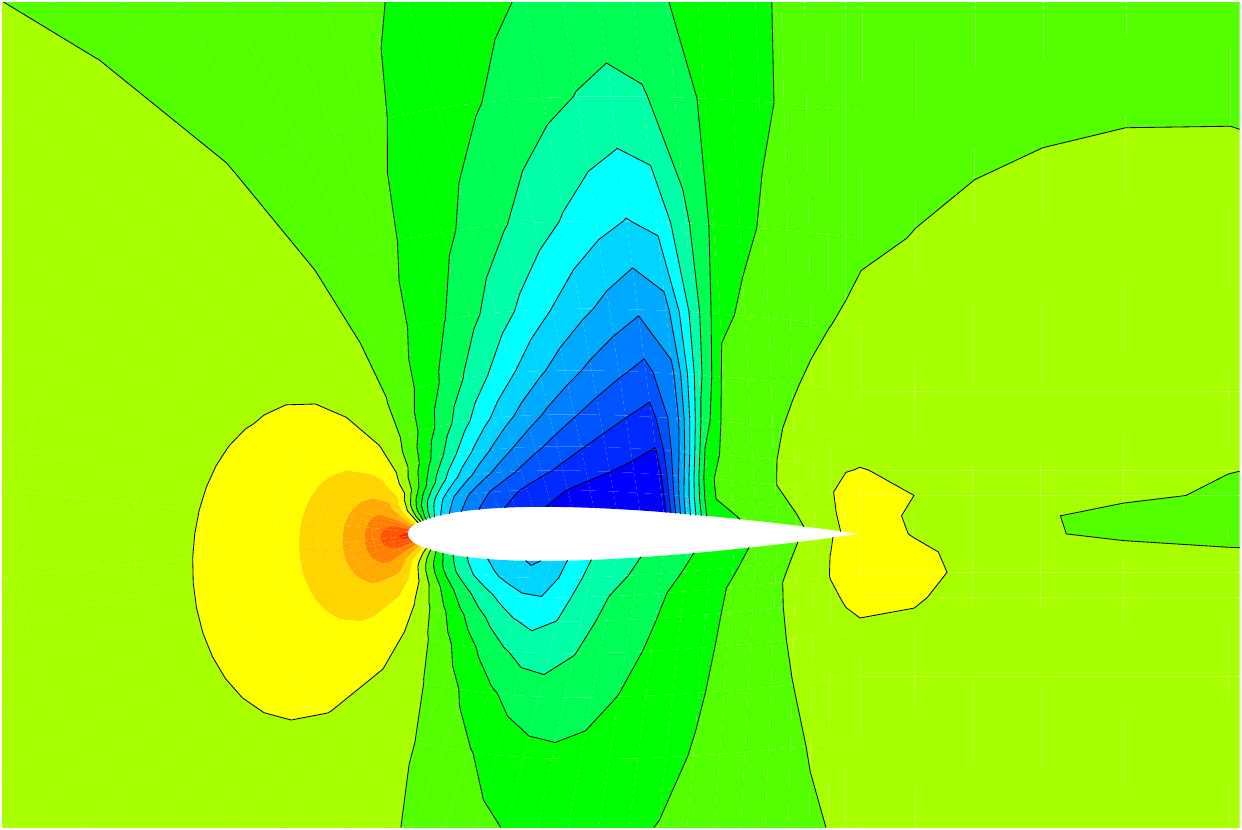}}
 \subfigure[2ndFV on mesh-2]{
 \includegraphics[width=4.8cm]{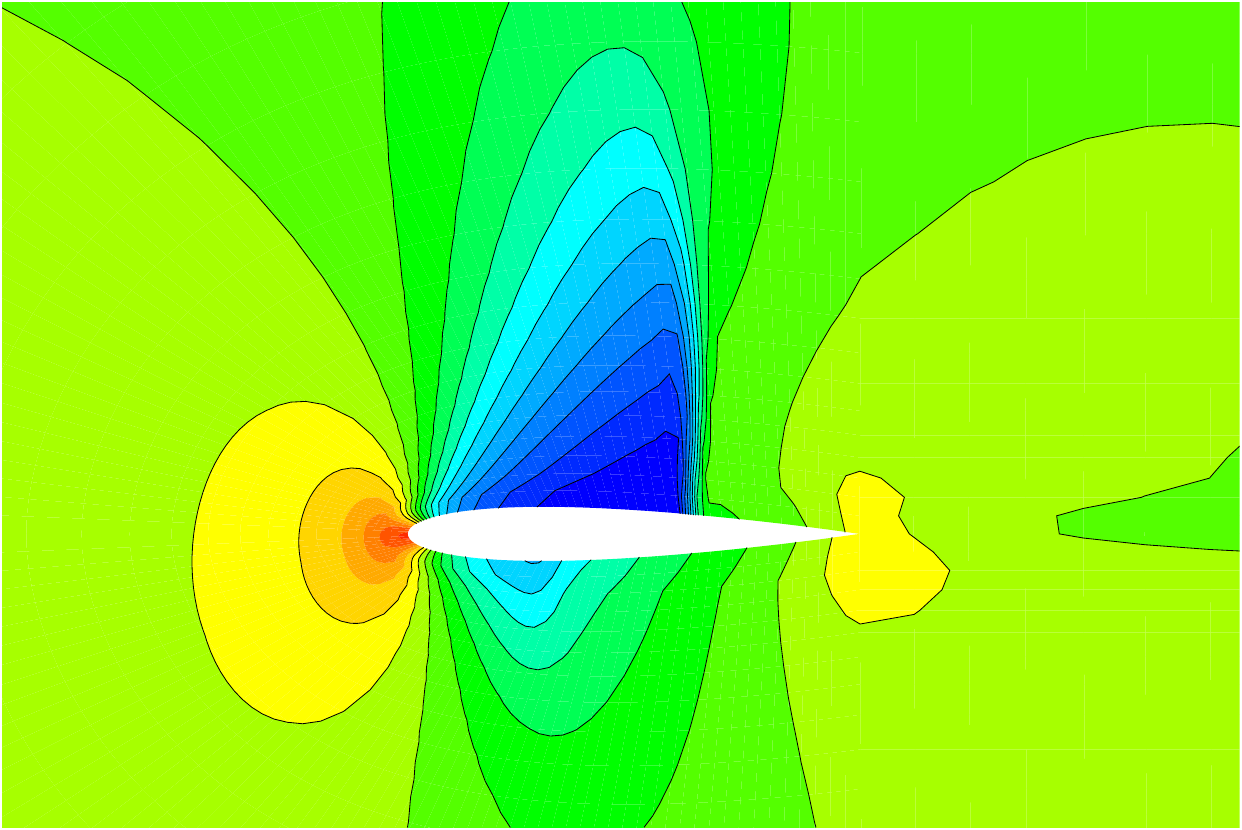}}
 \subfigure[DGp1 on mesh-1]{
 \includegraphics[width=4.8cm]{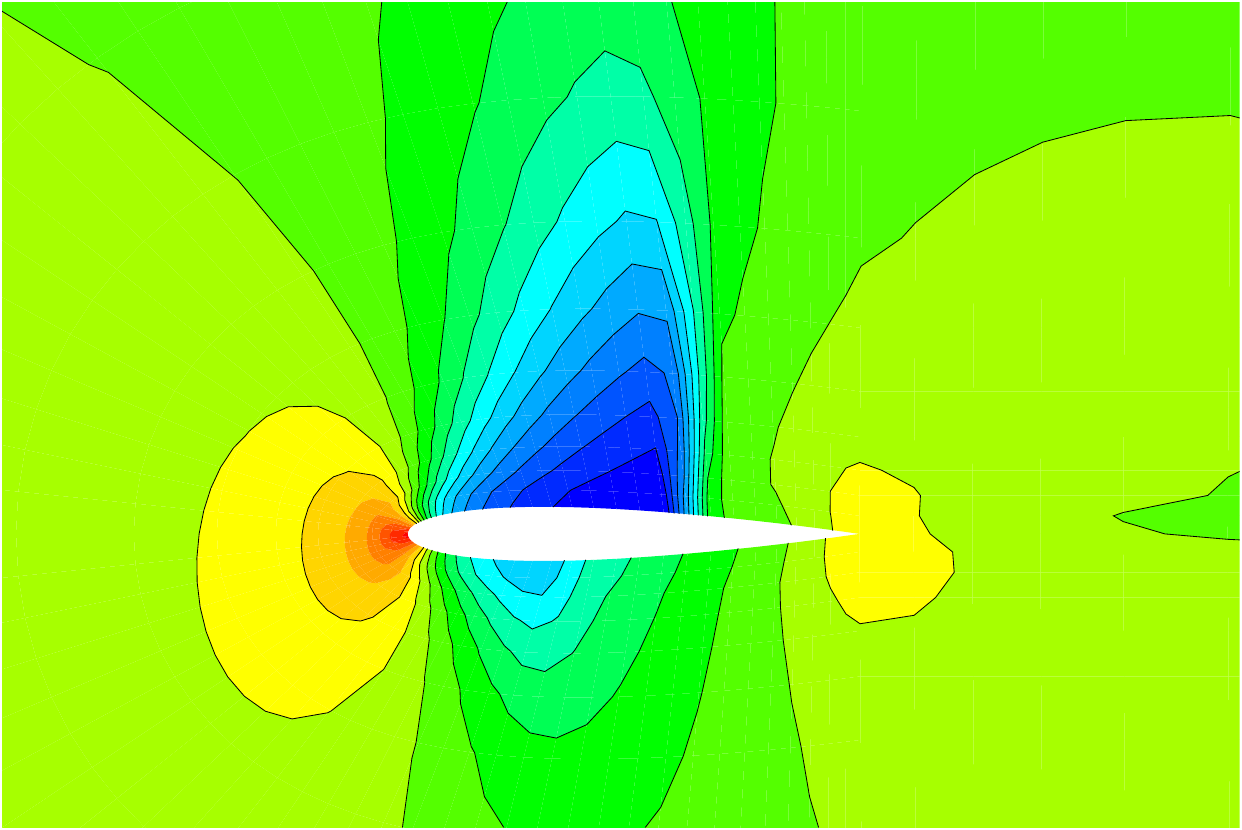}}
 \subfigure[DGp1 on mesh-2]{
 \includegraphics[width=4.8cm]{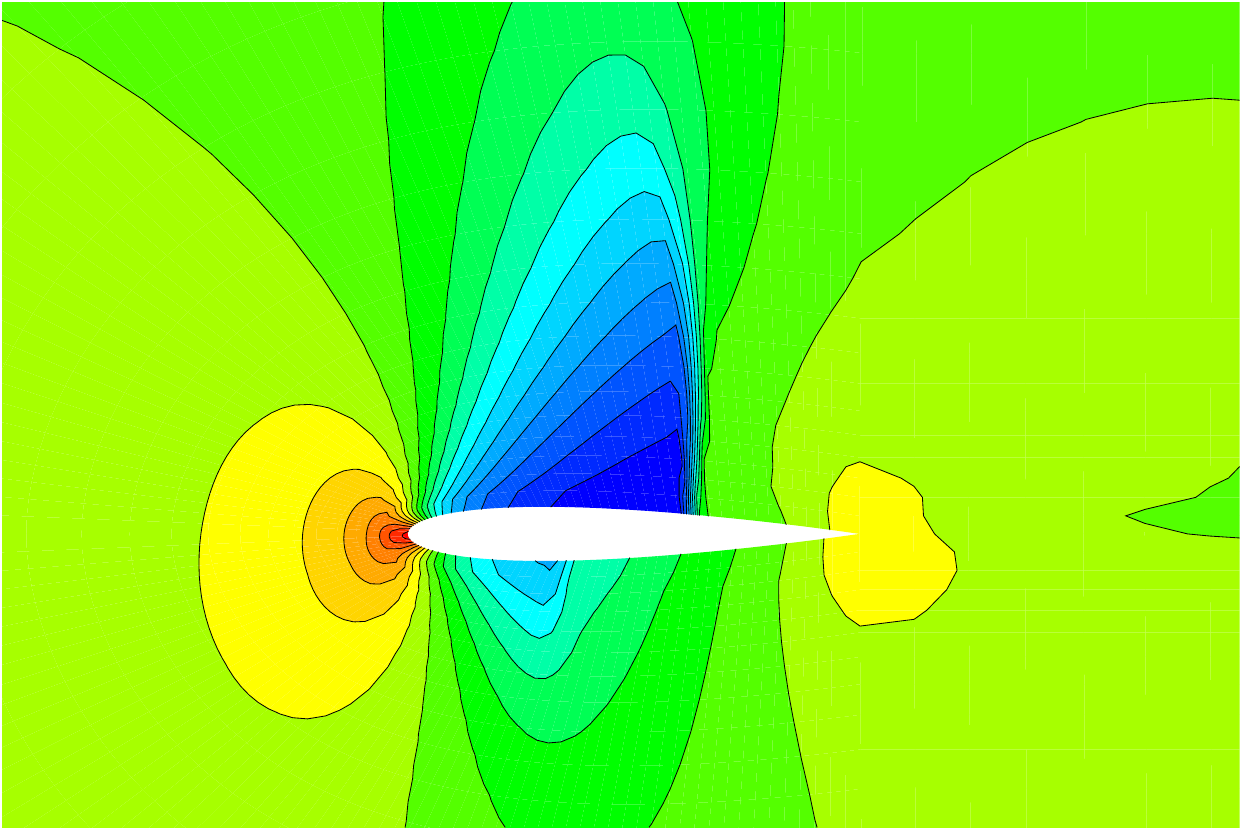}}
 \subfigure[DGp2 on mesh-1]{
 \includegraphics[width=4.8cm]{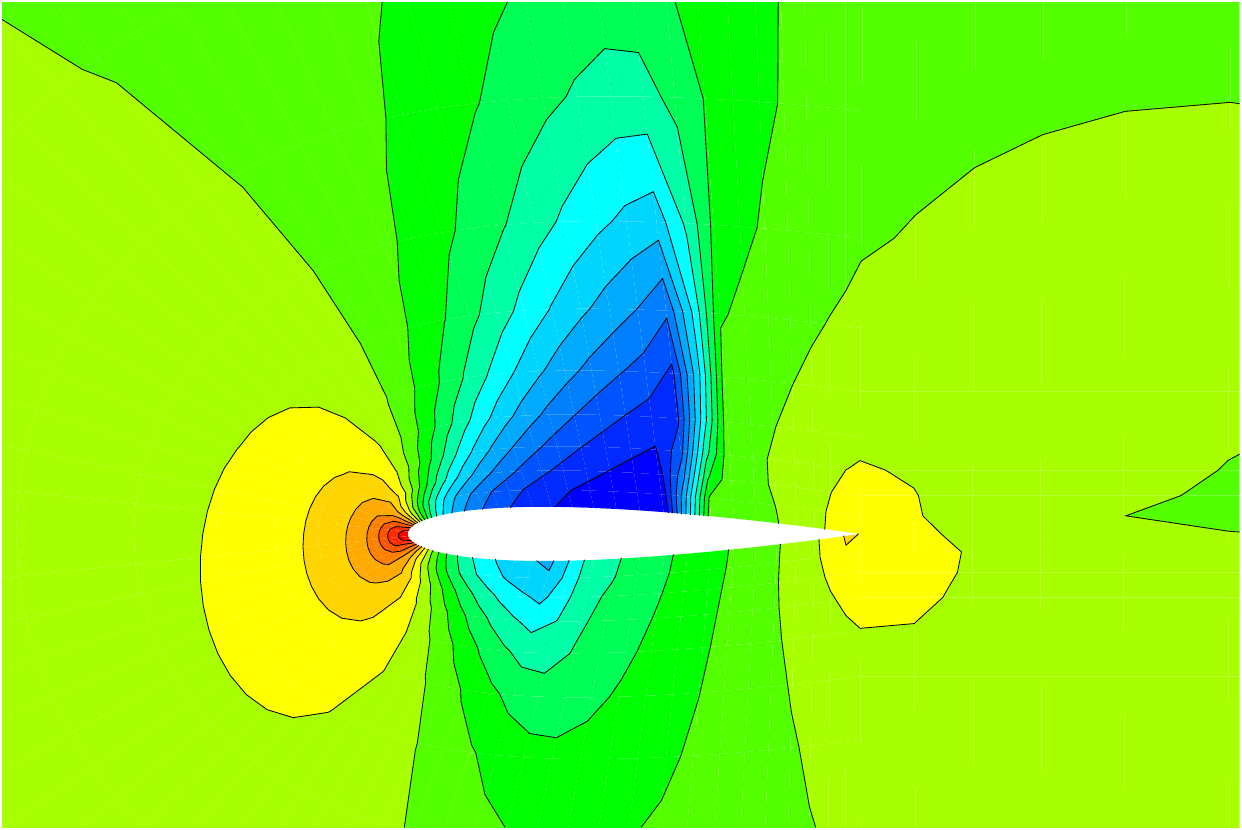}}
 \caption{Solutions (density) obtained by different solvers on various grids of the NACA-0012 airfoil test case.}
 \label{fig:naca0012flow}
\end{figure}

\begin{figure}[htbp]
  \centering
  \includegraphics[width=9.5cm]{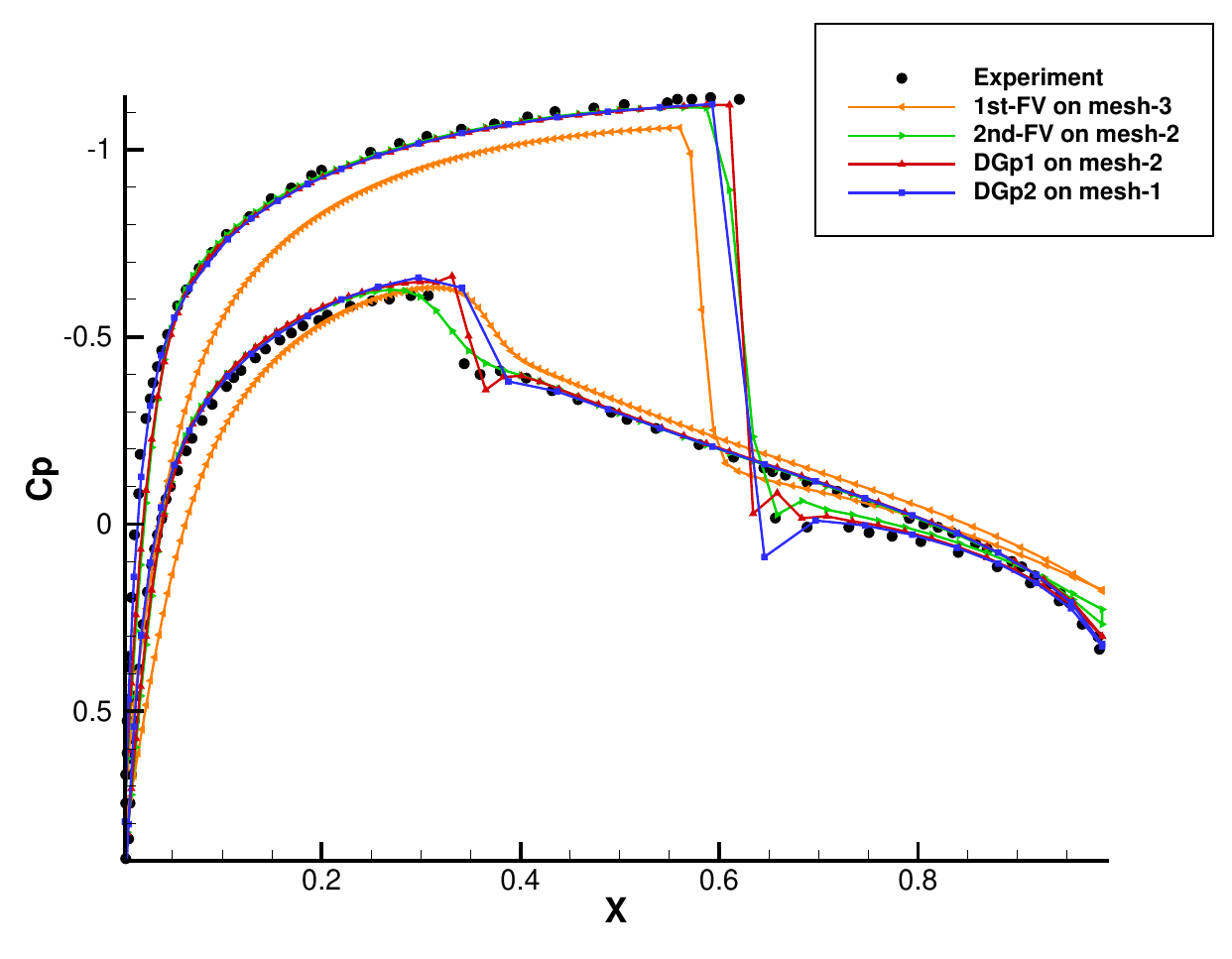}
  \caption{NACA-0012 $C_p$ distributions between CFD data and experimental data.}
  \label{fig:naca0012cp}
\end{figure}

In order to test the performance of different CFD solvers, 3 grids (mesh-1, mesh-2, mesh-3) are used in this test case (as shown in Fig. \ref{fig:naca0012mesh}), the mesh-1 is a C-type grid with 3072 quadrilateral elements, and the mesh-2, mesh-3 come from the global refinement of the mesh-1 once and twice, respectively. The iterative DoFs of different CFD solvers on each grid are presented in Table \ref{tab:naca0012-DoFs}, and the corresponding flow results obtained by different CFD solvers (1stFV, 2ndFV, DGp1, DGp2) are shown in Fig. \ref{fig:naca0012flow}. The $C_p$ distributions obtained by various CFD solvers are compared with experimental data \cite{norstrud1985AGARD} in Fig. \ref{fig:naca0012cp}. It can be observed from Fig. \ref{fig:naca0012flow},\ref{fig:naca0012cp} that similar and acceptable flow results can be reached through DGp2 on mesh-1, DGp1 and 2ndFV on mesh-2, and 1stFV on mesh-3.

\begin{table}[htbp]
\centering \caption{DoFs of various CFD solvers on each grid in the NACA-0012 test case}
\setlength{\tabcolsep}{8mm}
\begin{center}
\begin{tabular}{c|c|c|c}
\toprule
  & mesh-1 & mesh-2  & mesh-3 \\ \hline
  FVMs  &12,288 &49,152 &196,608 \\ \hline
  DGp1 &36,864 &147,456 &/   \\ \hline
  DGp2 &73,728 &/       &/  \\
\bottomrule
\end{tabular}
\end{center}
\label{tab:naca0012-DoFs}
\end{table}

\begin{figure}[htbp]
  \centering
  \includegraphics[width=10.0cm]{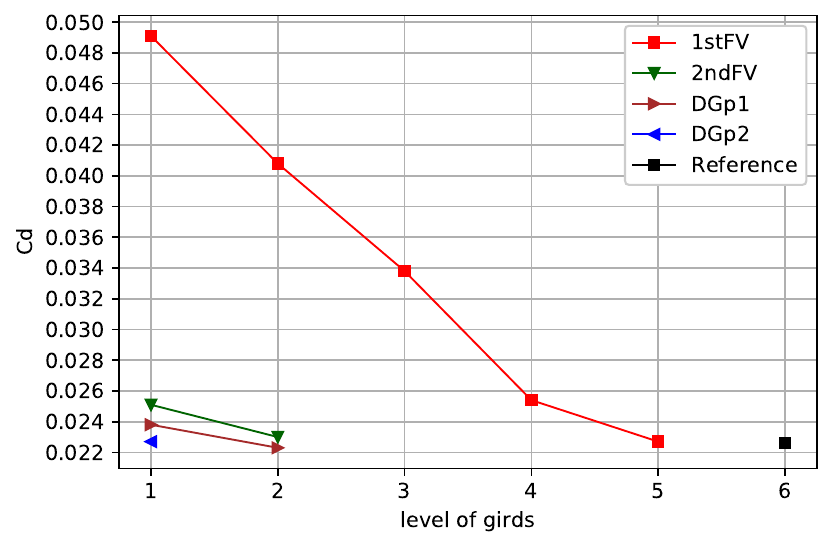}
  \caption{$C_d$ convergence of the NACA-0012 airfoil test case.}
  \label{fig:naca0012-obj}
\end{figure}

Fig. \ref{fig:naca0012-obj} shows the objective $C_d$ convergence of various CFD solvers, and the grid of mesh-4 and mesh-5 in Fig. \ref{fig:naca0012-obj} come from the global refinement of the grid of mesh-3 once and twice, respectively. It can be observed that the objective $C_d$ computed by the FVMs converges with the grid refinement, however, the FVM requires a rather fine gird (mesh-5) to achieve a similar accuracy of $C_d$ as compared with DGp1 on mesh-2 or DGp2 on mesh-1.

\begin{figure}[htbp]
  \centering
  \subfigure[sensitivity convergence of 1st-order FVMs]{
  \includegraphics[width=6.8cm]{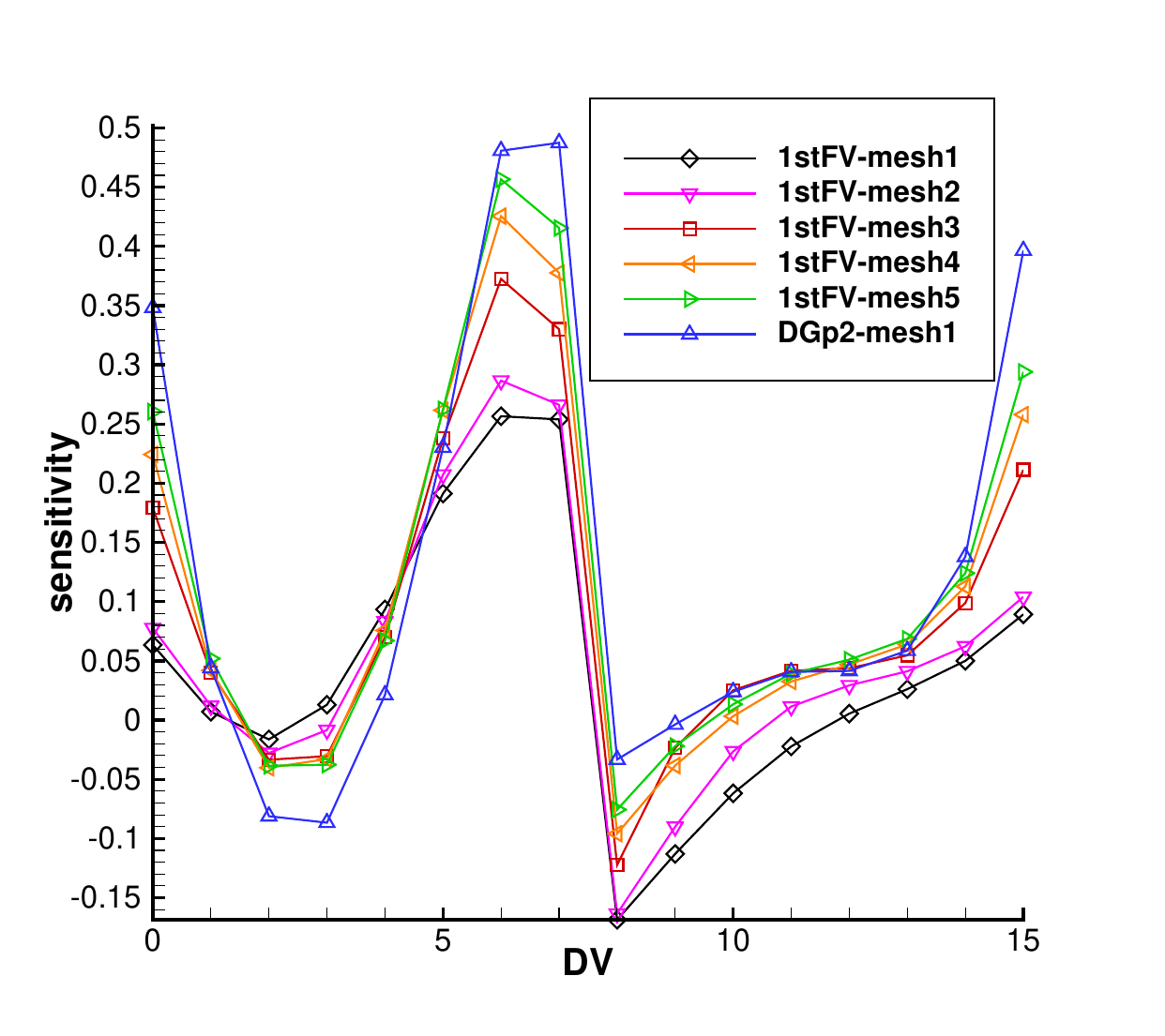}}
  \subfigure[sensitivity comparison under coequal DoFs]{
  \includegraphics[width=6.8cm]{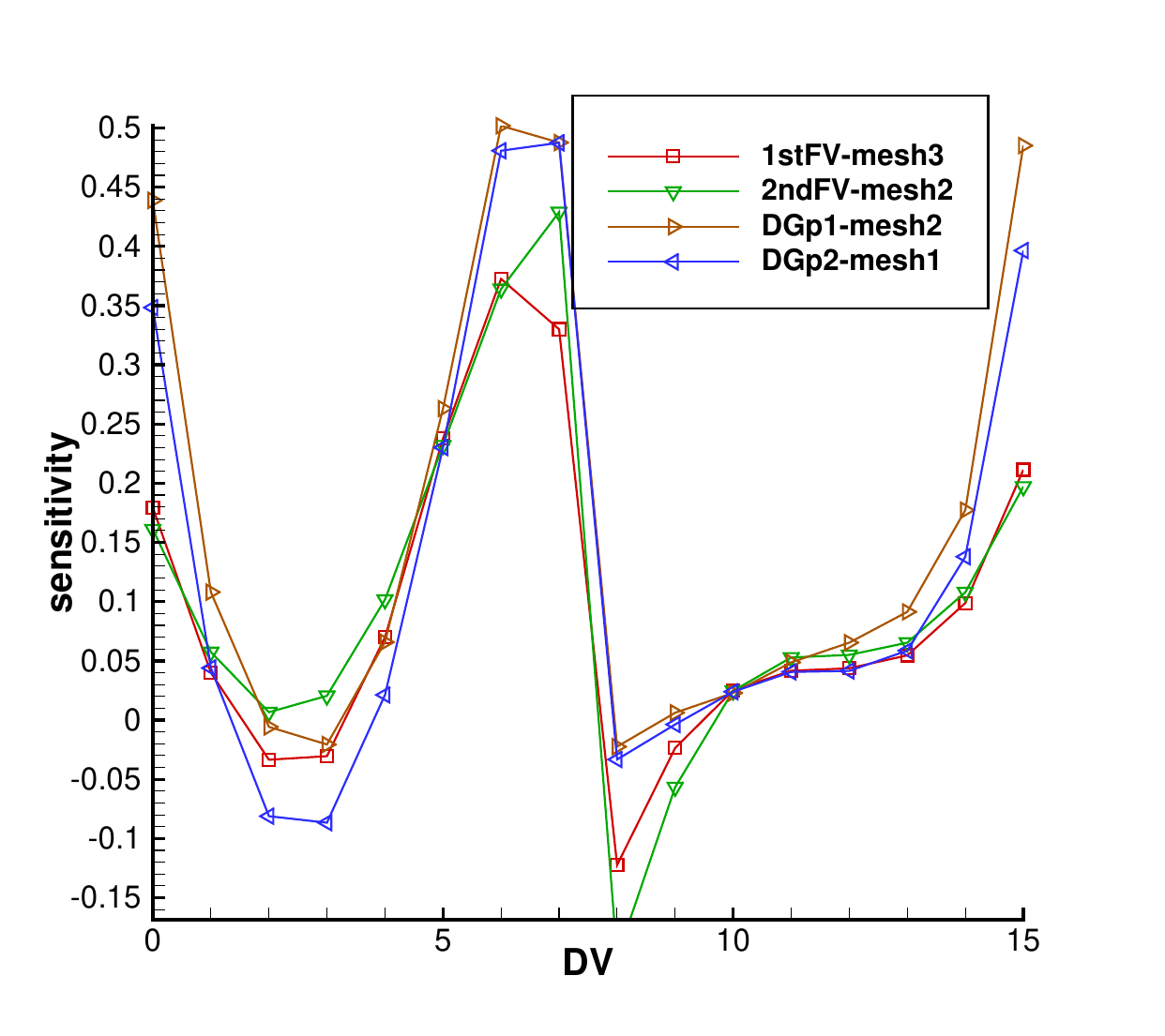}}
  \caption{Drag sensitivity of the NACA-0012 airfoil test case.}
  \label{fig:naca0012grad}
\end{figure}

The results of the corresponding adjoint gradients (sensitivity) computed by different CFD solvers are shown in Fig. \ref{fig:naca0012grad}. It can be observed from Fig. \ref{fig:naca0012grad}(a) that the sensitivity obtained by 1st-order FVMs gradually approaches the sensitivity produced by DGp2 on mesh-1 with the mesh refinement, but still fails to reach the sensitivity produced by DGp2 on mesh-1 at certain design variables even on mesh-5. Fig. \ref{fig:naca0012grad}(b) shows the sensitivity comparison under coequal DoFs. It can be found that the sensitivity of DGp1 and DGp2 are similar, and both are more precise than that of 1stFV on mesh-3 and 2ndFV on mesh-2.

\begin{figure}[htbp]
 \centering
 \subfigure[1stFV on mesh-3]{
 \includegraphics[width=6.0cm]{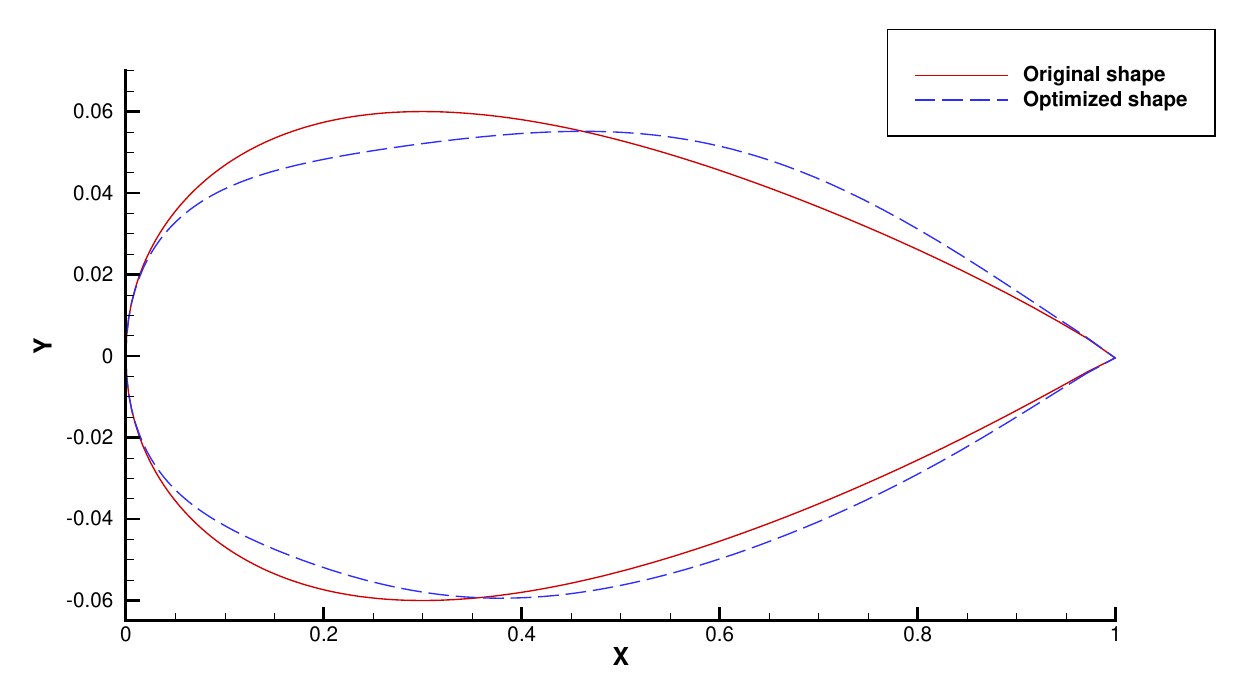}}
 \subfigure[2ndFV on mesh-2]{
 \includegraphics[width=6.0cm]{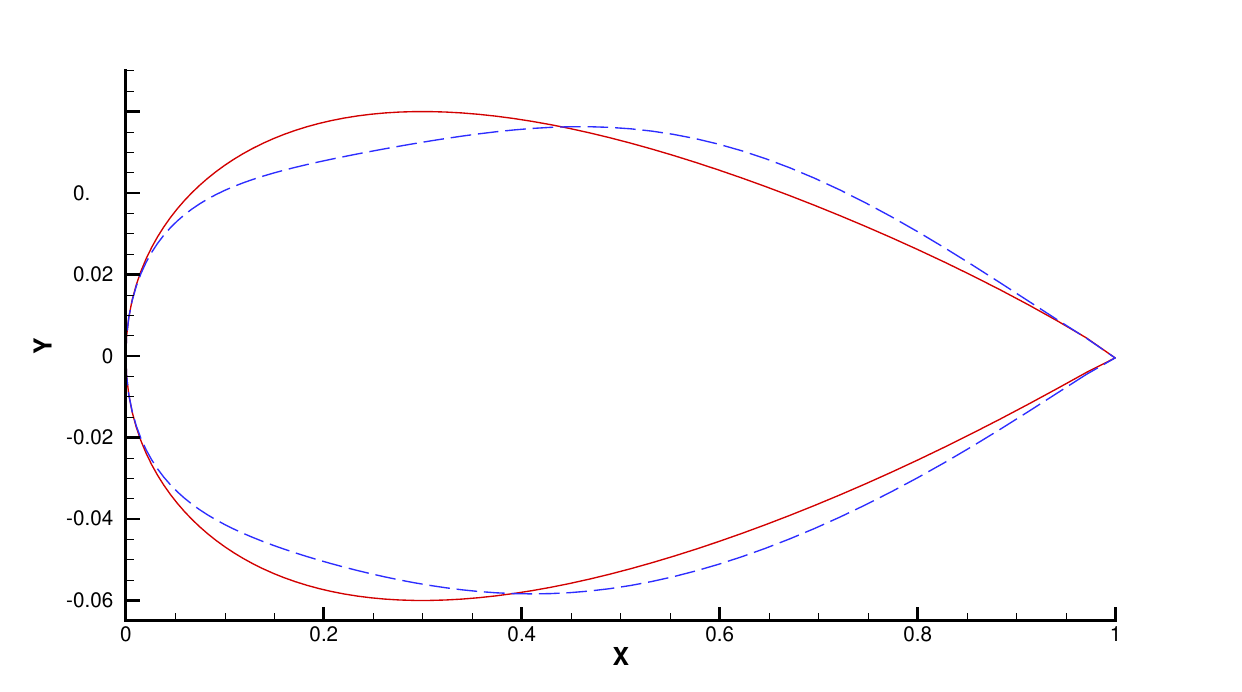}}
 \subfigure[DGp1 on mesh-2]{
 \includegraphics[width=6.0cm]{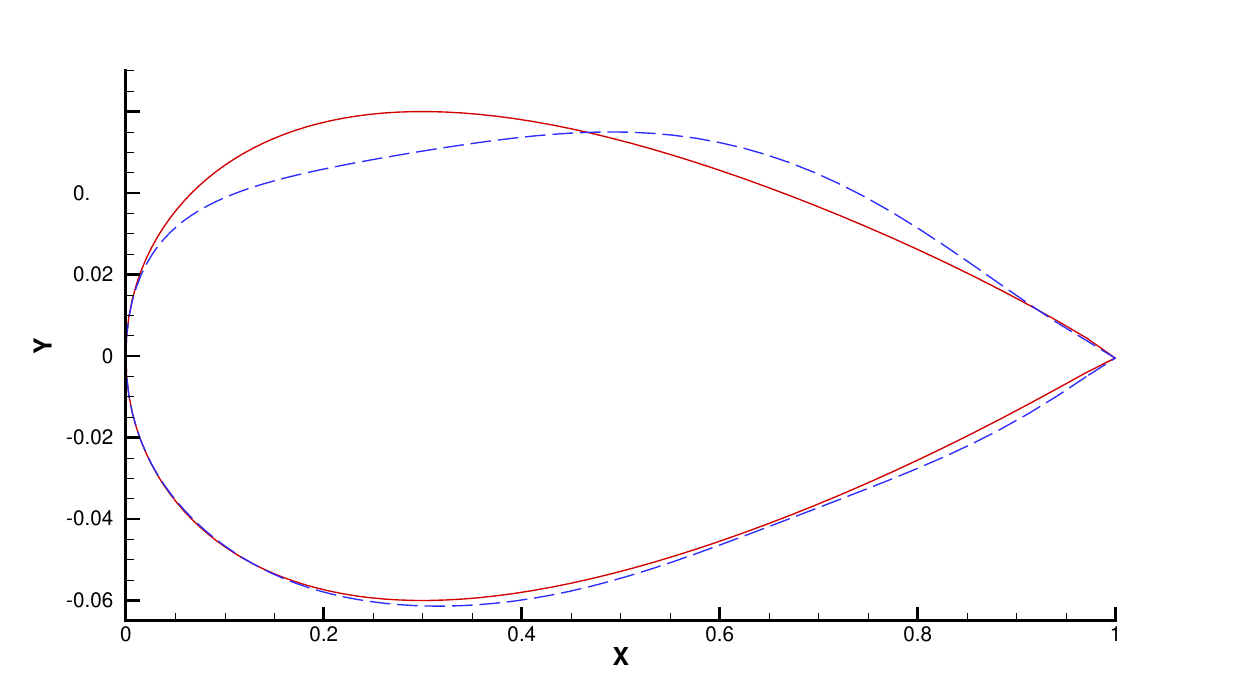}}
 \subfigure[DGp2 on mesh-1]{
 \includegraphics[width=6.0cm]{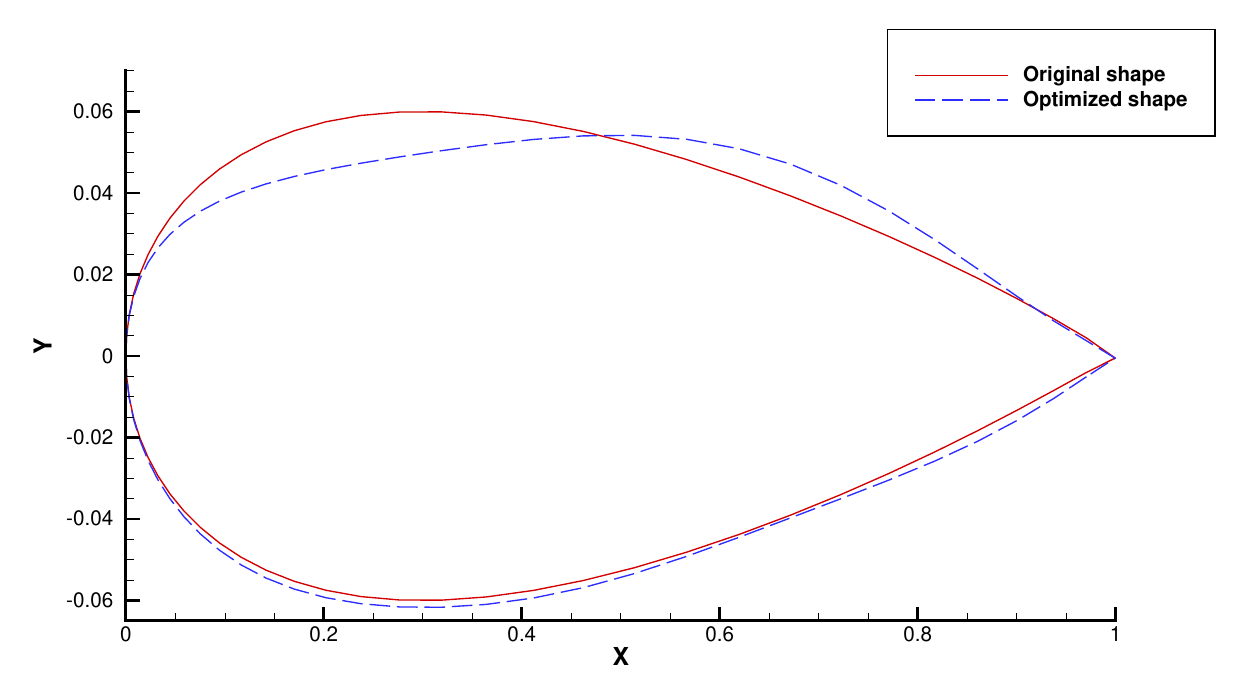}}
 \caption{Original and optimized shape of the NACA-0012 airfoil test case.}
 \label{fig:naca0012shape-opt}
\end{figure}

\begin{figure}[htbp]
 \centering
 \subfigure[1stFV on mesh-3]{
 \includegraphics[width=6.0cm]{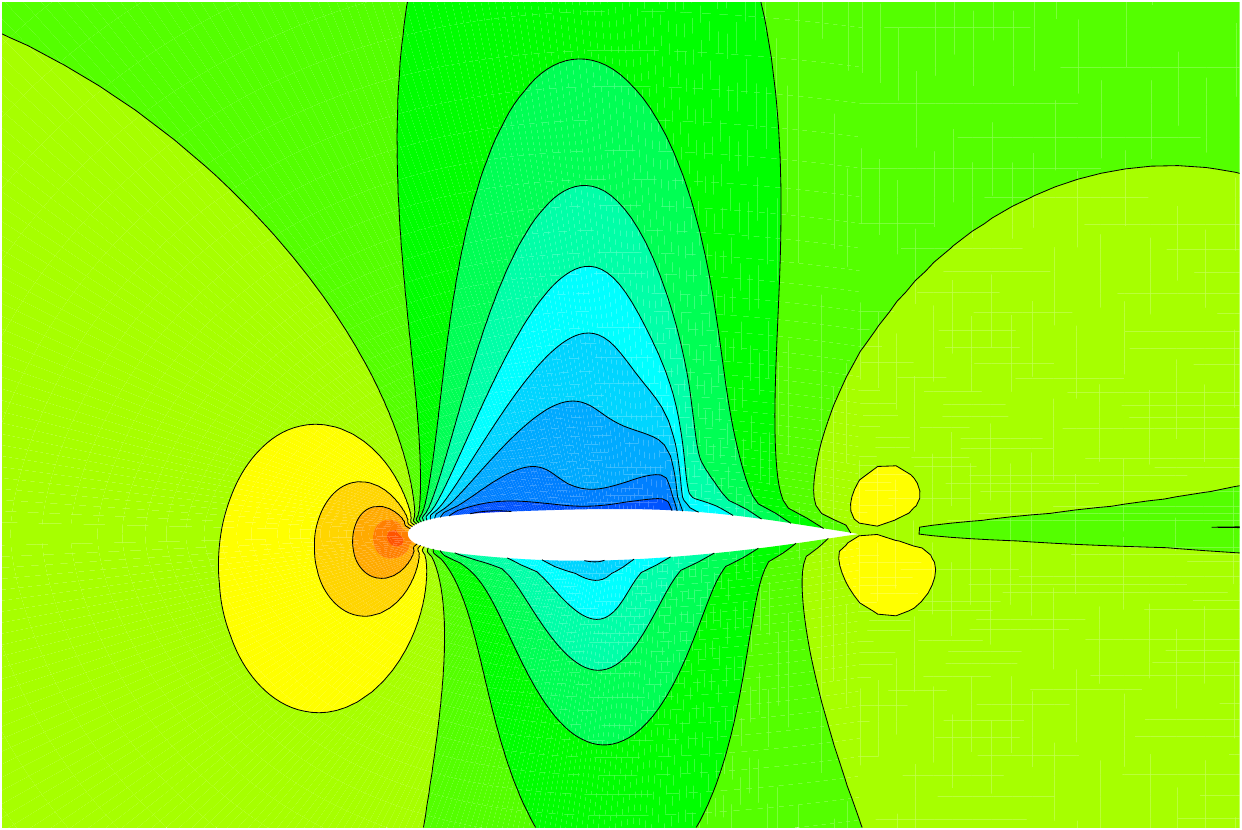}}
 \subfigure[2ndFV on mesh-2]{
 \includegraphics[width=6.0cm]{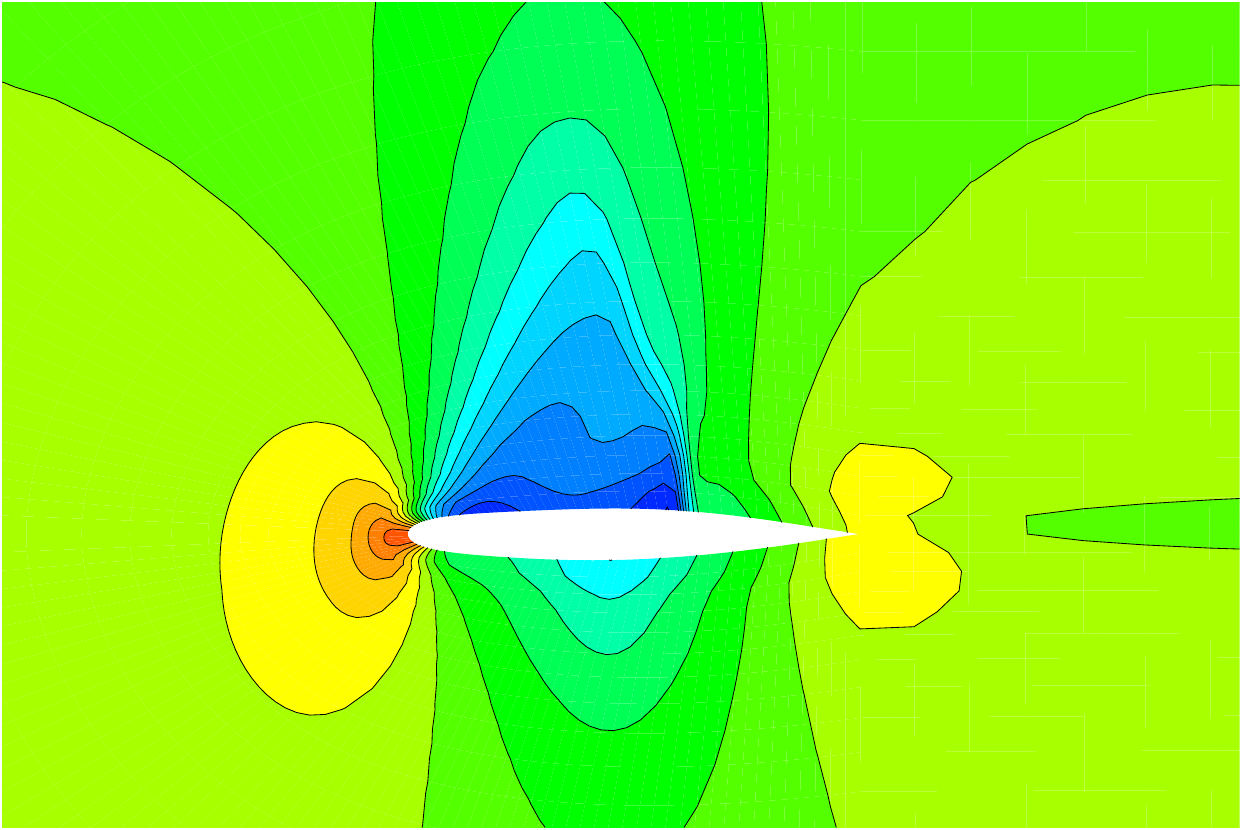}}
 \subfigure[DGp1 on mesh-2]{
 \includegraphics[width=6.0cm]{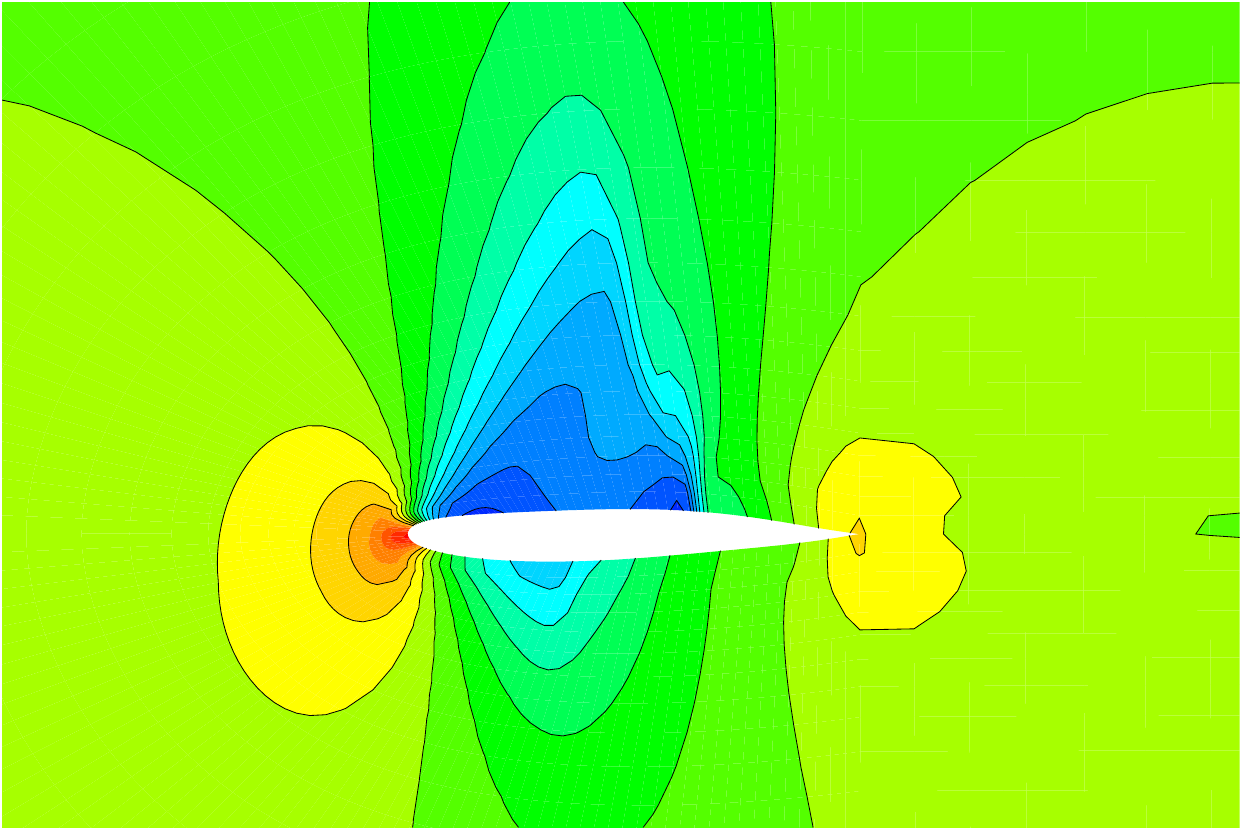}}
 \subfigure[DGp2 on mesh-1]{
 \includegraphics[width=6.0cm]{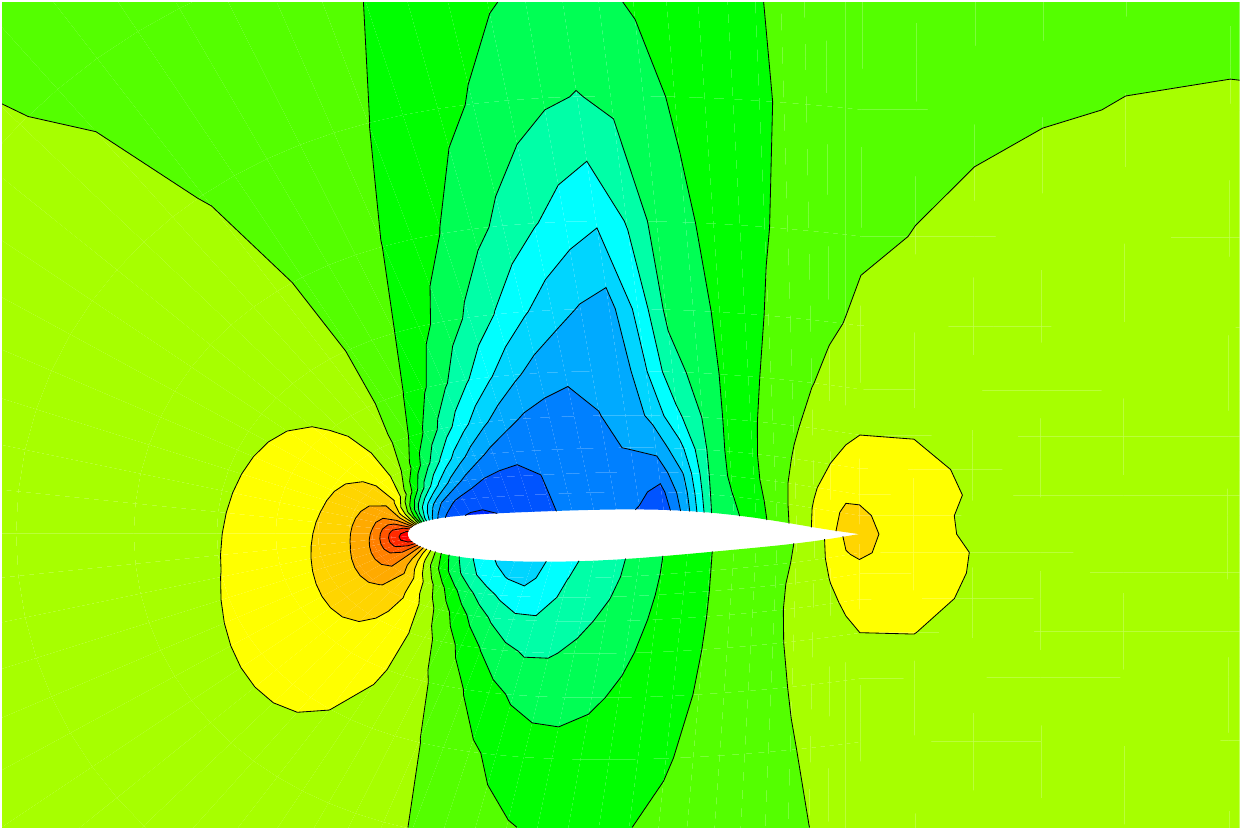}}
 \caption{Final density of different solvers and grids of the NACA-0012 airfoil test case.}
 \label{fig:naca0012flow-opt}
\end{figure}

\begin{figure}[htbp]
 \centering
 \subfigure[1stFV on mesh-3]{
 \includegraphics[width=6.0cm]{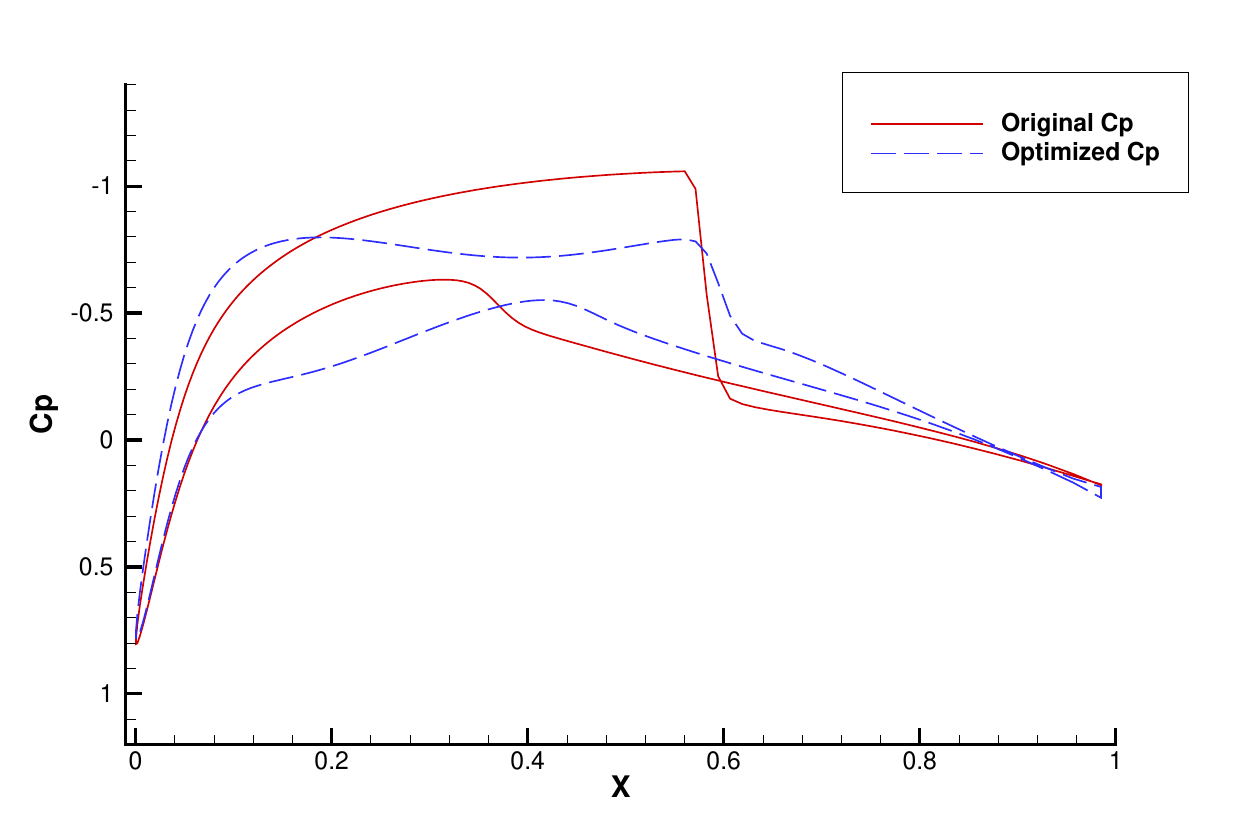}}
 \subfigure[2ndFV on mesh-2]{
 \includegraphics[width=6.0cm]{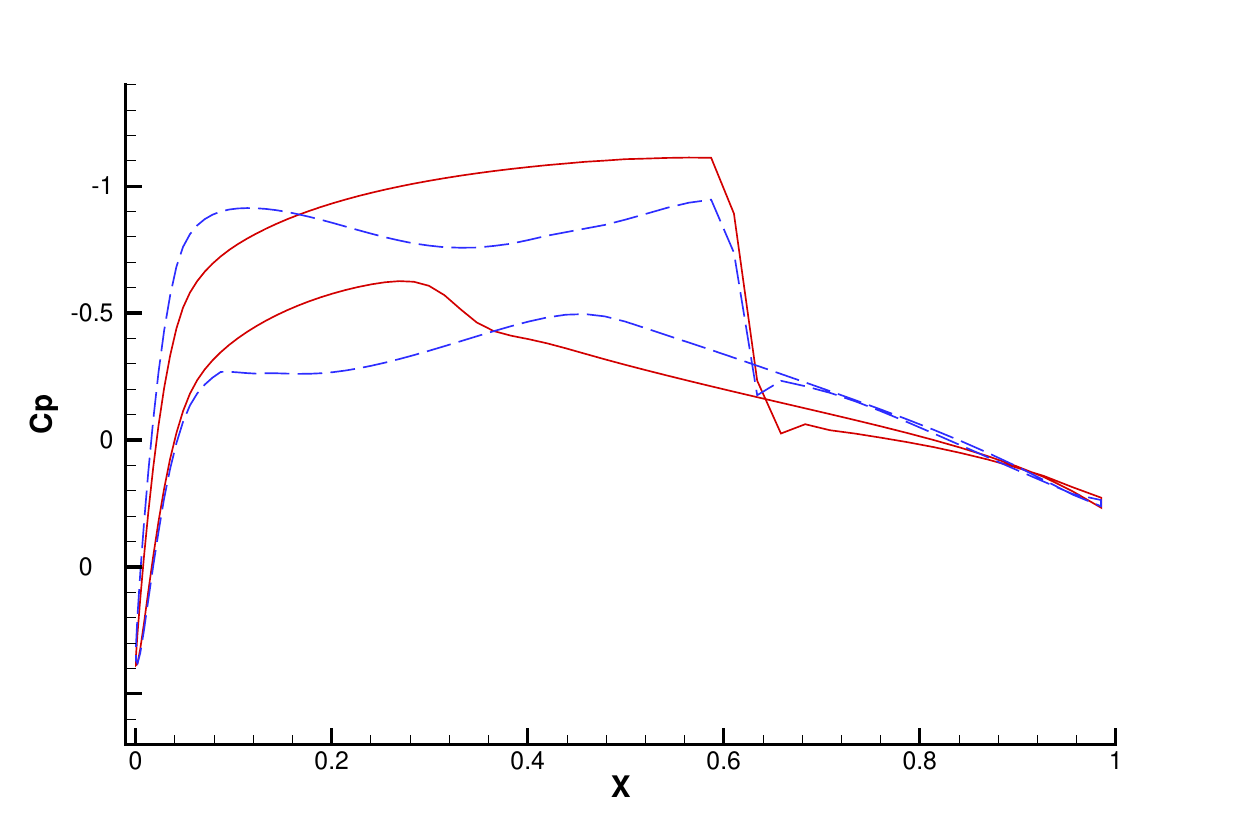}}
 \subfigure[DGp1 on mesh-2]{
 \includegraphics[width=6.0cm]{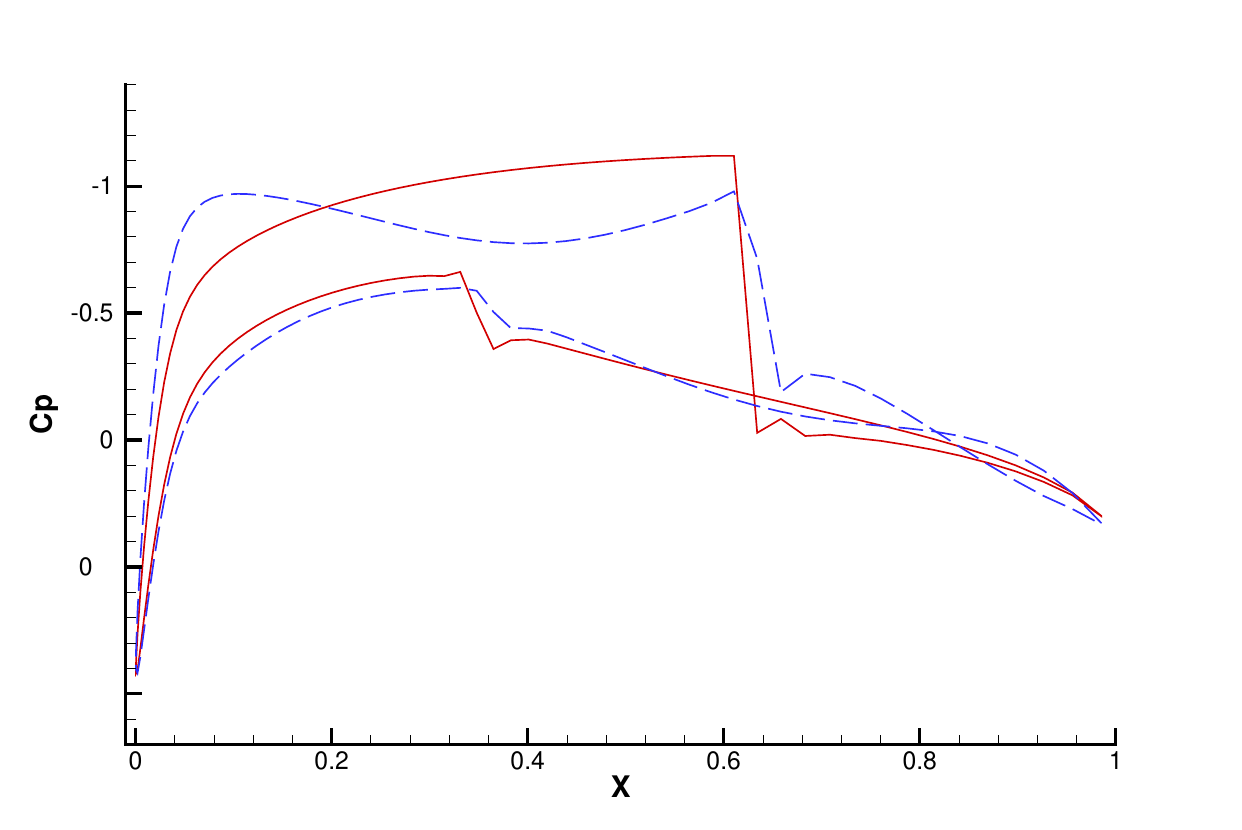}}
 \subfigure[DGp2 on mesh-1]{
 \includegraphics[width=6.0cm]{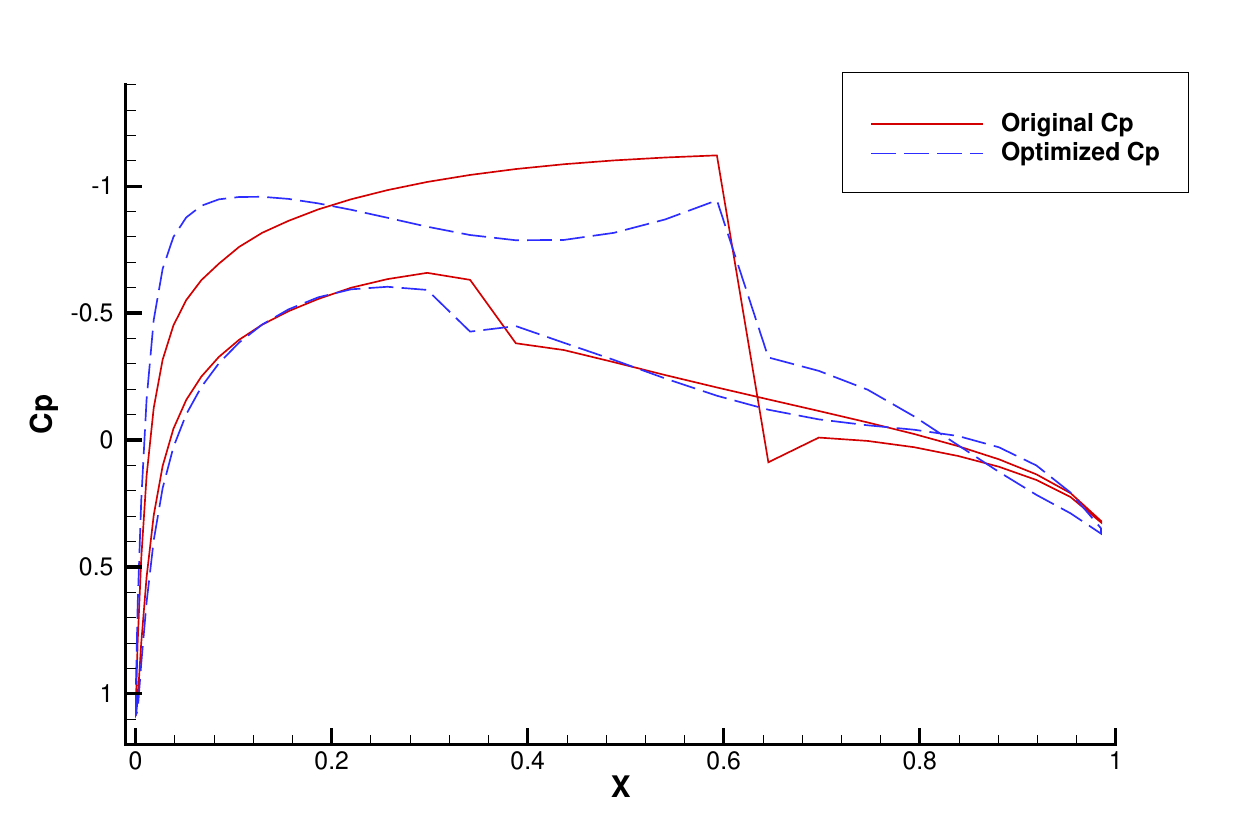}}
 \caption{Original and optimized Cp distribution of the NACA-0012 airfoil test case.}
 \label{fig:naca0012Cp-opt}
\end{figure}

The distinct of gradient calculation is further reflected in the optimized results as shown in Fig. \ref{fig:naca0012shape-opt}, \ref{fig:naca0012flow-opt} and \ref{fig:naca0012Cp-opt}, only the optimized results of DGp2 on mesh-1, DGp1 on mesh-2, 2ndFV on mesh-2, and 1stFV on mesh-3 are provided. At the same time, the iteration steps of optimization and the serial CPU run time in different cases are presented in Table \ref{tab:naca0012-time}, and the data of the aerodynamic performance improvement optimized by different CFD solvers are provided in Table \ref{tab:naca0012-opt}.

\begin{table}[htbp]
\centering \caption{Iteration step and serial CPU time (min) with/without DoFs remapping in different cases of the NACA-0012 airfoil test case}
\begin{center}
\begin{tabular}{c|c|c|c|c|c|c|c|c|c}
\toprule
 \multirow{2}*  & \multicolumn{3}{c|}{Iteration step} & \multicolumn{3}{c|}{\small{CPU time (no remapping)}} &\multicolumn{3}{c}{\small{CPU time (remapping)}}\\ \cline{2-10}
  &\small{mesh-1} &\small{mesh-2} &\small{mesh-3} &\small{mesh-1} &\small{mesh-2} &\small{mesh-3} &\small{mesh-1} &\small{mesh-2} &\small{mesh-3} \\ \hline
  1stFV &13 &8 &12 &0.50 &2.28 &34.00 &0.43 &1.57 &15.22 \\ \hline
  2ndFV &13 &17 &/ &5.97 &16.73 &/ &4.68 &8.58 &/  \\ \hline
  DGp1 &21 &19 &/ &12.68 &72.48 &/ &7.75 &23.48 &/ \\ \hline
  DGp2 &24 &/ &/ &48.53 &/ &/ &19.10 &/ &/ \\
\bottomrule
\end{tabular}
\end{center}
\label{tab:naca0012-time}
\end{table}

\begin{table}[htbp]
\centering \caption{Performance improvement of various CFD solvers in the NACA-0012 airfoil test case}
\begin{center}
\begin{tabular}{c c c c c c c c c c}
\toprule
     &\small{$C_{d0}$} &\small{$C_{d1}$} &\small{$\Delta C_{d}$} &\small{$C_{l0}$} &\small{$C_{l1}$} &\small{$\Delta C_{l}$} &\small{$A_{0}$} &\small{$A_{1}$} &\small{$\Delta A$} \\
\midrule
  \small{Reference\cite{yano2012case}} &\small{2.26e-2} &/ &/ &\small{3.52e-1} &/ &/ &/ &/ &/ \\
  \small{1stFV on mesh-3} &\small{3.38e-2} &\small{2.12e-2} &\small{-37.3$\%$} &\small{2.81e-1} &\small{2.81e-1} &\small{+0$\%$} &\small{8.24e-2} &\small{8.24e-2} &\small{+0$\%$} \\
  \small{2ndFV on mesh-2} &\small{2.30e-2} &\small{1.24e-2} &\small{-46.1$\%$} &\small{3.48e-1} &\small{3.48e-1} &\small{+0$\%$} &\small{8.24e-2} &\small{8.24e-2} &\small{+0$\%$} \\
  \small{DGp1 on mesh-2} &\small{2.23e-2} &\small{5.80e-3} &\small{-74.0$\%$} &\small{3.51e-1} &\small{3.52e-1} &\small{+0$\%$} &\small{8.24e-2} &\small{8.24e-2} &\small{+0$\%$} \\
  \small{DGp2 on mesh-1} &\small{2.27e-2} &\small{5.22e-3} &\small{-77.0$\%$} &\small{3.53e-1} &\small{3.53e-1} &\small{+0$\%$} &\small{8.23e-2} &\small{8.23e-2} &\small{+0$\%$} \\
\bottomrule
\end{tabular}
\end{center}
\label{tab:naca0012-opt}
\end{table}

It can be found from Table \ref{tab:naca0012-time} that under similar DoFs in CFD solvers, the CPU time cost by 2nd-order FV on mesh-2 is least, DGp2 on mesh-1 and 1st-order FV on mesh-3 is similar, and all less than that of DGp1 on mesh-2, while in most cases, the application of DoFs remapping techniques can reduce the CPU time of the optimization by around $50\%\sim70\%$; from Fig. \ref{fig:naca0012shape-opt}-\ref{fig:naca0012Cp-opt} and Table \ref{tab:naca0012-opt}, it can be found that under coequal computational costs, DGp1 on mesh-2 and DGp2 on mesh-1 produce similar results of the optimized shape, flow distribution, and the Cp distribution, and both can achieve drag reduction of more than $70\%$. The optimized results produced by DGp1 and DGp2 are significantly superior to that produced by FVMs. The reason might be related to that the more sensitive adjoint gradients computed by DGMs extend the capability to explore superior optimal solutions during the optimization process.

\subsection{RAE-2822 airfoil drag minimization}
\qquad We now consider the RAE-2822 drag minimization, the target is still to reduce the drag $C_d$ without decreasing the lift $C_l$ and the area $A$ of the airfoil. The free stream is set to $\text{Ma}=0.726, \text{AOA}=2.44^o$. Still only the changes of y-coordinates are taken into consideration, the number of the design variables is set to 32, and their ranges of variation are set no more than $15\%$.

\begin{figure}[htbp]
 \centering
 \subfigure[mesh-1]{
 \includegraphics[width=4.8cm]{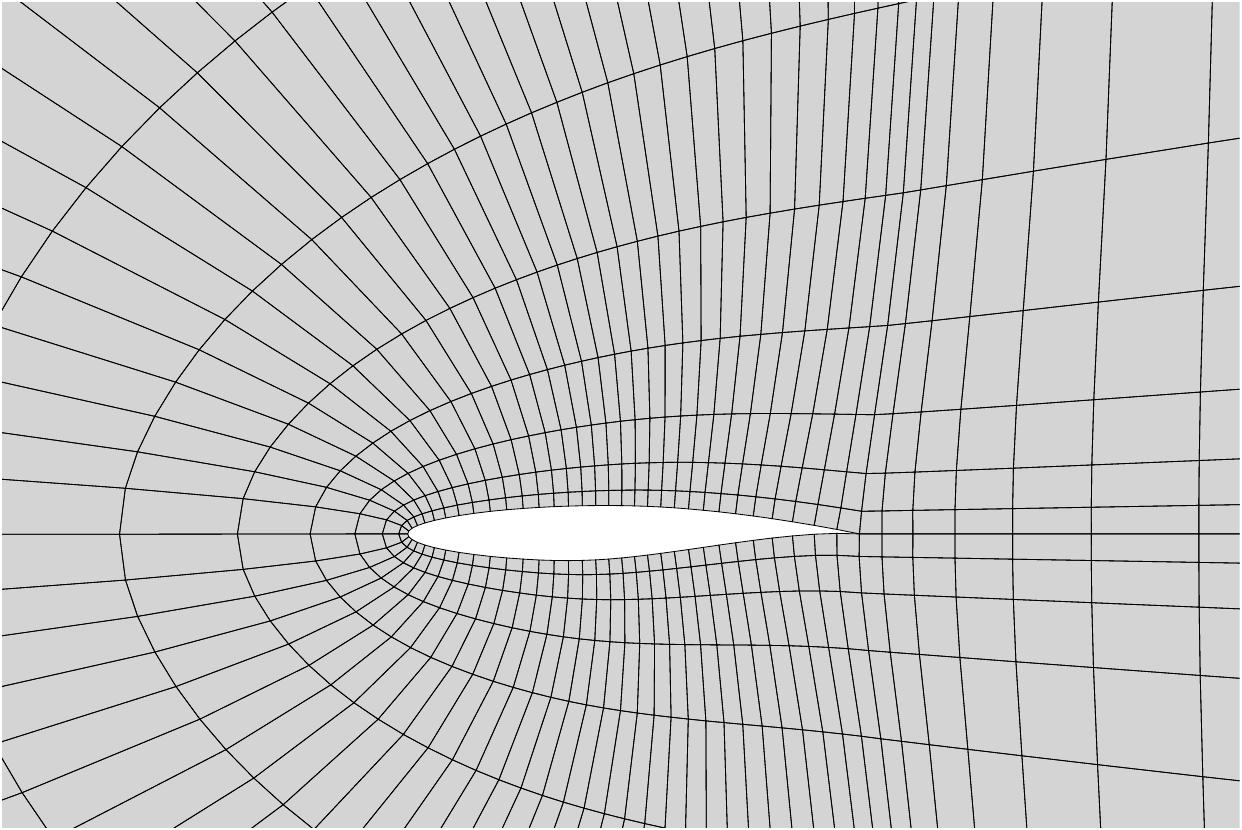}}
 \subfigure[mesh-2]{
 \includegraphics[width=4.8cm]{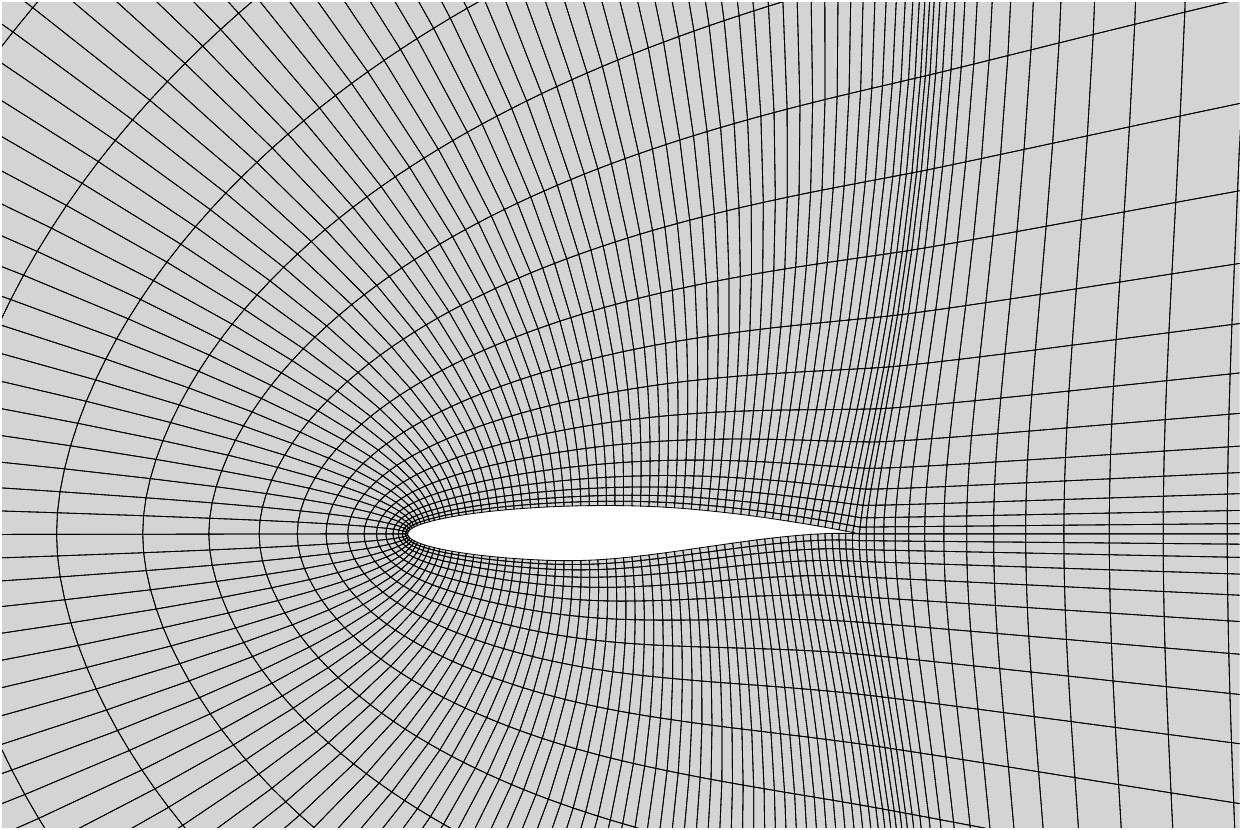}}
 \subfigure[mesh-3]{
 \includegraphics[width=4.8cm]{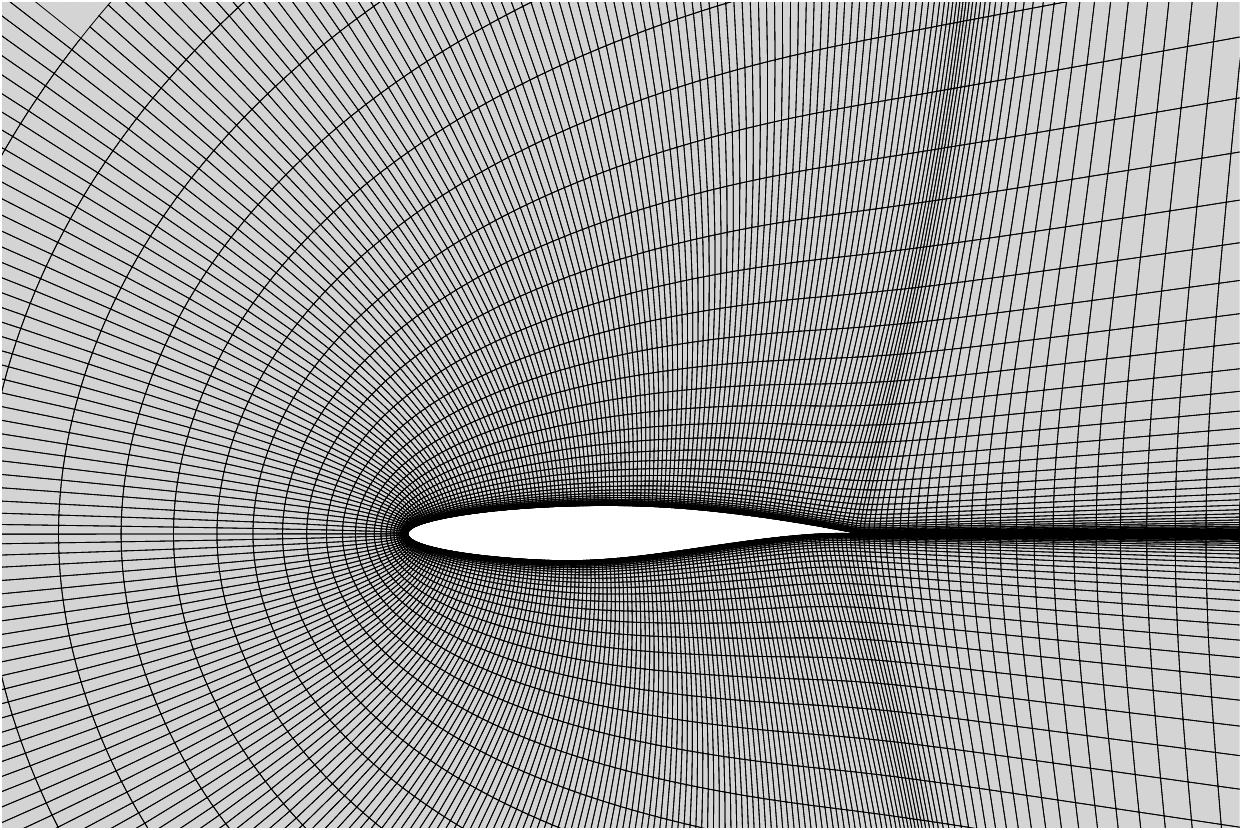}}
 \caption{Computational grids with different levels of the RAE-2822 airfoil test case.}
 \label{fig:rae2822mesh}
\end{figure}

\begin{figure}[htbp]
 \centering
 \subfigure[1stFV on mesh-3]{
 \includegraphics[width=4.8cm]{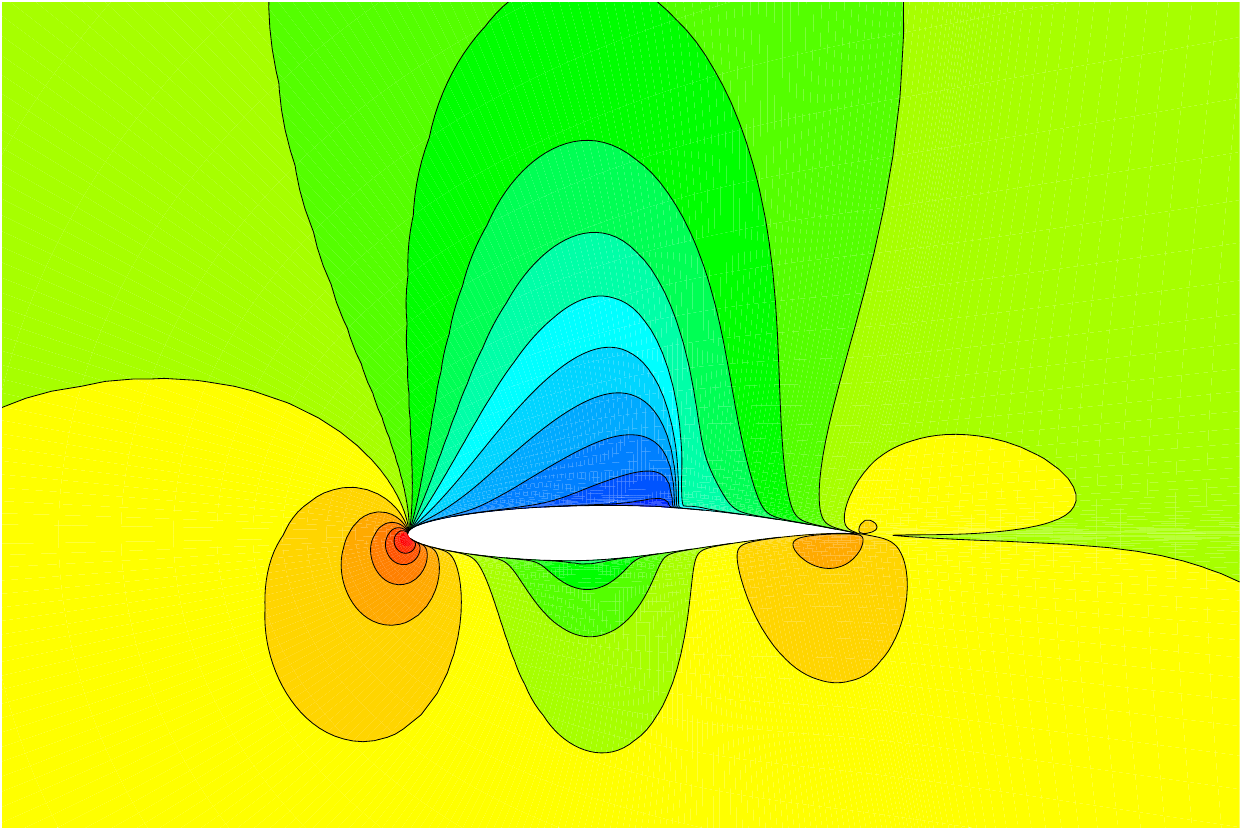}}
 \subfigure[2ndFV on mesh-1]{
 \includegraphics[width=4.8cm]{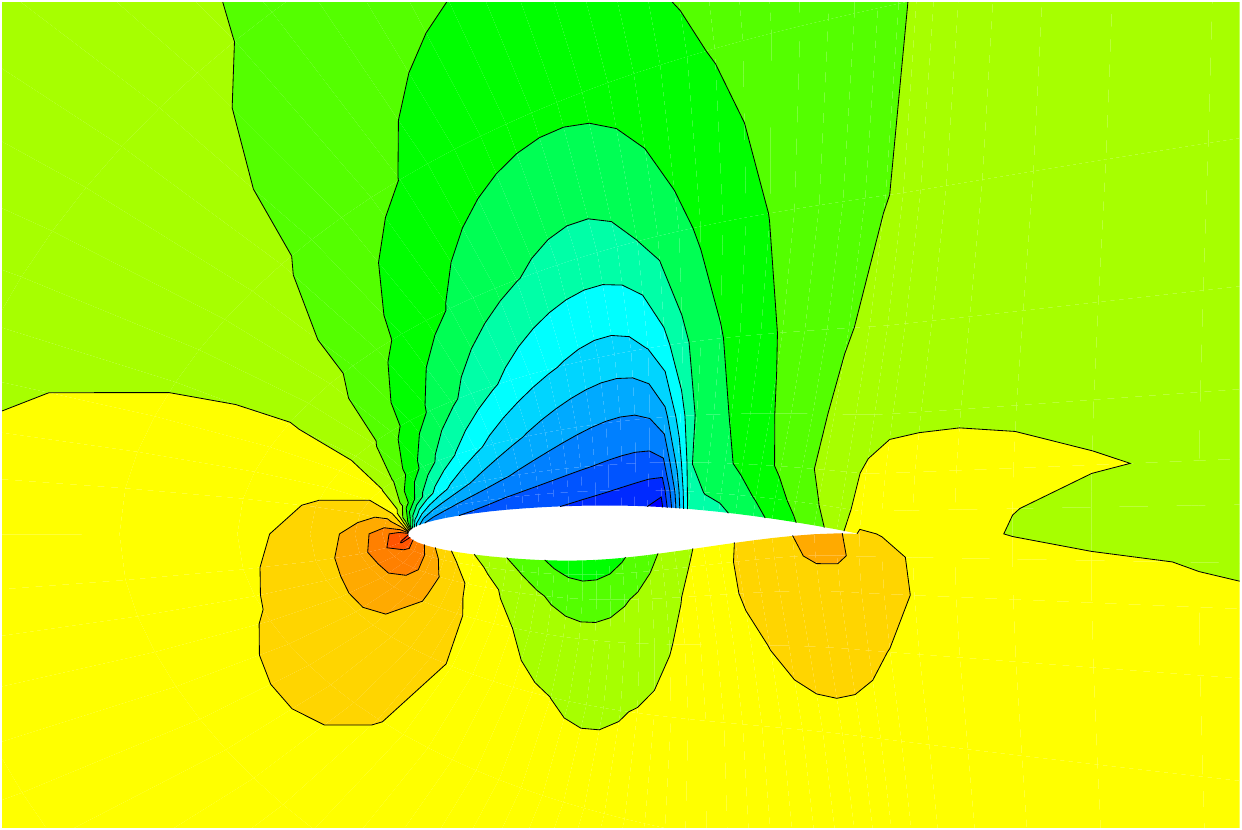}}
 \subfigure[2ndFV on mesh-2]{
 \includegraphics[width=4.8cm]{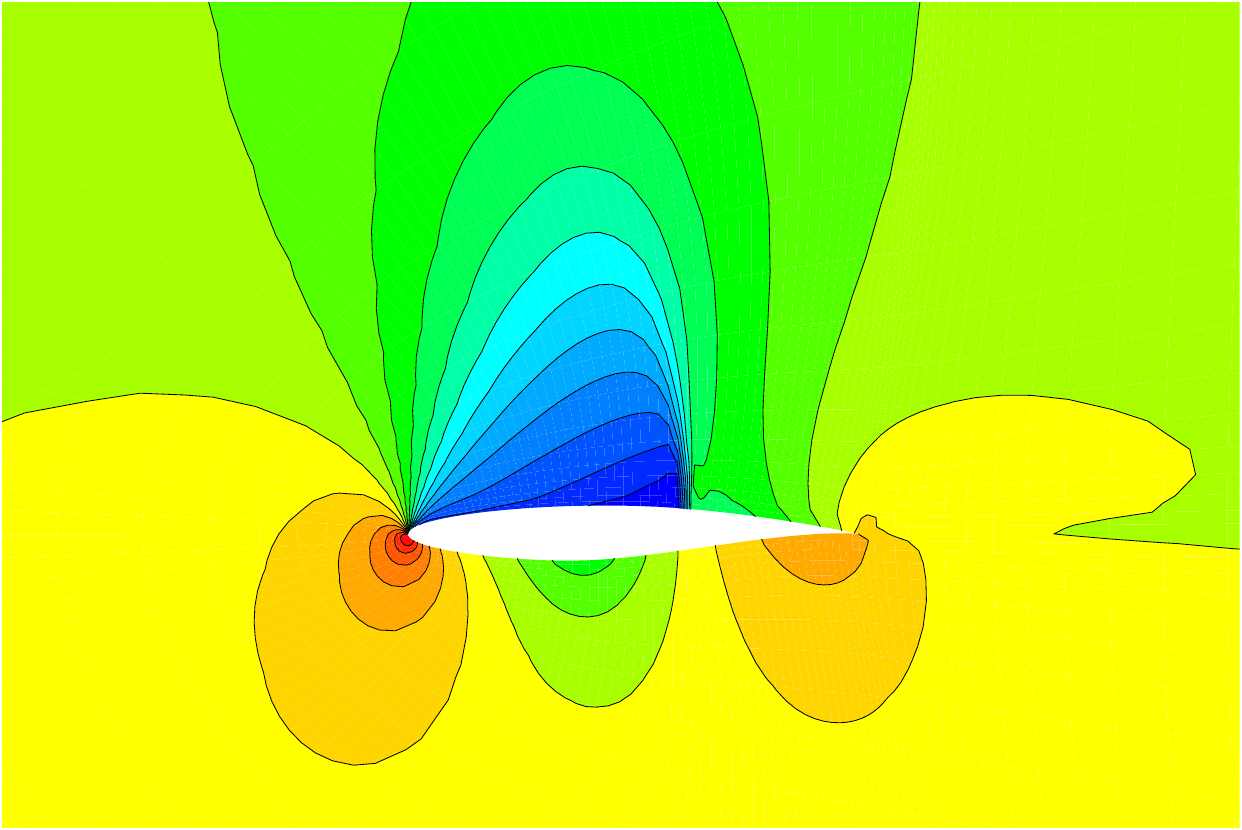}}
 \subfigure[DGp1 on mesh-1]{
 \includegraphics[width=4.8cm]{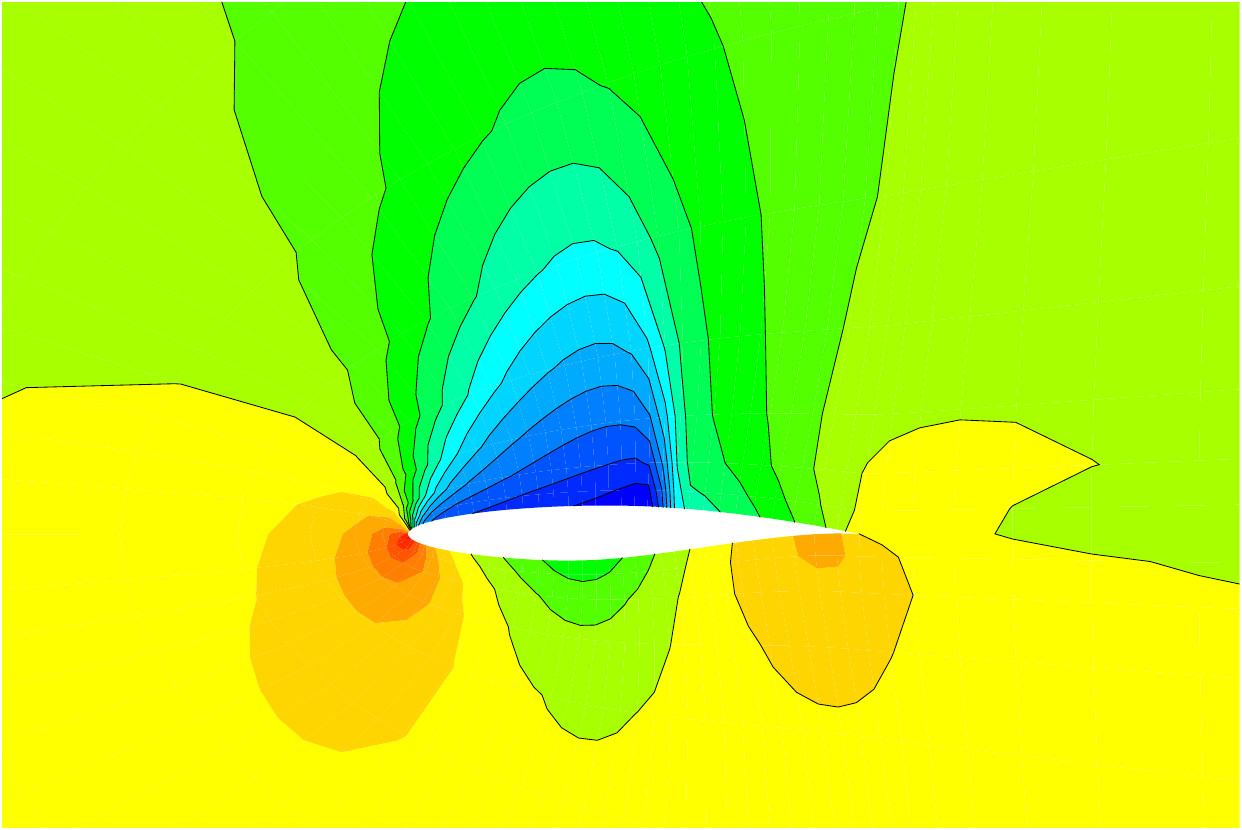}}
 \subfigure[DGp1 on mesh-2]{
 \includegraphics[width=4.8cm]{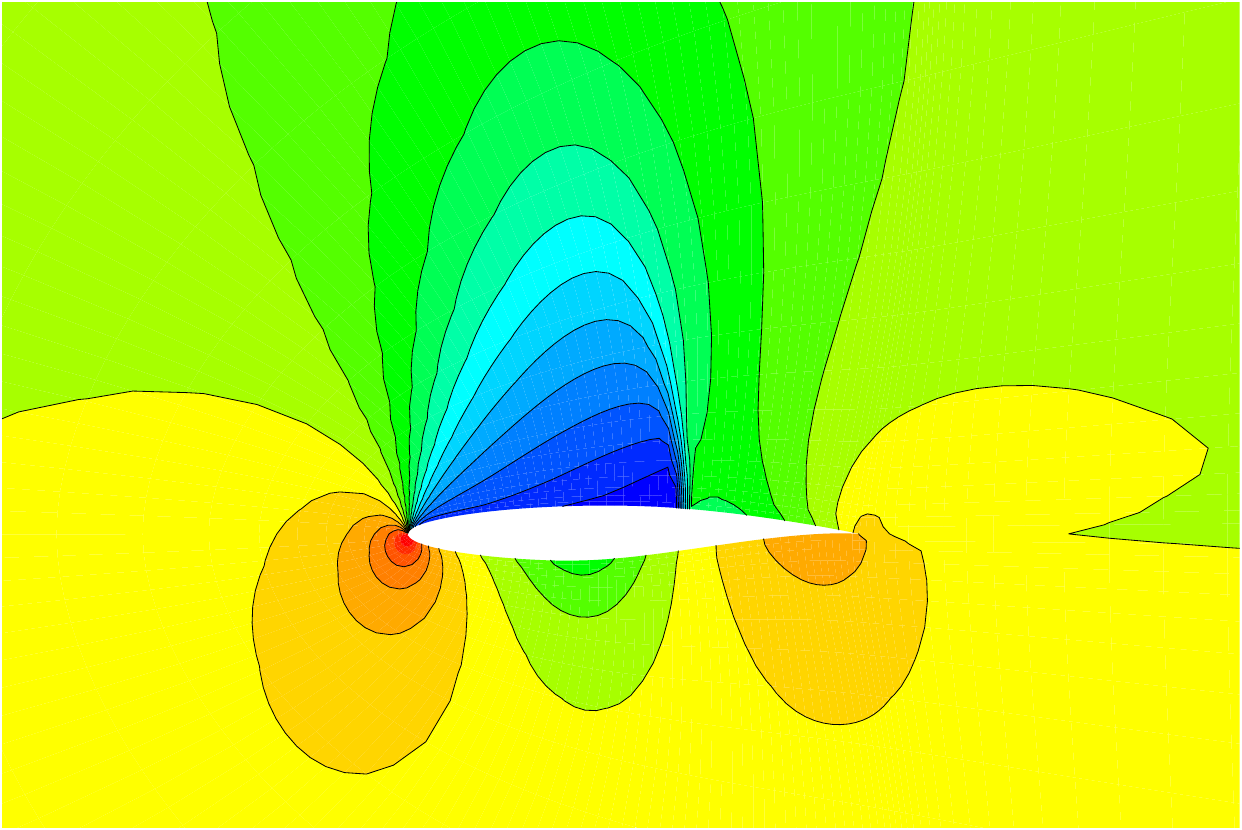}}
 \subfigure[DGp2 on mesh-1]{
 \includegraphics[width=4.8cm]{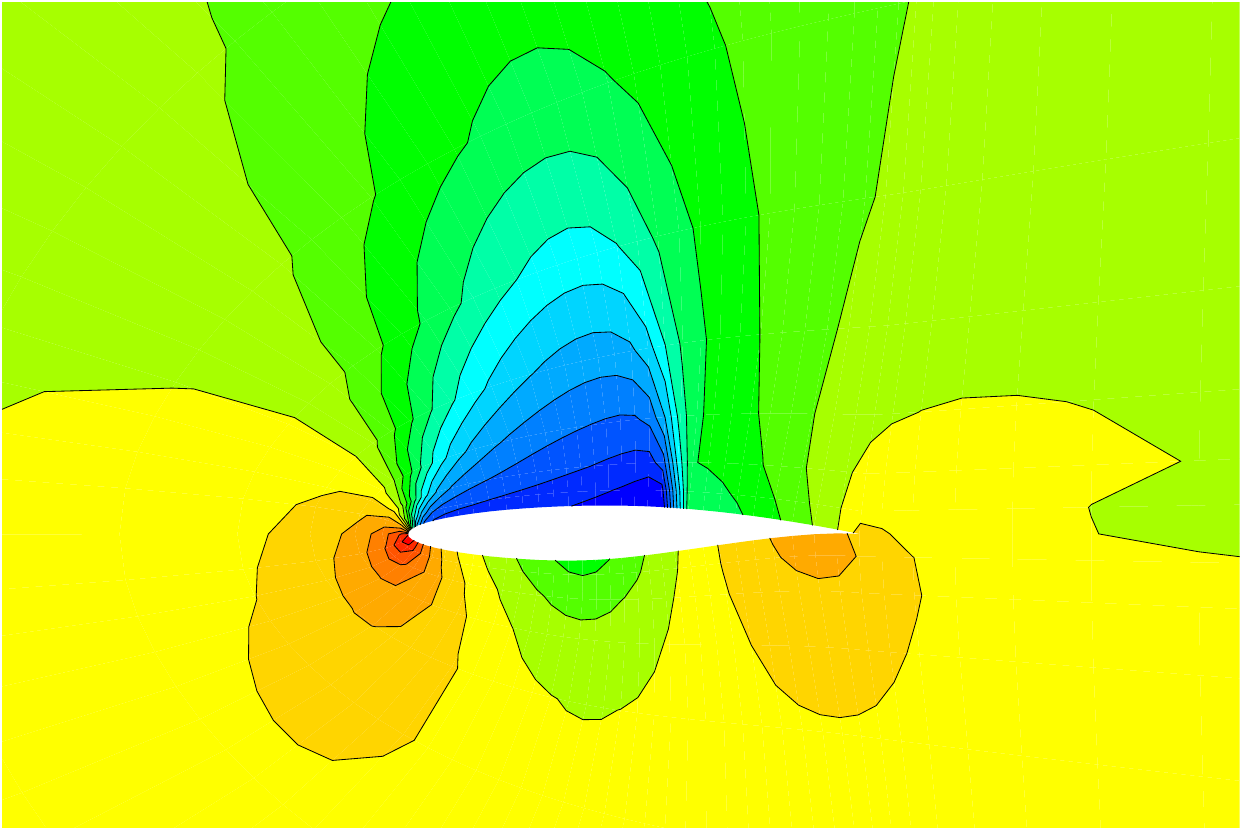}}
 \caption{Solutions (density) obtained by different solvers on various grids of the RAE-2822 airfoil test case.}
 \label{fig:rae2822flow}
\end{figure}

\begin{figure}[htbp]
  \centering
  \includegraphics[width=9.5cm]{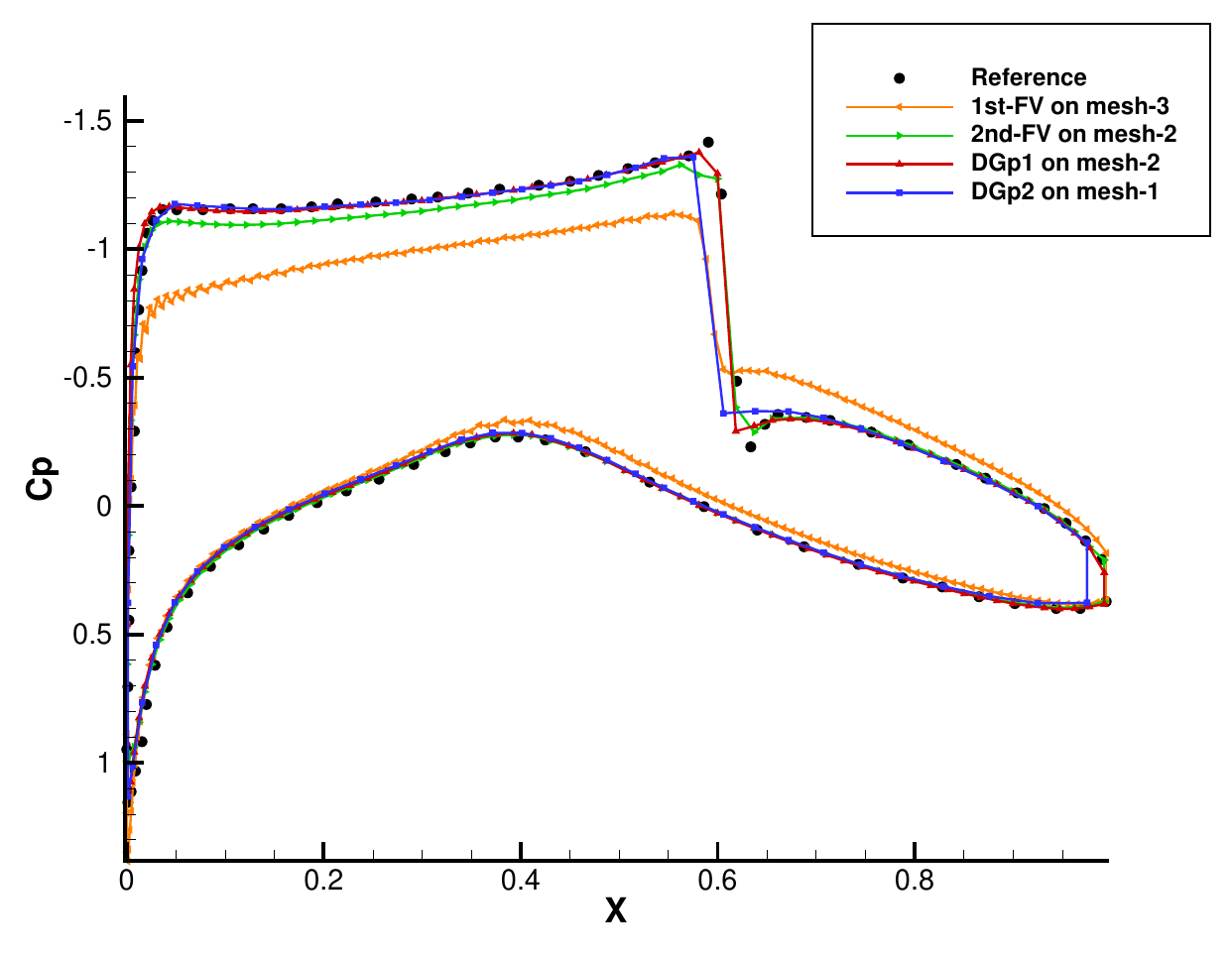}
  \caption{RAE-2822 $C_p$ distributions between various CFD solvers.}
  \label{fig:rae2822cp}
\end{figure}

Similarly, 3 grids (mesh-1, mesh-2, mesh-3) are used in this test case to test the performance of different CFD solvers (as shown in Fig. \ref{fig:rae2822mesh}), the iterative DoFs of different CFD solvers on each grid is presented in Table \ref{tab:rae2822-DoFs}, and the corresponding flow results obtained by various CFD solvers (1stFV, 2ndFV, DGp1, DGp2) are presented in Fig. \ref{fig:rae2822flow}, the $C_p$ distributions obtained by various CFD solvers are compared with reference data \cite{wang2019adjoint} in Fig. \ref{fig:rae2822cp}. It can be observed that similar and acceptable flow results can be reached through DGp2 on mesh-1, DGp1 and 2ndFV on mesh-2, and 1stFV on mesh-3.

\begin{table}[htbp]
\centering \caption{DoFs of various CFD solvers on each grid in the RAE-2822 test case}
\setlength{\tabcolsep}{8mm}
\begin{center}
\begin{tabular}{c|c|c|c}
\toprule
  & mesh-1 & mesh-2  & mesh-3 \\ \hline
  FVMs  &12,288 &49,152 &196,608 \\ \hline
  DGp1 &36,864 &147,456 &/   \\ \hline
  DGp2 &73,728 &/       &/  \\
\bottomrule
\end{tabular}
\end{center}
\label{tab:rae2822-DoFs}
\end{table}

\begin{figure}[htbp]
  \centering
  \includegraphics[width=10.0cm]{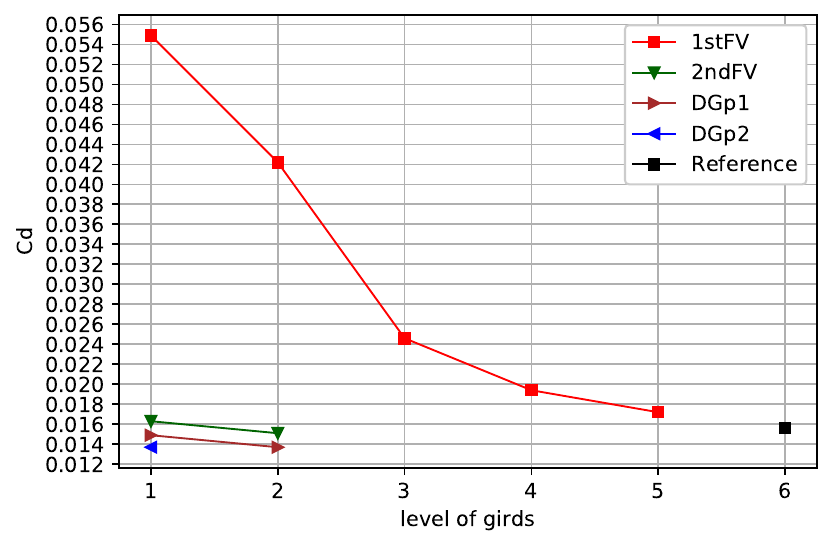}
  \caption{$C_d$ convergence of the RAE-2822 airfoil test case.}
  \label{fig:rae2822-obj}
\end{figure}

\begin{figure}[htbp]
  \centering
  \subfigure[sensitivity convergence of FVMs]{
  \includegraphics[width=6.8cm]{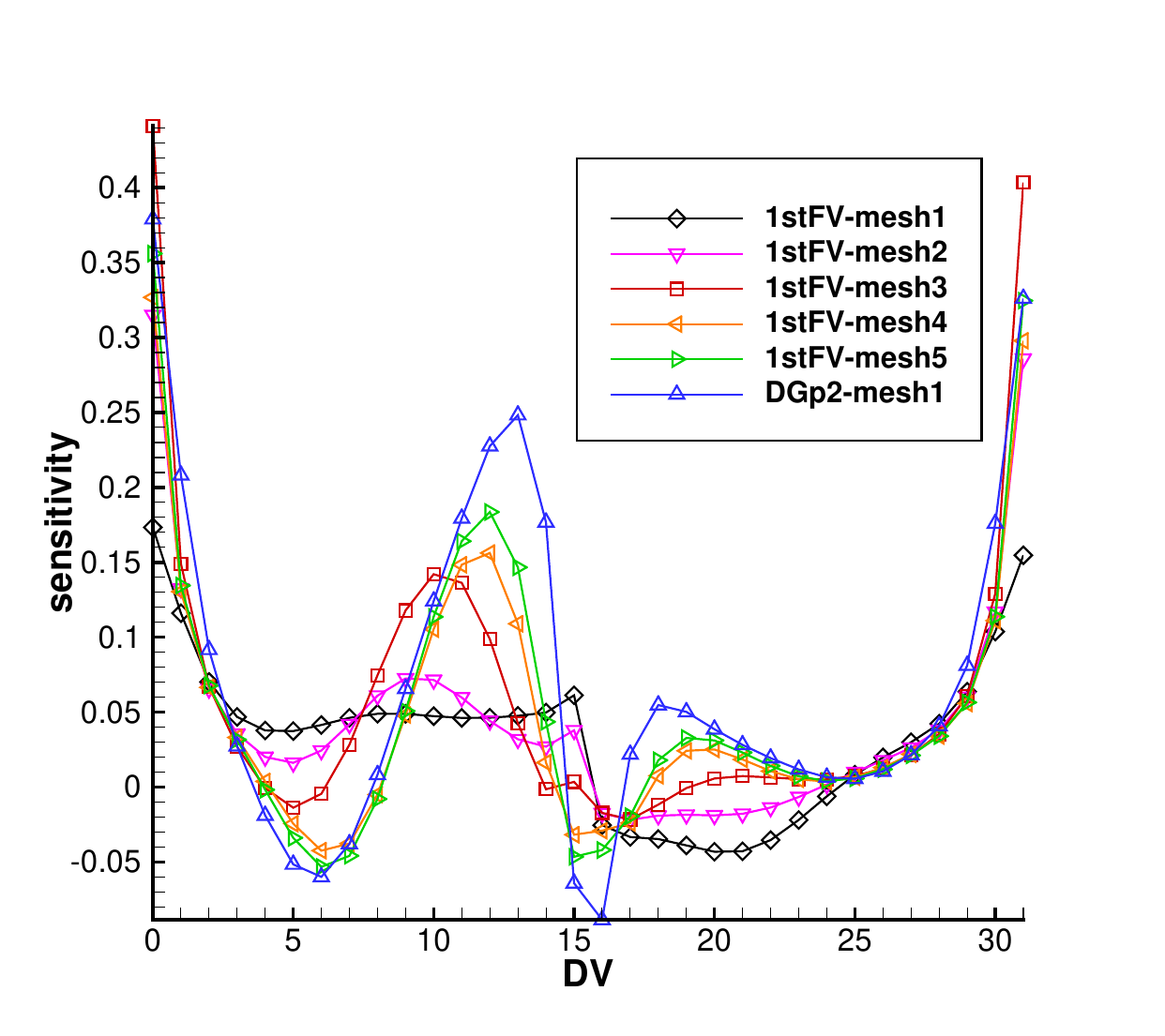}}
  \subfigure[sensitivity comparison under coequal DoFs]{
  \includegraphics[width=6.8cm]{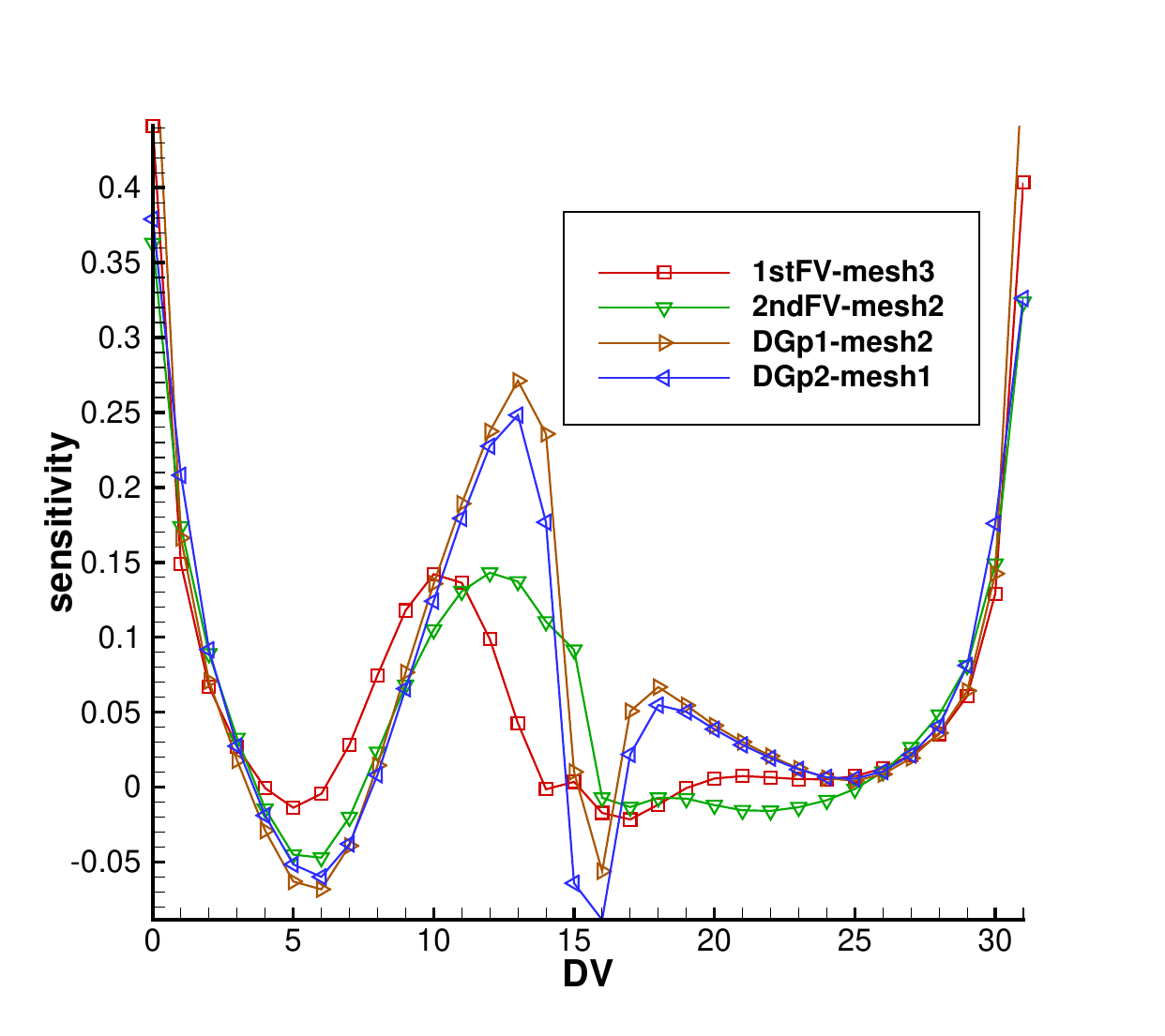}}
  \caption{Drag sensitivity of the RAE-2822 airfoil test case}
  \label{fig:rae2822grad}
\end{figure}

The results of the $C_d$ convergence and the adjoint gradients (sensitivity) computed by different CFD solvers are shown in Fig. \ref{fig:rae2822-obj} and \ref{fig:rae2822grad}, respectively. Similar phenomenon can be observed from the Fig. \ref{fig:rae2822-obj} and \ref{fig:rae2822grad} that the $C_d$ computed by 1st-order FVMs converges with the gird refinement, but requires a rather fine gird (mesh-5) to achieve a similar accuracy of $C_d$ as compared with DGp1 on mesh-2 or DGp2 on mesh-1; the adjoint gradients computed by the DGMs on coarse grids (mesh-1 and mesh-2) are similar, and more sensitive to certain design variables than that computed by the FVMs.

\begin{figure}[htbp]
 \centering
 \subfigure[1stFV on mesh-3]{
 \includegraphics[width=6.0cm]{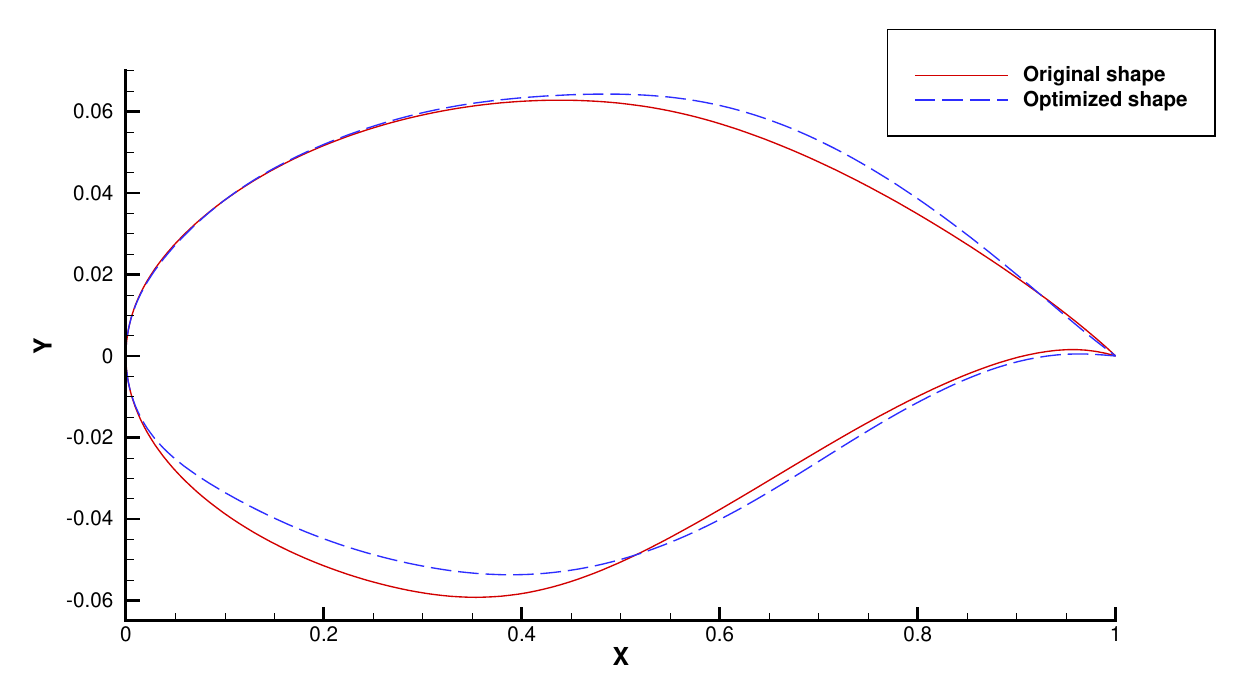}}
 \subfigure[2ndFV on mesh-2]{
 \includegraphics[width=6.0cm]{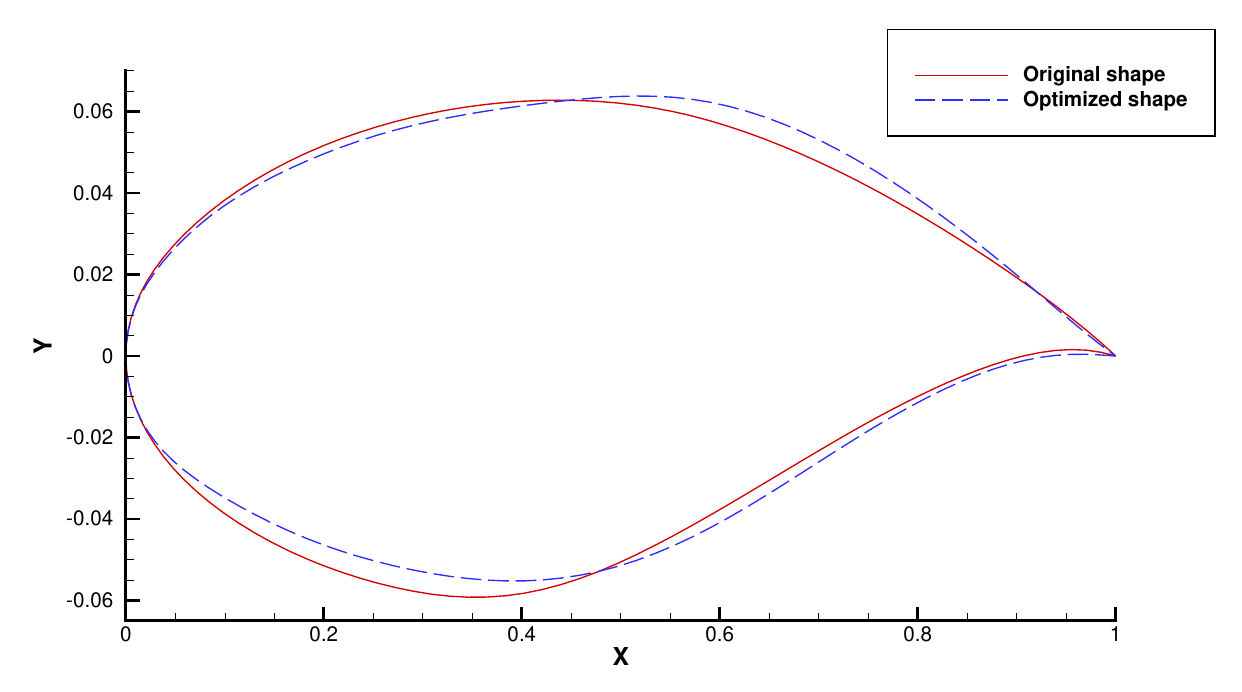}}
 \subfigure[DGp1 on mesh-2]{
 \includegraphics[width=6.0cm]{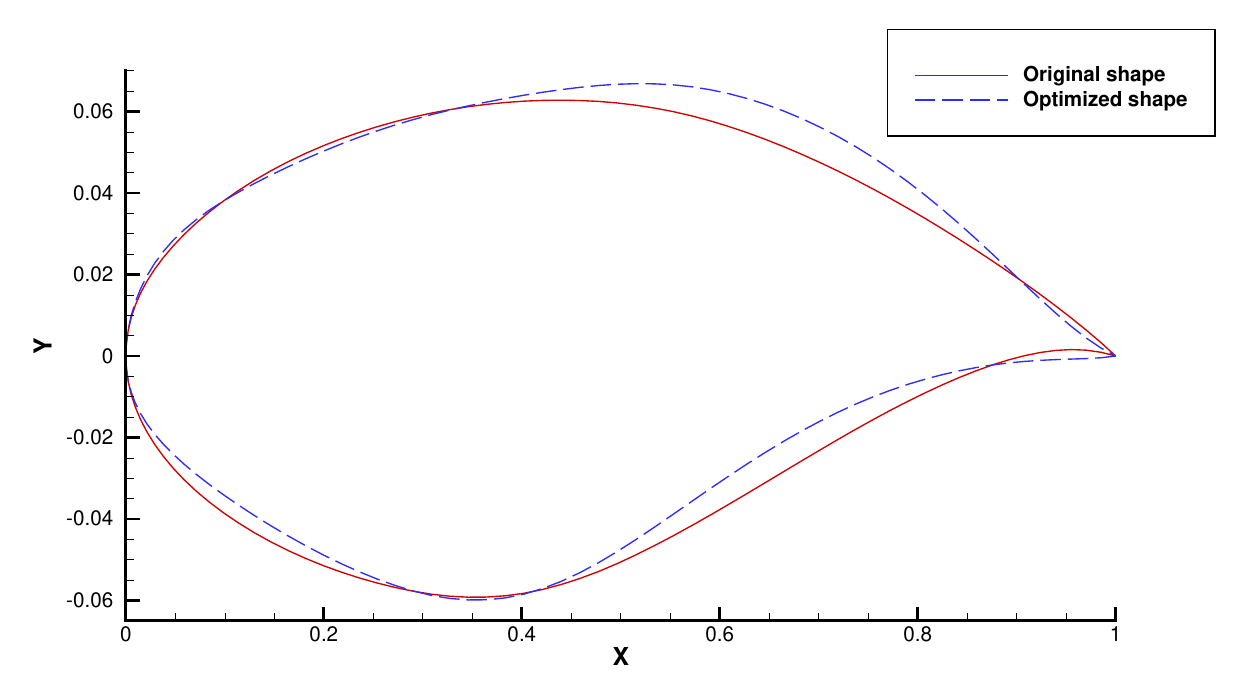}}
 \subfigure[DGp2 on mesh-1]{
 \includegraphics[width=6.0cm]{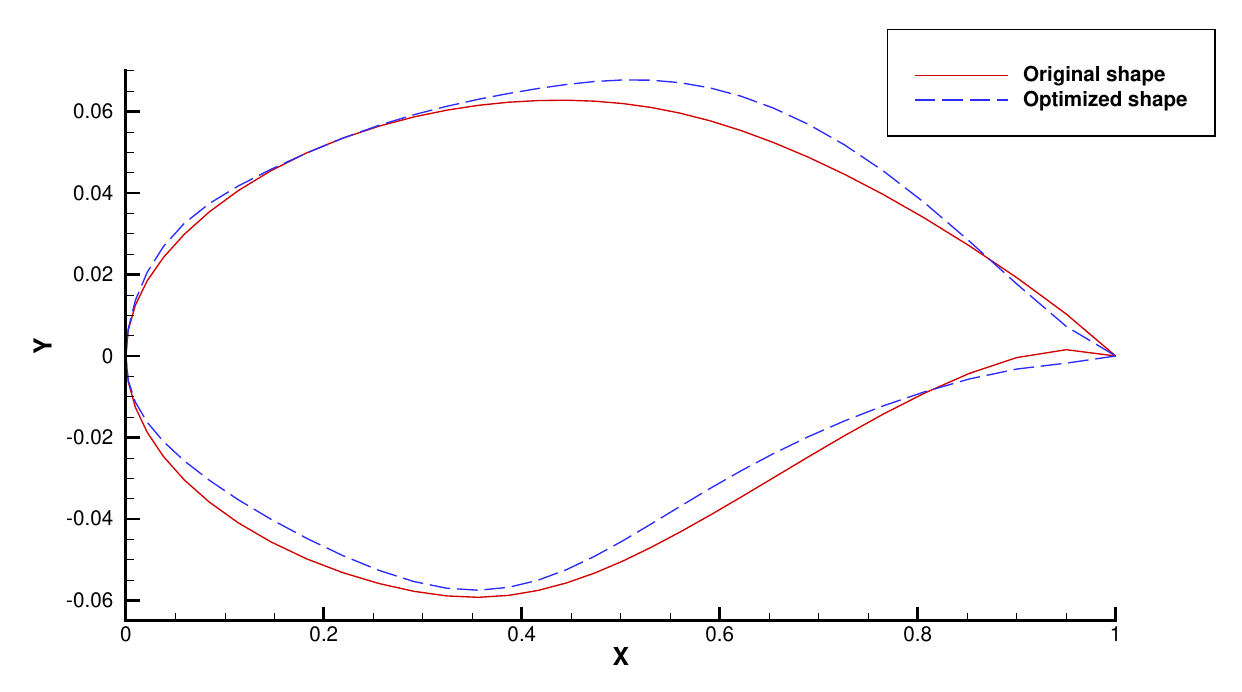}}
 \caption{Original and optimized shape of the RAE-2822 airfoil test case.}
 \label{fig:rae2822shape-opt}
\end{figure}

\begin{figure}[htbp]
 \centering
 \subfigure[1stFV on mesh-3]{
 \includegraphics[width=6.0cm]{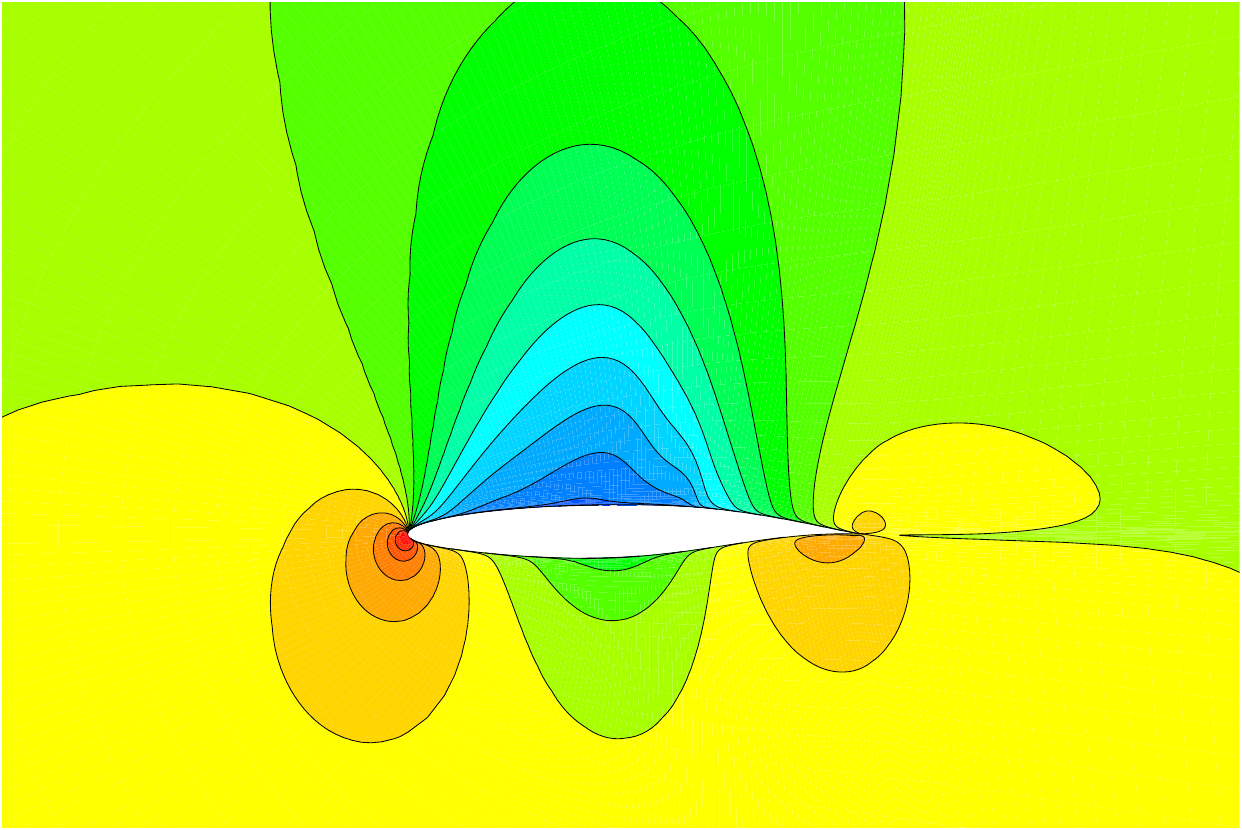}}
 \subfigure[2ndFV on mesh-2]{
 \includegraphics[width=6.0cm]{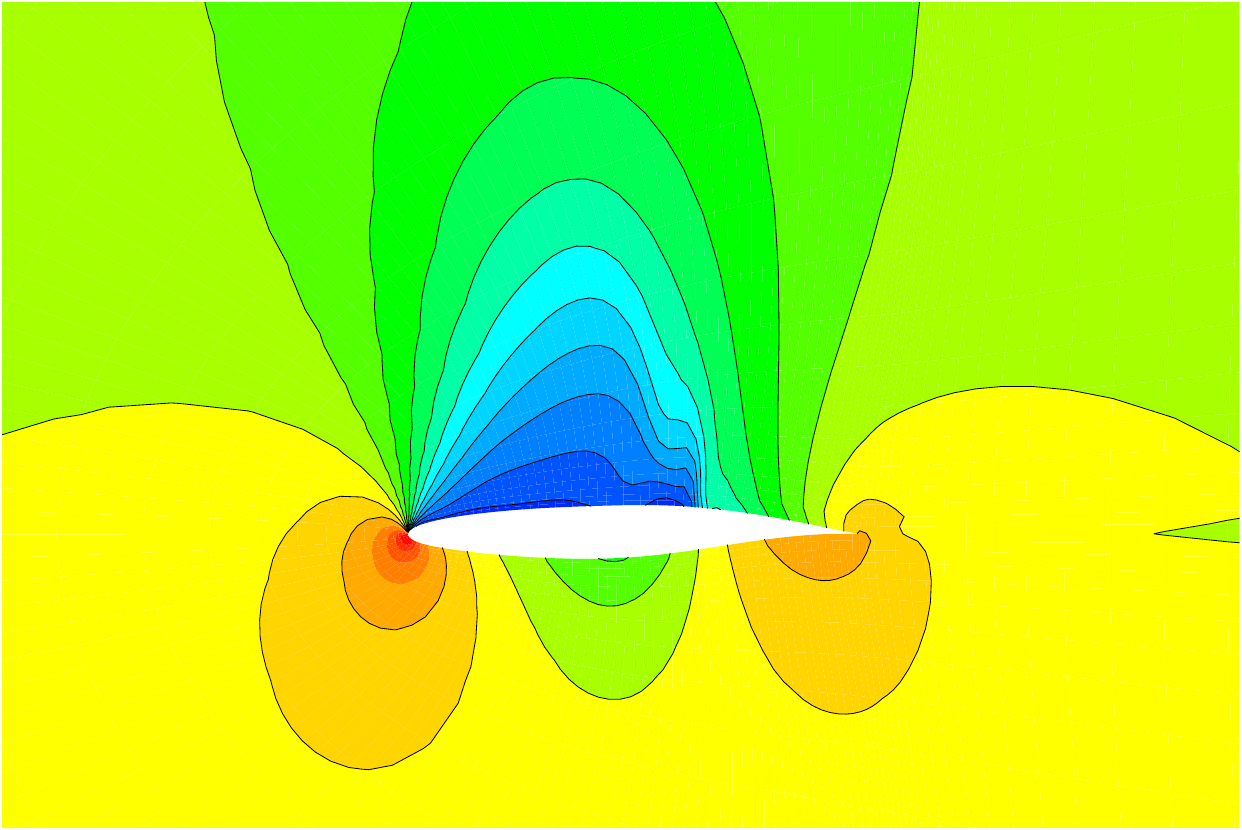}}
 \subfigure[DGp1 on mesh-2]{
 \includegraphics[width=6.0cm]{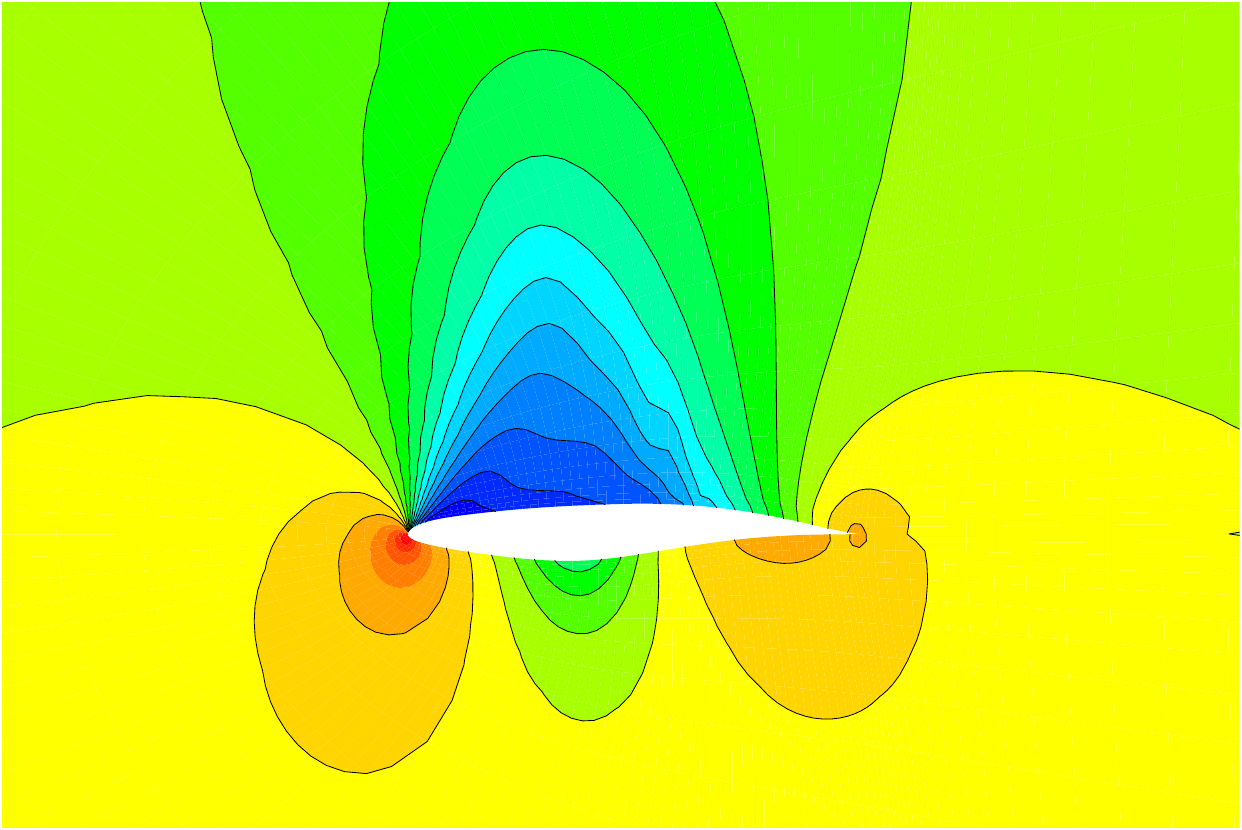}}
 \subfigure[DGp2 on mesh-1]{
 \includegraphics[width=6.0cm]{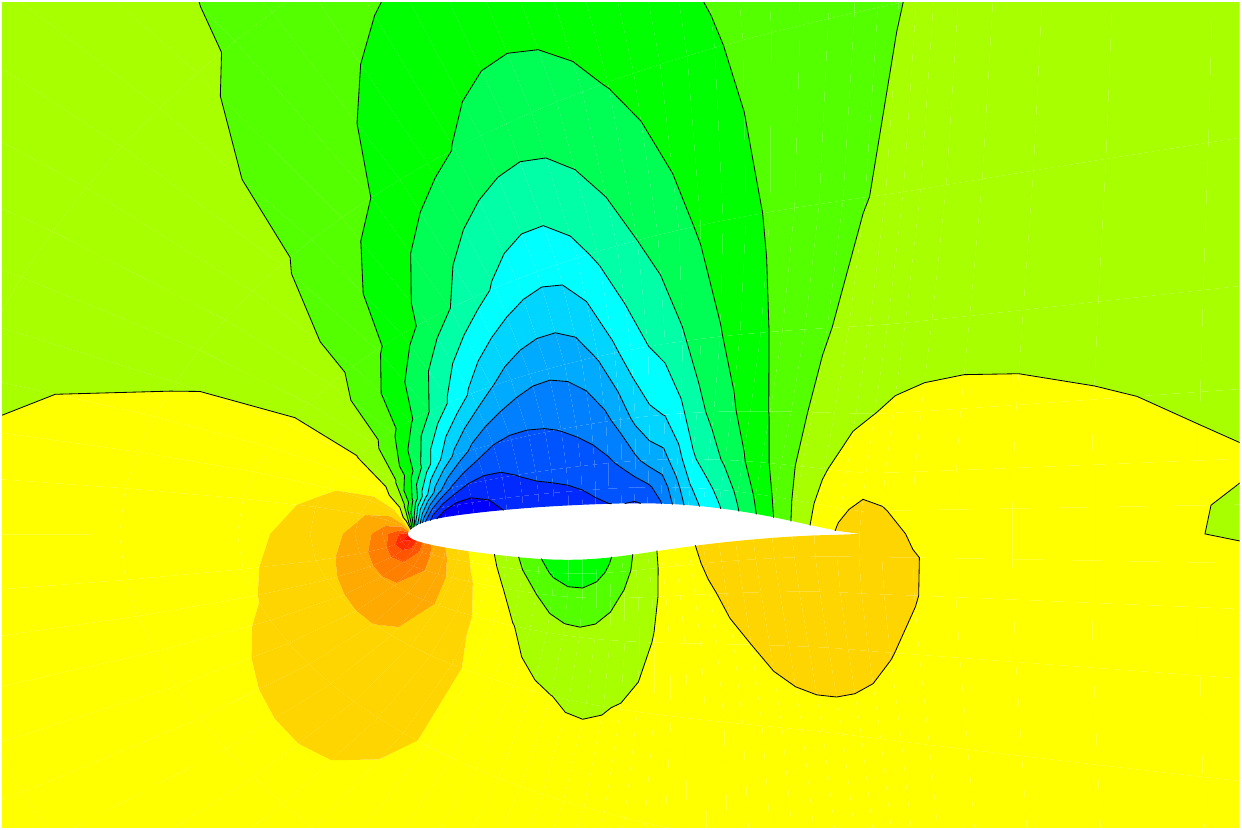}}
 \caption{Final density of different solvers and grids of the RAE-2822 airfoil test case.}
 \label{fig:rae2822flow-opt}
\end{figure}

\begin{figure}[htbp]
 \centering
 \subfigure[1stFV on mesh-3]{
 \includegraphics[width=6.0cm]{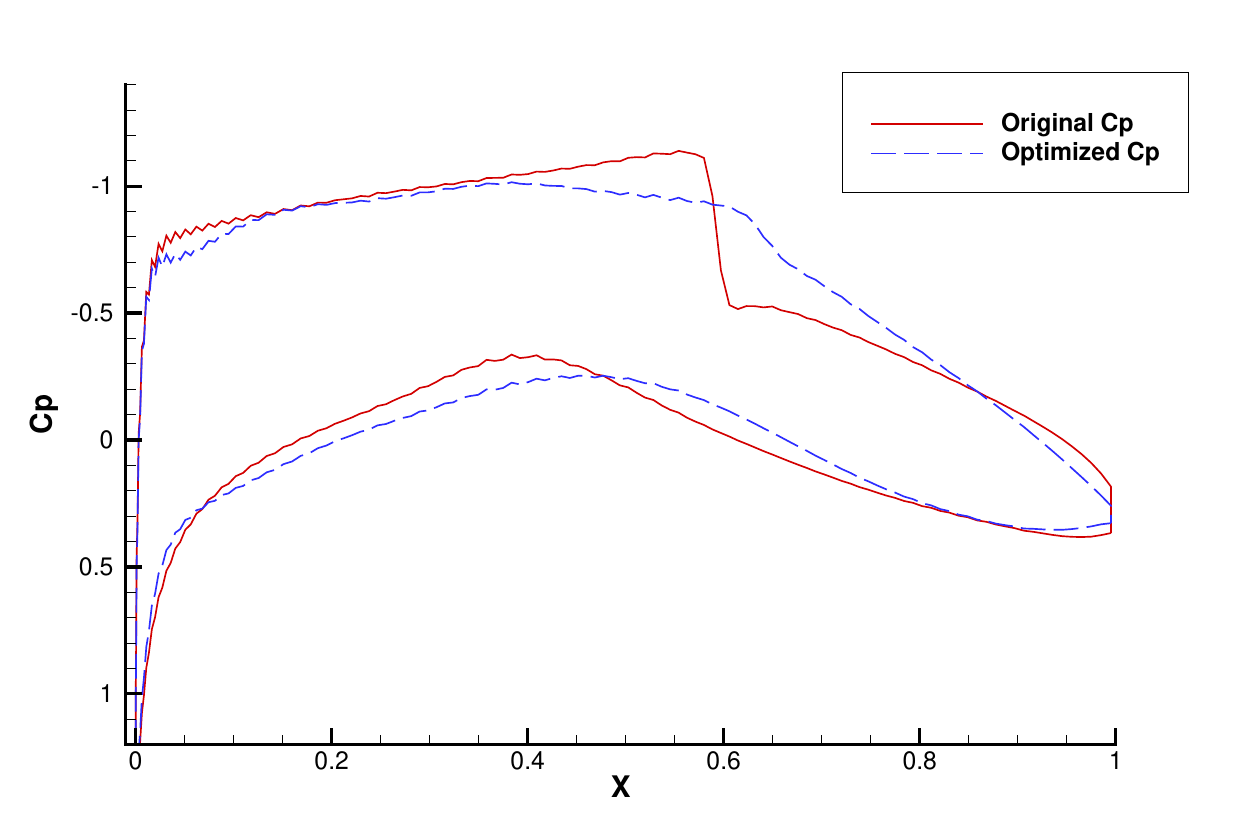}}
 \subfigure[2ndFV on mesh-2]{
 \includegraphics[width=6.0cm]{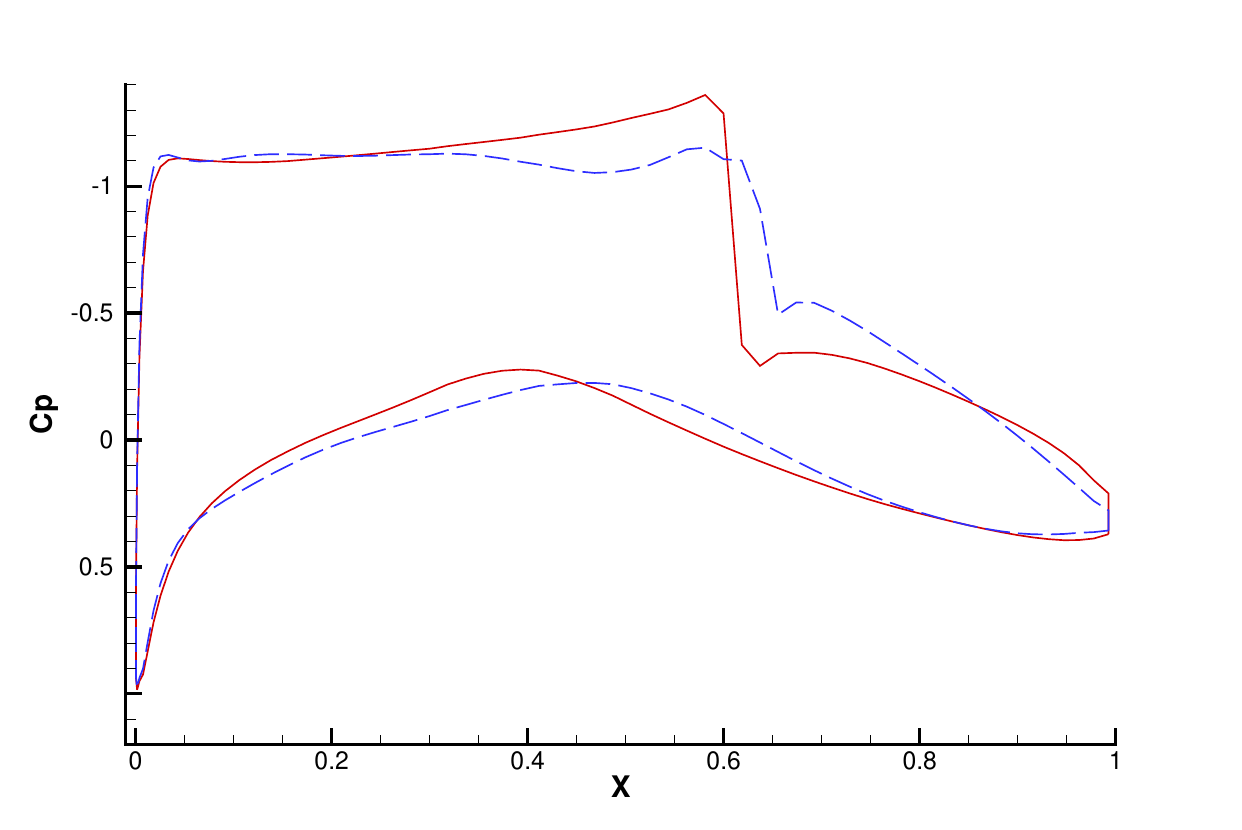}}
 \subfigure[DGp1 on mesh-2]{
 \includegraphics[width=6.0cm]{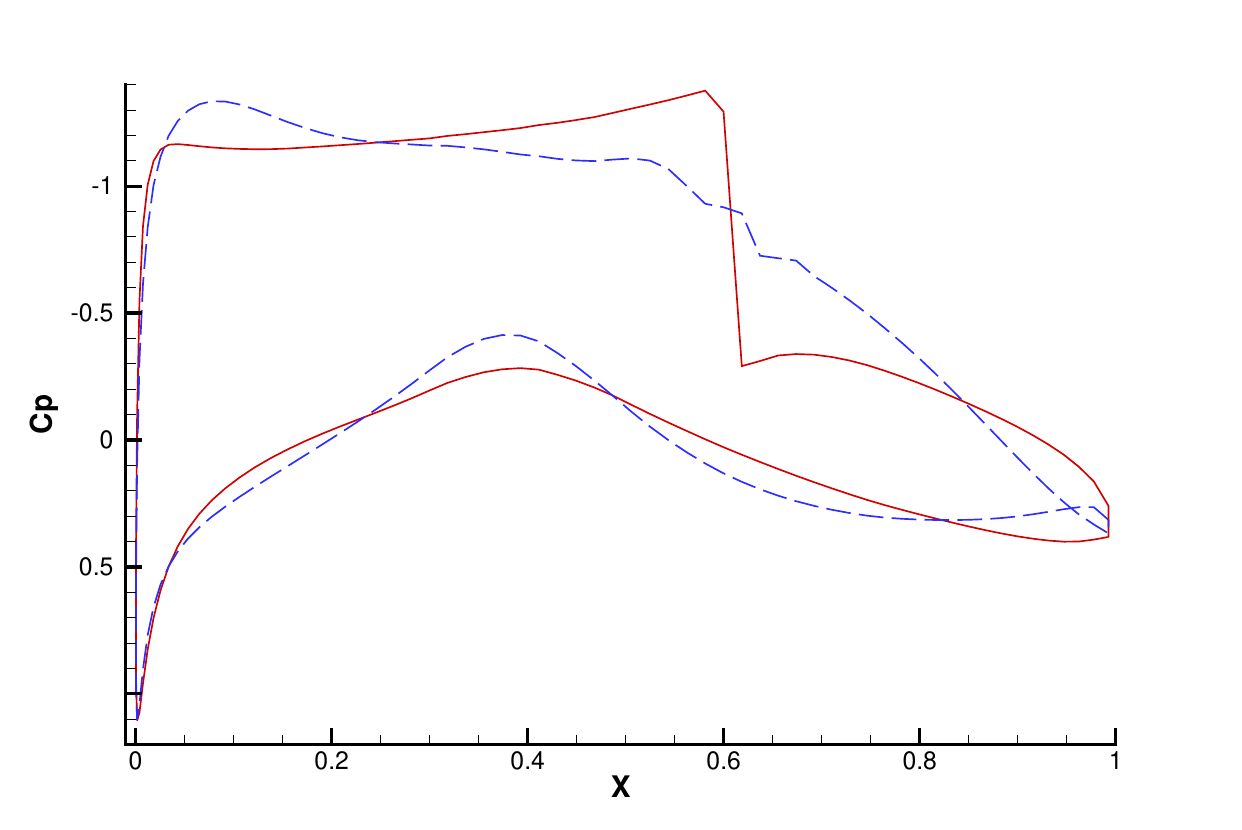}}
 \subfigure[DGp2 on mesh-1]{
 \includegraphics[width=6.0cm]{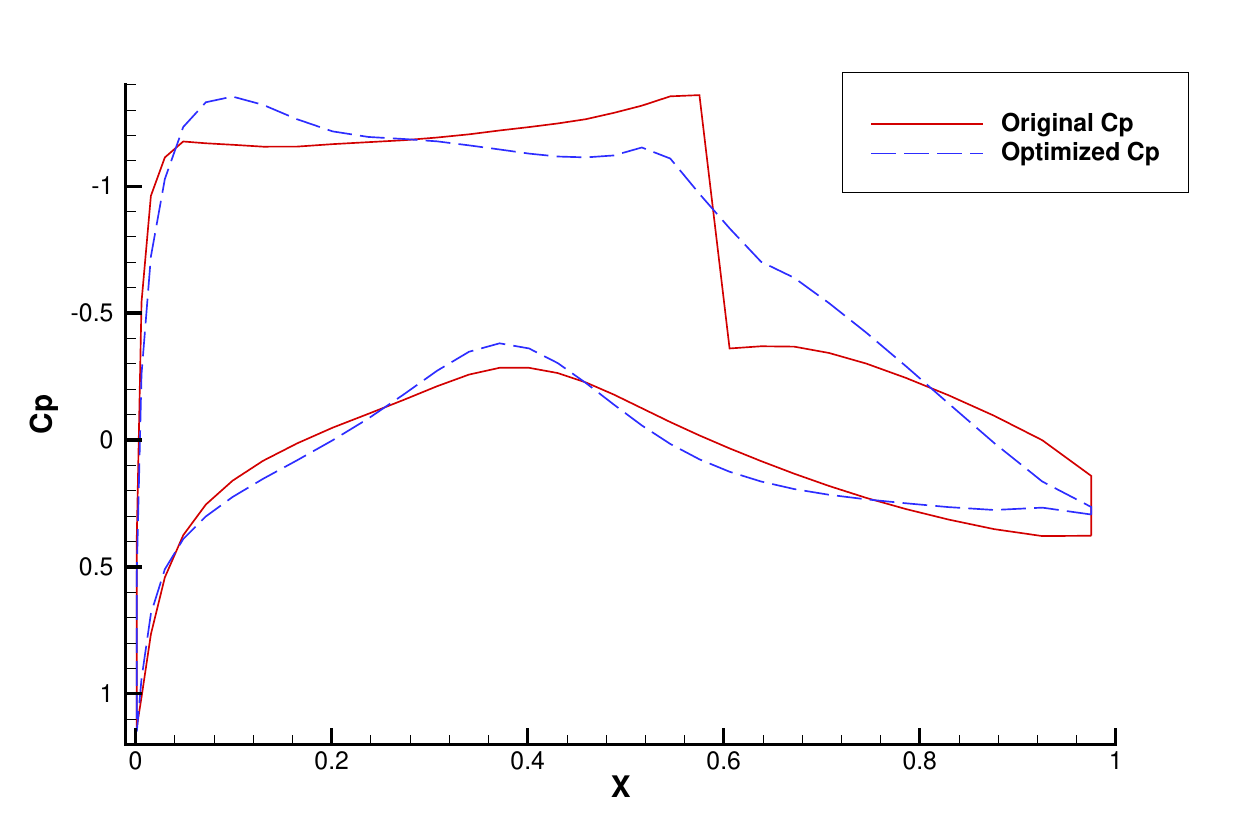}}
 \caption{Original and optimized Cp distribution of the RAE-2822 airfoil test case.}
 \label{fig:rae2822Cp-opt}
\end{figure}

The optimized results are shown in Fig. \ref{fig:rae2822shape-opt}, \ref{fig:rae2822flow-opt} and \ref{fig:rae2822Cp-opt}, similarly only the optimized results of DGp2 on mesh-1, DGp1 and 2ndFV on mesh-2, 1stFV on mesh-3 are provided. At the same time, the iteration steps of optimization and the serial CPU run time in different cases are presented in Table \ref{tab:rae2822-time} for efficiency comparison, and the data of the aerodynamic performance improvement optimized by different CFD solvers are provided in Table \ref{tab:rae2822-opt}.

\begin{table}[htbp]
\centering \caption{Iteration step and serial CPU time (min) with/without DoFs remapping in different cases of the RAE-2822 airfoil test case}
\begin{center}
\begin{tabular}{c|c|c|c|c|c|c|c|c|c}
\toprule
 \multirow{2}*  & \multicolumn{3}{c|}{Iteration step} & \multicolumn{3}{c|}{\small{CPU time (no remapping)}} &\multicolumn{3}{c}{\small{CPU time (remapping)}}\\ \cline{2-10}
  &\small{mesh-1} &\small{mesh-2} &\small{mesh-3} &\small{mesh-1} &\small{mesh-2} &\small{mesh-3} &\small{mesh-1} &\small{mesh-2} &\small{mesh-3} \\ \hline
  1stFV &12 &13 &19 &1.18 &3.78 &35.73 &1.10 &2.43 &14.27 \\ \hline
  2ndFV &18 &33 &/ &3.82 &15.07 &/ &1.25 &10.30 &/\\ \hline
  DGp1 &16 &32 &/ &20.58 &120.75 &/ &10.53 &58.35 &/ \\ \hline
  DGp2 &22 &/ &/ &43.23 &/ &/ &16.93 &/ &/ \\
\bottomrule
\end{tabular}
\end{center}
\label{tab:rae2822-time}
\end{table}

\begin{table}[htbp]
\centering \caption{Performance improvement of various CFD solvers in the RAE-2822 airfoil test case}
\begin{center}
\begin{tabular}{c c c c c c c c c c}
\toprule
     &\small{$C_{d0}$} &\small{$C_{d1}$} &\small{$\Delta C_{d}$} &\small{$C_{l0}$} &\small{$C_{l1}$} &\small{$\Delta C_{l}$} &\small{$A_{0}$} &\small{$A_{1}$} &\small{$\Delta A$} \\
\midrule
  \small{Reference\cite{wang2019adjoint}} &\small{1.56e-2} &/ &/ &\small{8.95e-1} &/ &/ &\small{7.79e-2} &/ &/ \\
  \small{1stFV on mesh-3} &\small{2.46e-2} &\small{2.20e-2} &\small{-10.6$\%$} &\small{7.46e-1} &\small{7.46e-1} &\small{+0$\%$} &\small{7.79e-2} &\small{7.79e-2} &\small{+0$\%$} \\
  \small{2ndFV on mesh-2} &\small{1.51e-2} &\small{1.06e-2} &\small{-29.8$\%$} &\small{8.72e-1} &\small{8.72e-1} &\small{+0$\%$} &\small{7.81e-2} &\small{7.81e-2} &\small{+0$\%$} \\
  \small{DGp1 on mesh-2} &\small{1.37e-2} &\small{6.42e-3} &\small{-53.1$\%$} &\small{8.84e-1} &\small{8.84e-1} &\small{+0$\%$} &\small{7.81e-2} &\small{7.81e-2} &\small{+0$\%$} \\
  \small{DGp2 on mesh-1} &\small{1.37e-2} &\small{7.26e-3} &\small{-47.0$\%$} &\small{8.67e-1} &\small{8.67e-1} &\small{+0$\%$} &\small{7.82e-2} &\small{7.84e-2} &\small{+0$\%$} \\
\bottomrule
\end{tabular}
\end{center}
\label{tab:rae2822-opt}
\end{table}

Very similar conclusions can be made from Table \ref{tab:rae2822-time}-\ref{tab:rae2822-opt} and Fig. \ref{fig:rae2822shape-opt}-\ref{fig:rae2822Cp-opt} that (a) under coequal computational costs, the optimized shape, flow distribution, and the Cp distributions produced by DGp2 on mesh-1 and DGp1 on mesh-2 are similar, and are significantly superior to that produced by FVMs. (b) DGMs can achieve around $50\%$ drag reduction, while the 1stFV is only $10\%$ and the 2ndFV is around $30\%$. (c) the the CPU time cost by DGp2 on mesh-1 and FVMs is similar, and much less than the time cost by DGp1 on mesh-2; (d) the application of DoFs remapping techniques in most cases can reduce the CPU time of the optimization by around $35\%\sim70\%$.

\subsection{Onera M6 wing drag minimization}
\qquad We then consider the 3D Onera M6 wing drag minimization, the target is to reduce the drag $C_d$ without decreasing the lift $C_l$ and the volume $V$ of the wing as well. The free stream is set to $\text{Ma}=0.84, \text{AOA}=3.06^o, \text{AOS}=0^o$. Only the changes of z-coordinates are taken into consideration, the number of the design variables is set to 162 ($N_x\times N_y\times N_z = 9\times 9\times 2$), and their ranges of variation are set no more than $20\%$.
%\begin{figure}[htbp]
%  \centering
%  \includegraphics[width=5.0cm]{figure/m6-FFD-Box.pdf}
%  \caption{The FFD box of the Onera M6 wing}
%  \label{fig:m6-FFD}
%\end{figure}

\begin{figure}[htbp]
  \centering
  \subfigure[mesh-1]{
  \includegraphics[width=6.0cm]{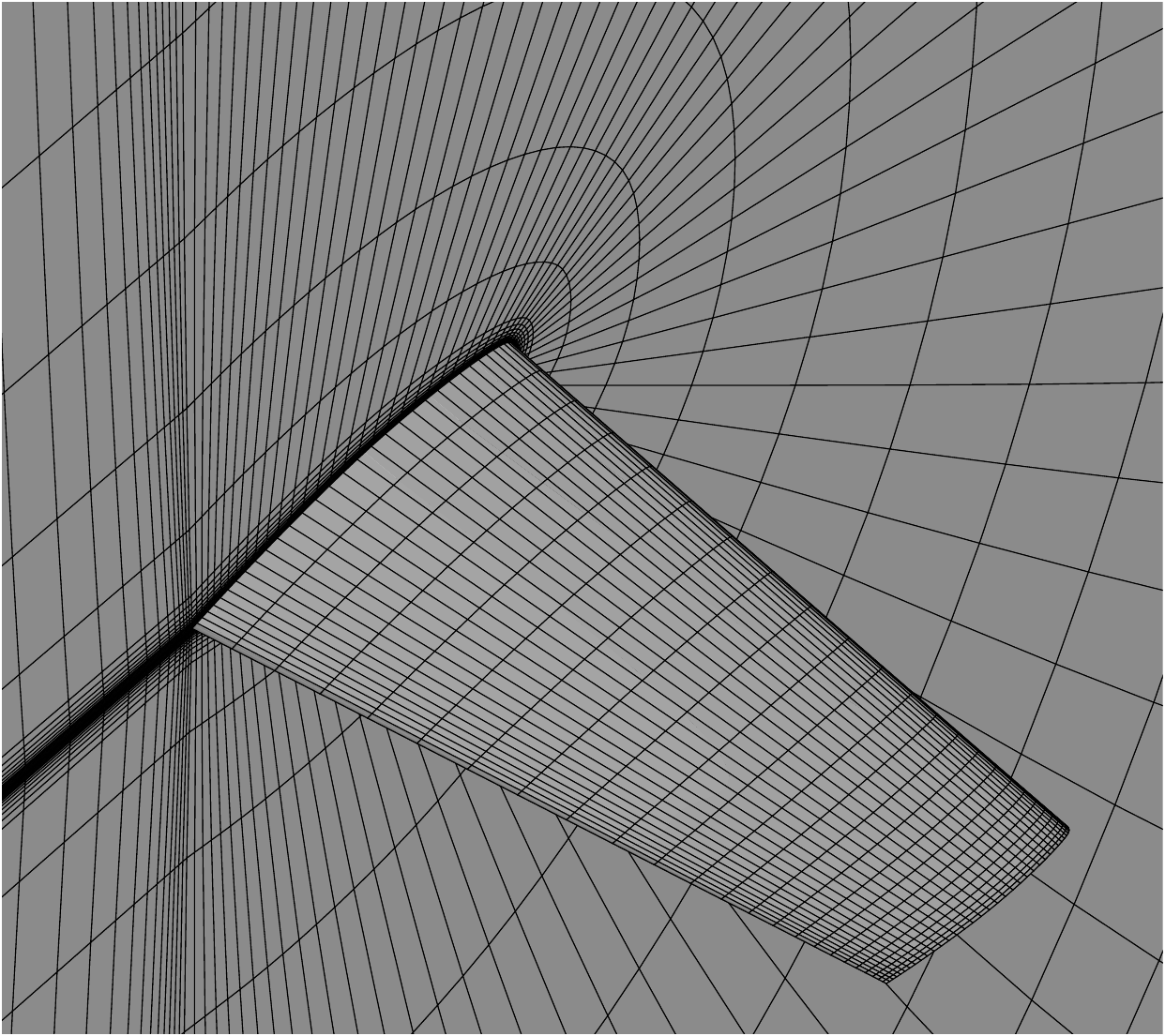}}
  \subfigure[mesh-2]{
  \includegraphics[width=6.0cm]{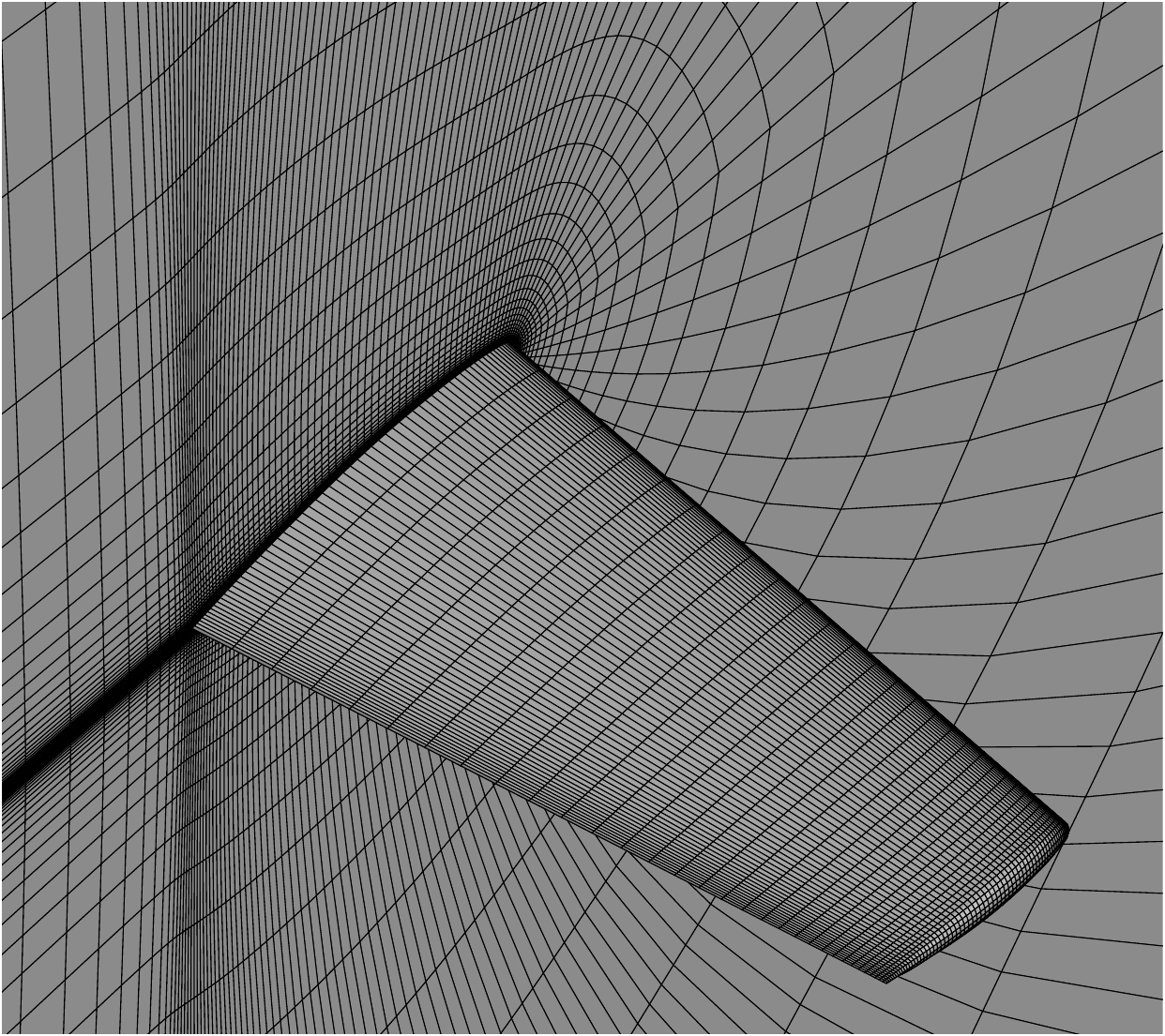}}
  \caption{The surface grid of the Onera M6 wing}
  \label{fig:m6mesh}
\end{figure}

\begin{figure}[htbp]
  \centering
  \subfigure[$95\%$ span]{
  \includegraphics[width=6.0cm]{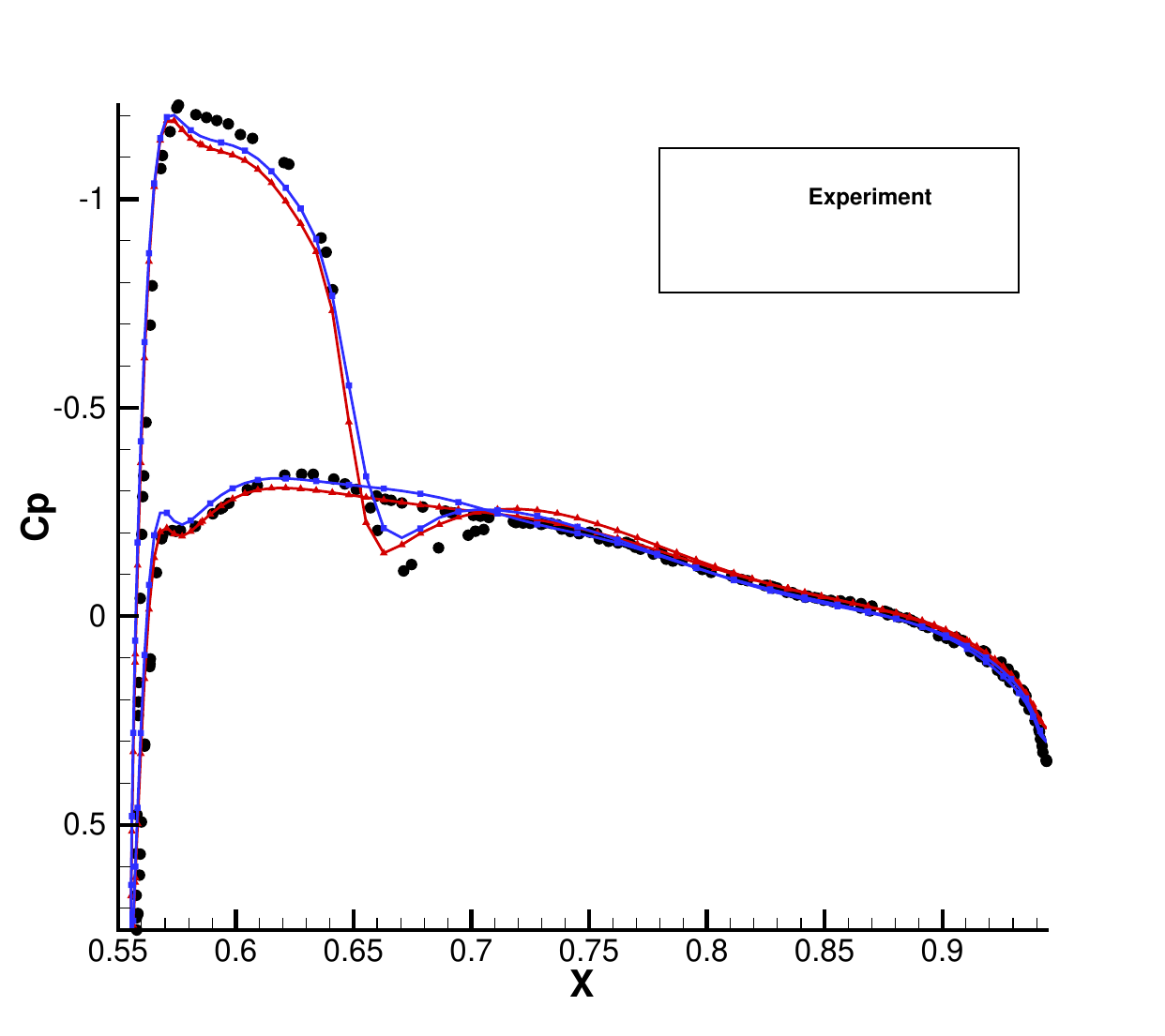}}
  \subfigure[$90\%$ span]{
  \includegraphics[width=6.0cm]{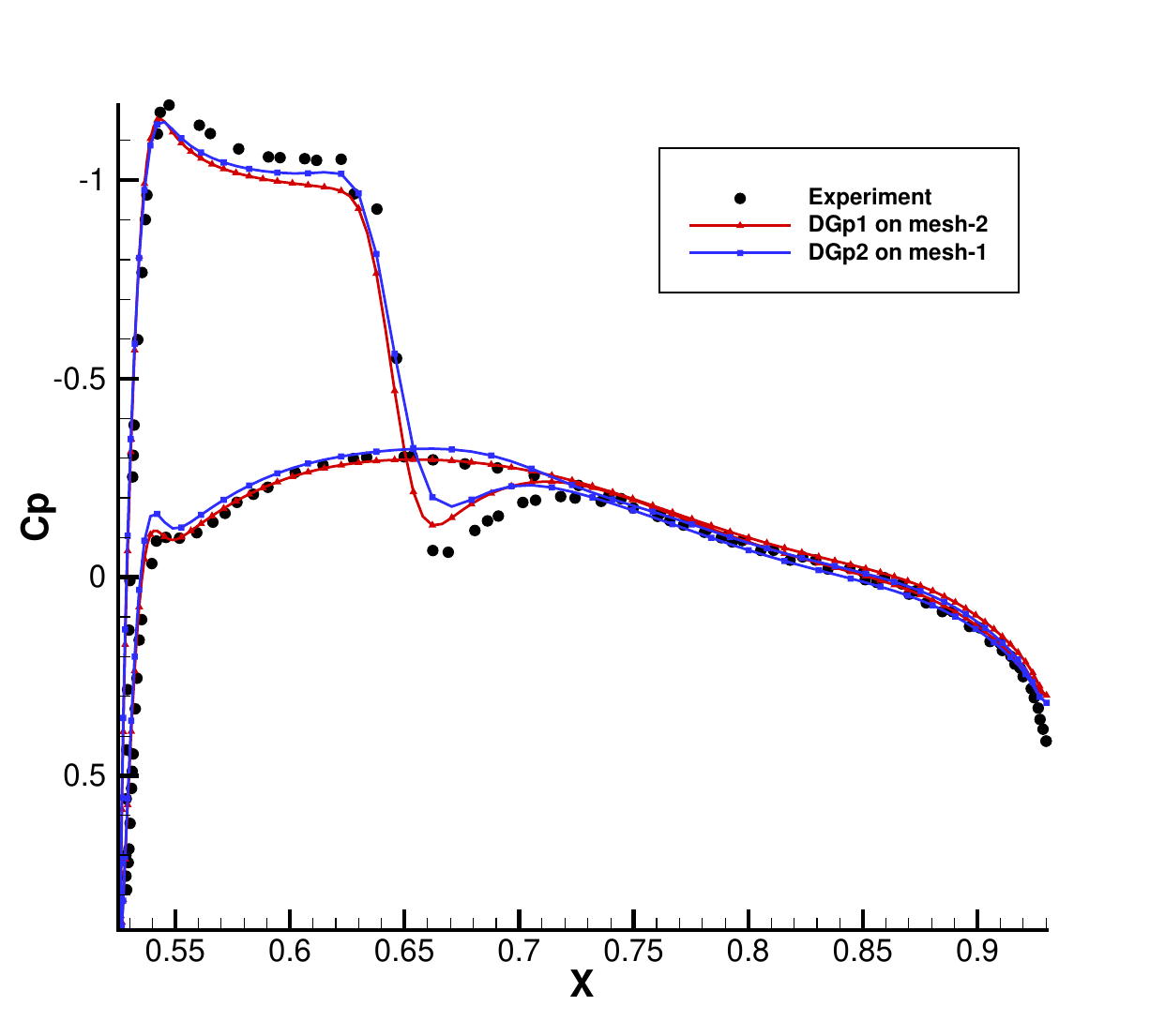}}
  \caption{Comparison of M6 $C_p$ distributions between CFD data and experimental data}
  \label{fig:m6Cp-exp}
\end{figure}

\begin{figure}[htbp]
 \centering
 \subfigure[Original density]{
 \includegraphics[width=4.8cm]{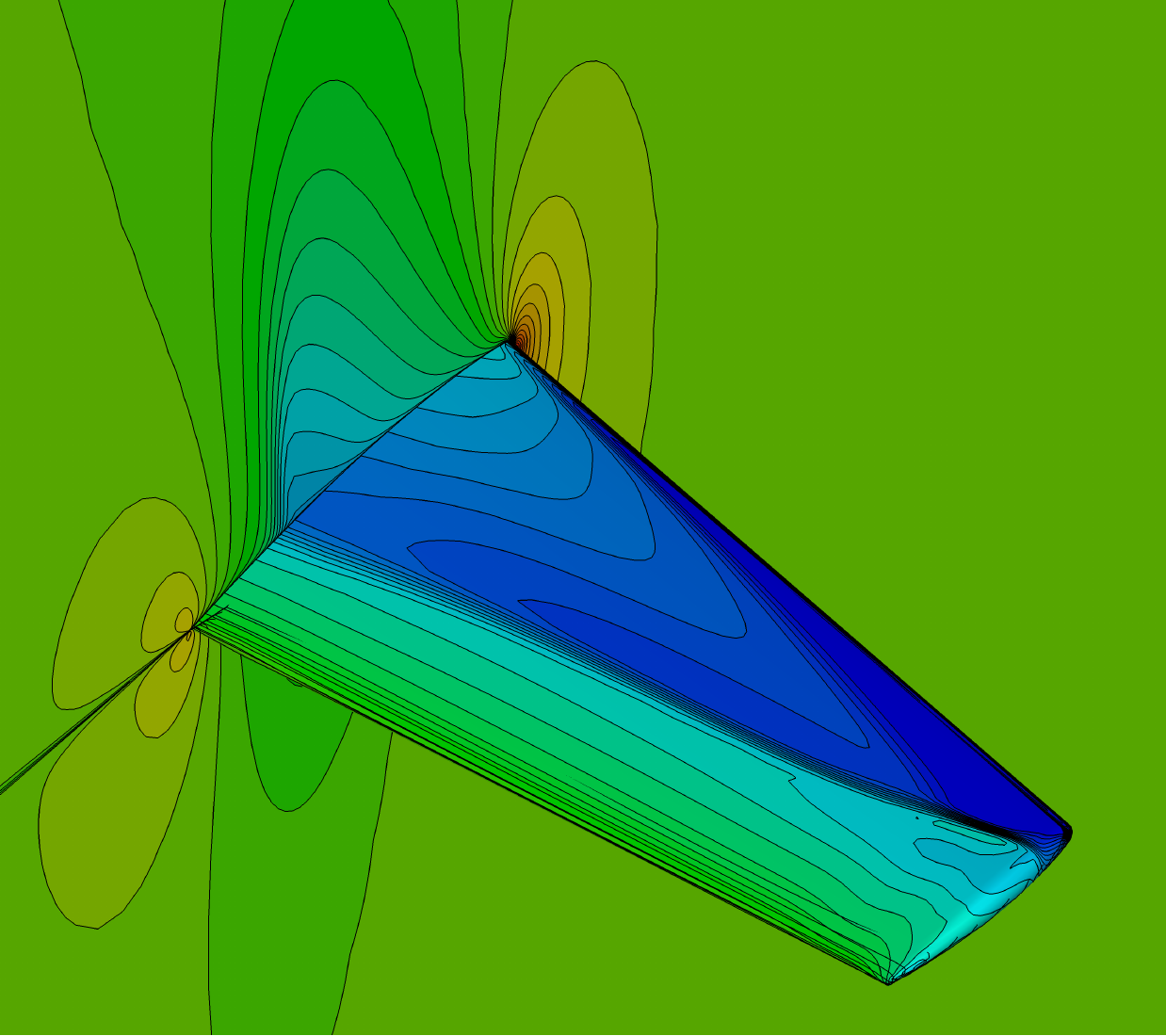}}
 \subfigure[Density after 5 steps]{
 \includegraphics[width=4.8cm]{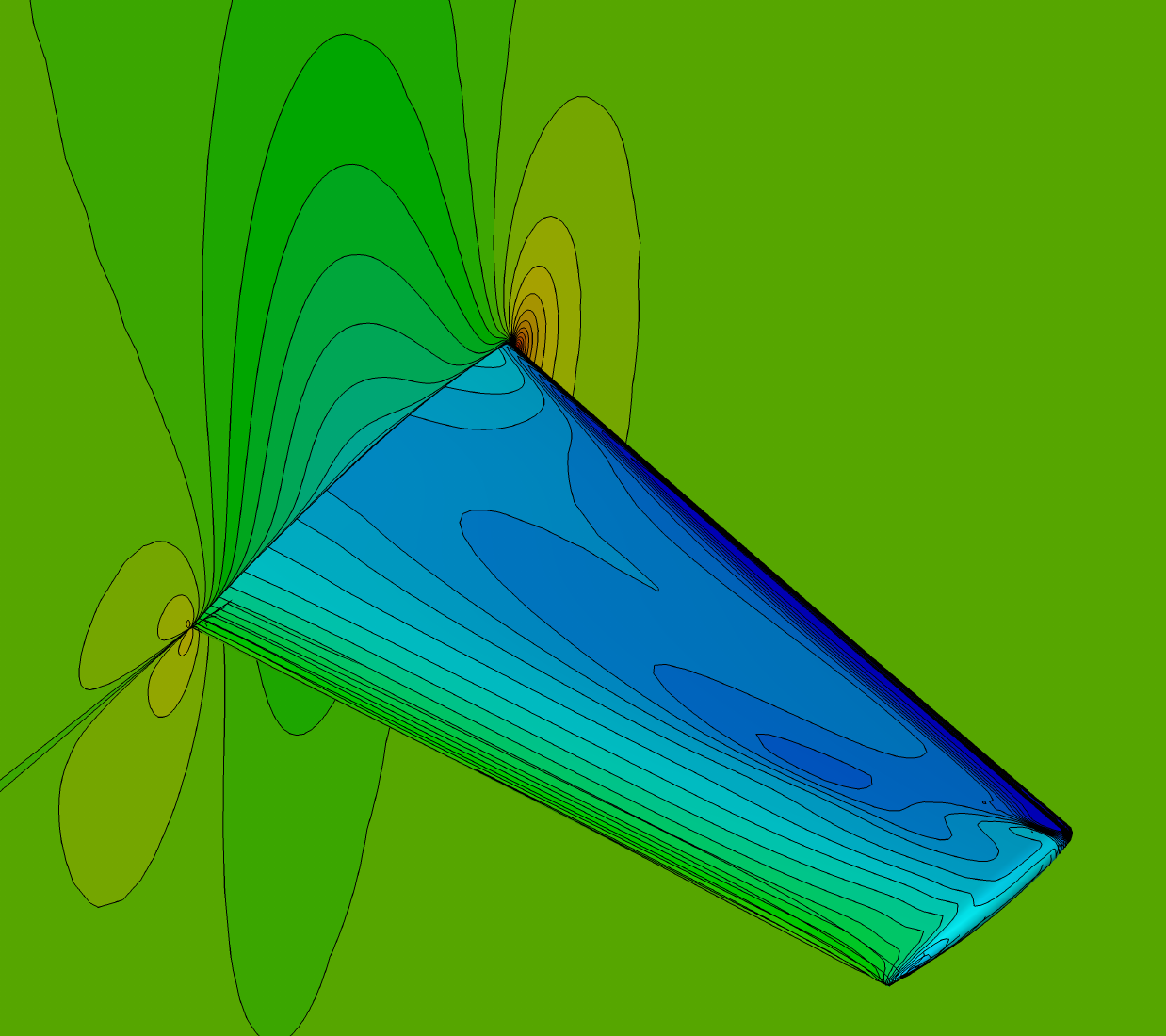}}
 \subfigure[Final density (25 steps)]{
 \includegraphics[width=4.8cm]{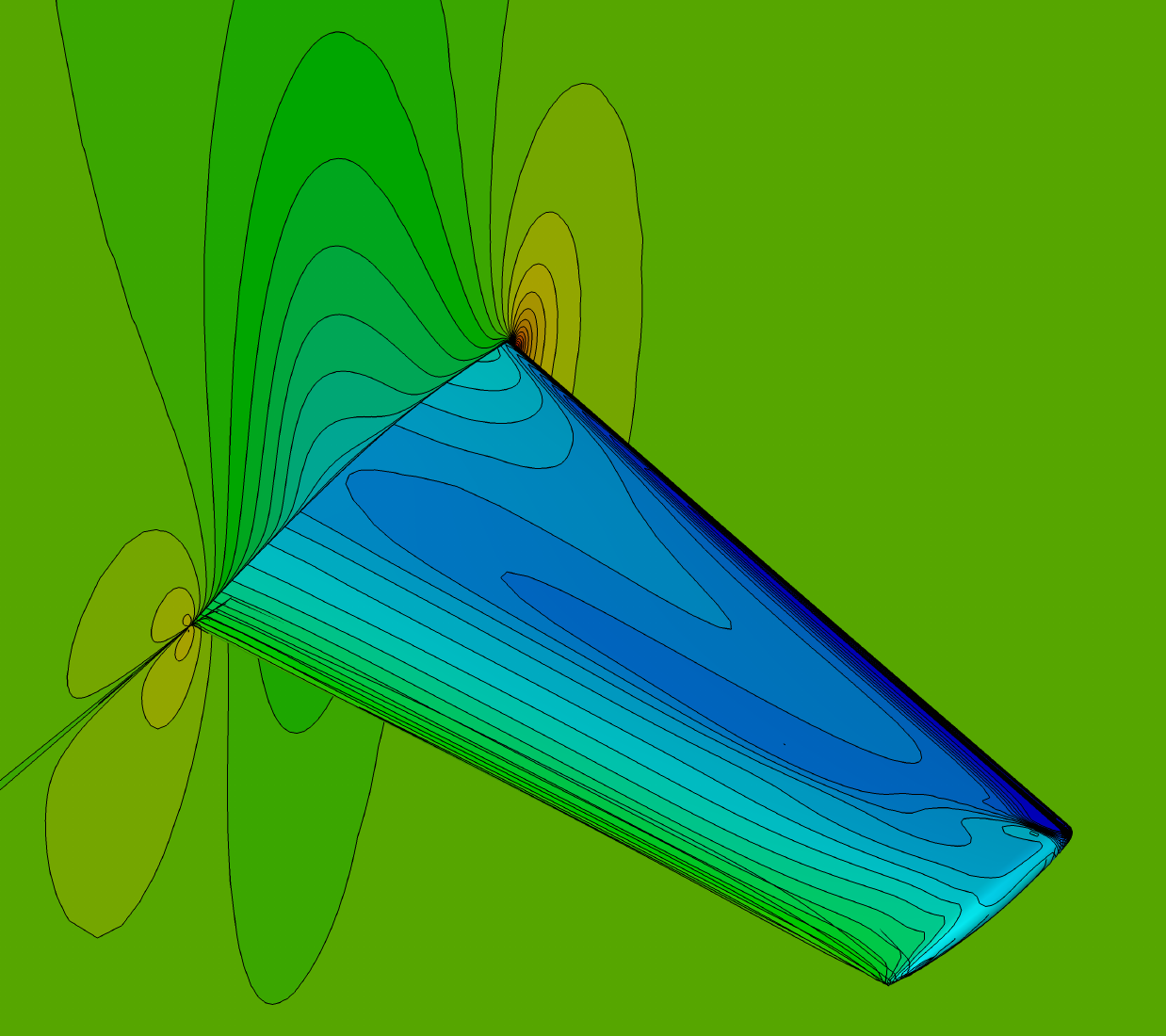}}
 \caption{Optimized density counters by DGp1 on mesh-2 of the Onera M6 wing test case}
 \label{fig:m6-flow1}
\end{figure}

\begin{figure}[htbp]
 \centering
 \subfigure[Original density]{
 \includegraphics[width=4.8cm]{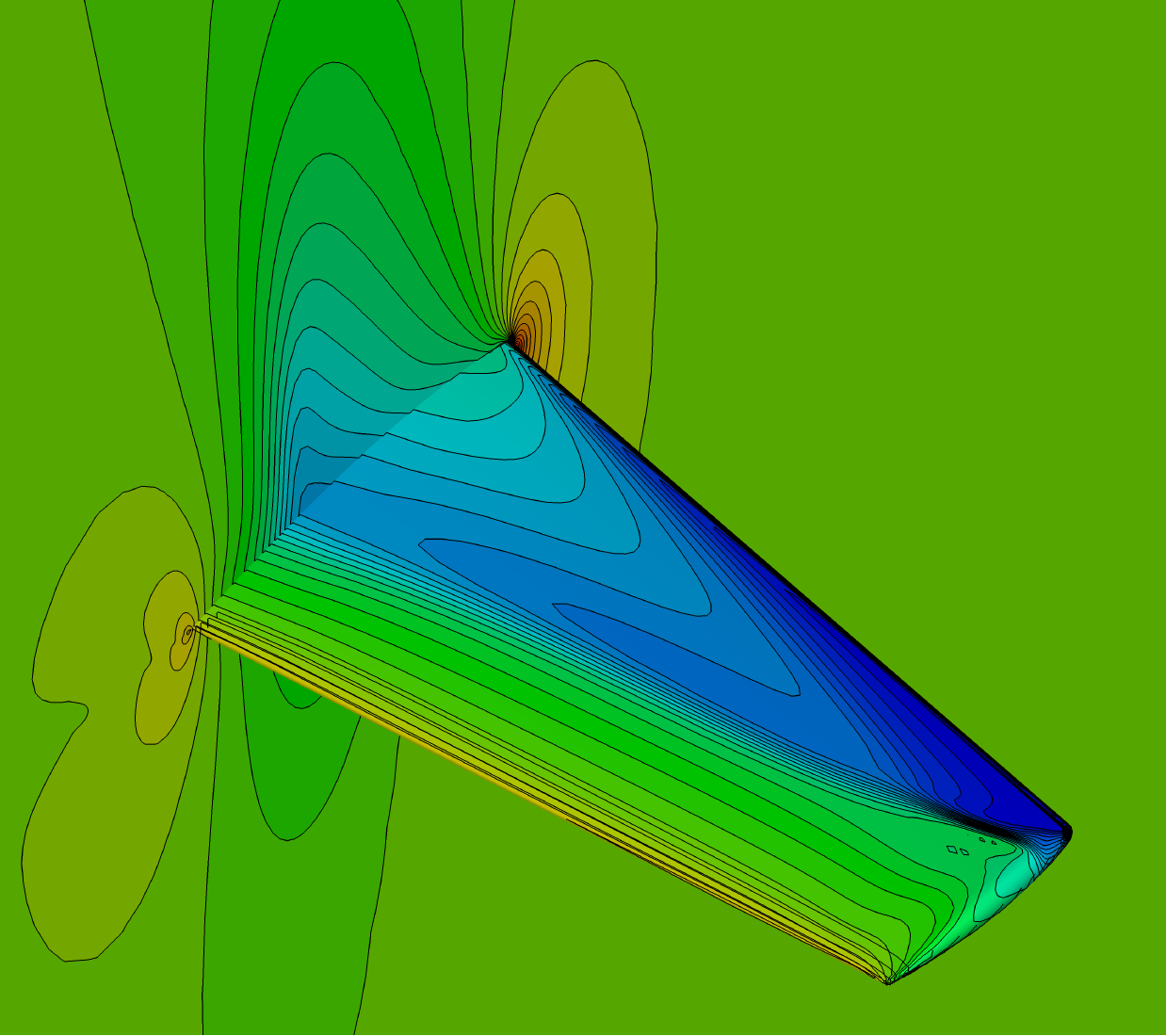}}
 \subfigure[Density after 5 steps]{
 \includegraphics[width=4.8cm]{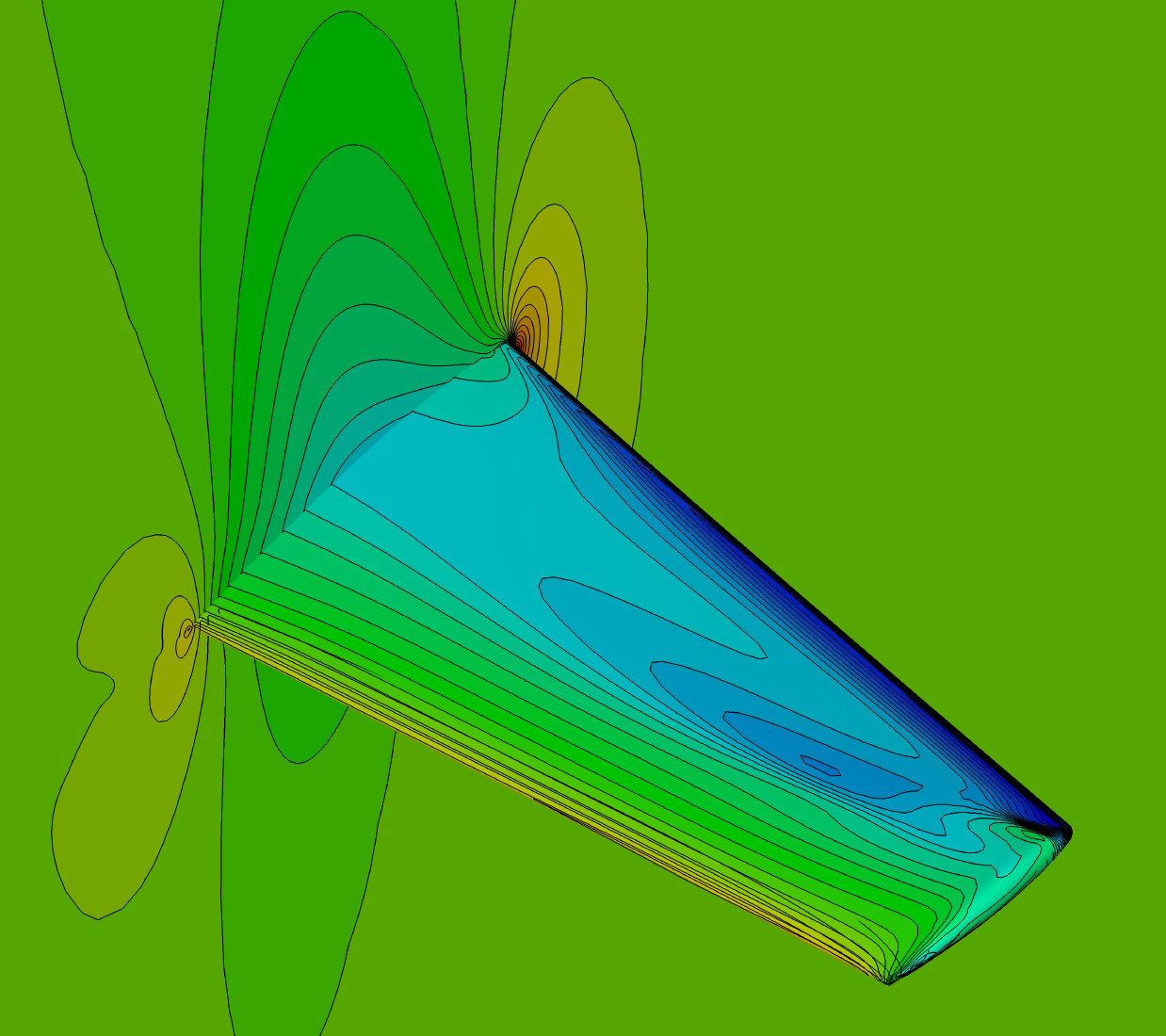}}
 \subfigure[Final density (20 steps)]{
 \includegraphics[width=4.8cm]{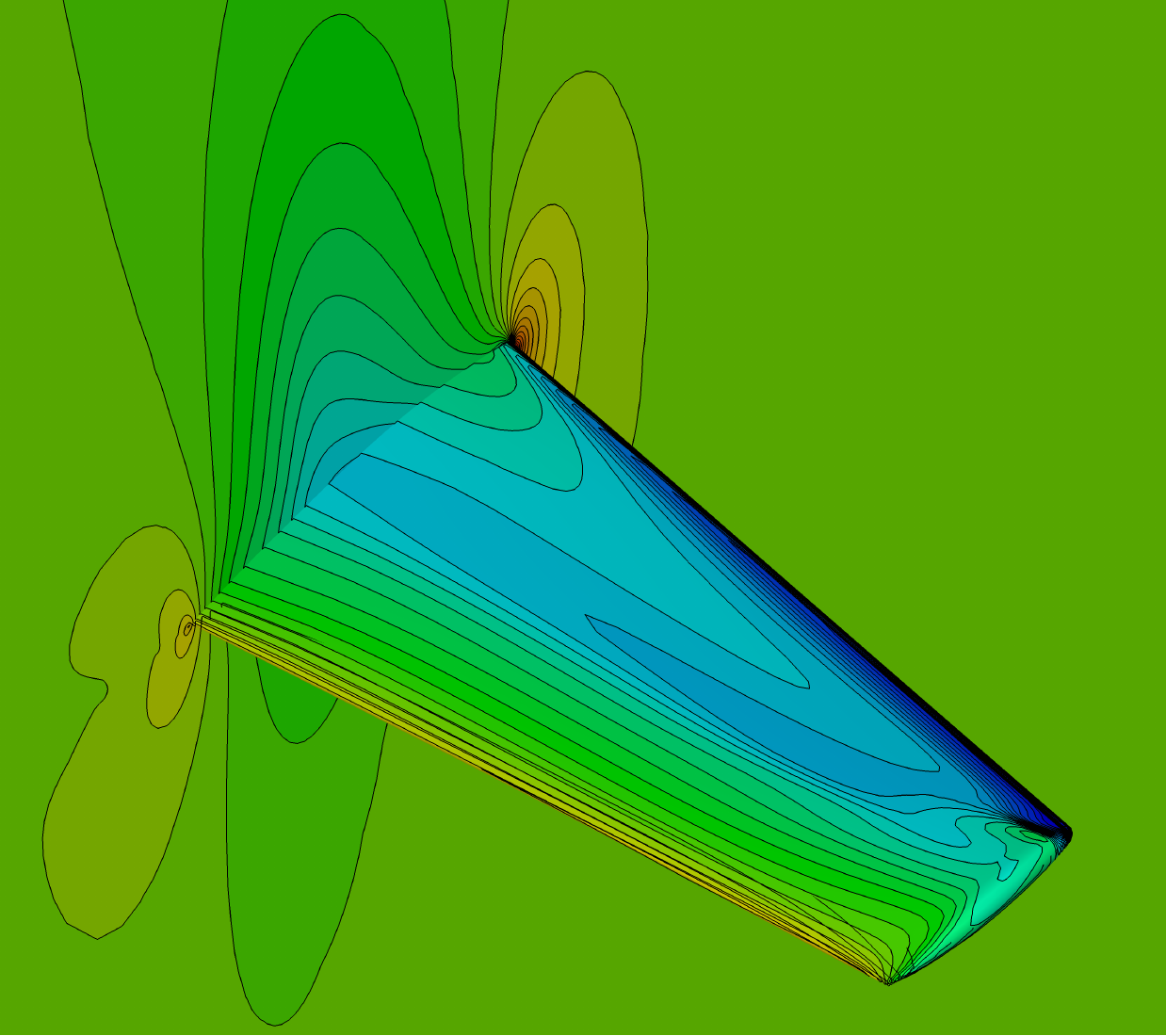}}
 \caption{Optimized density counters by DGp2 on mesh-1 of the Onera M6 wing test case}
 \label{fig:m6-flow2}
\end{figure}

\begin{figure}[htbp]
 \centering
 \subfigure[Original shock]{
 \includegraphics[width=4.8cm]{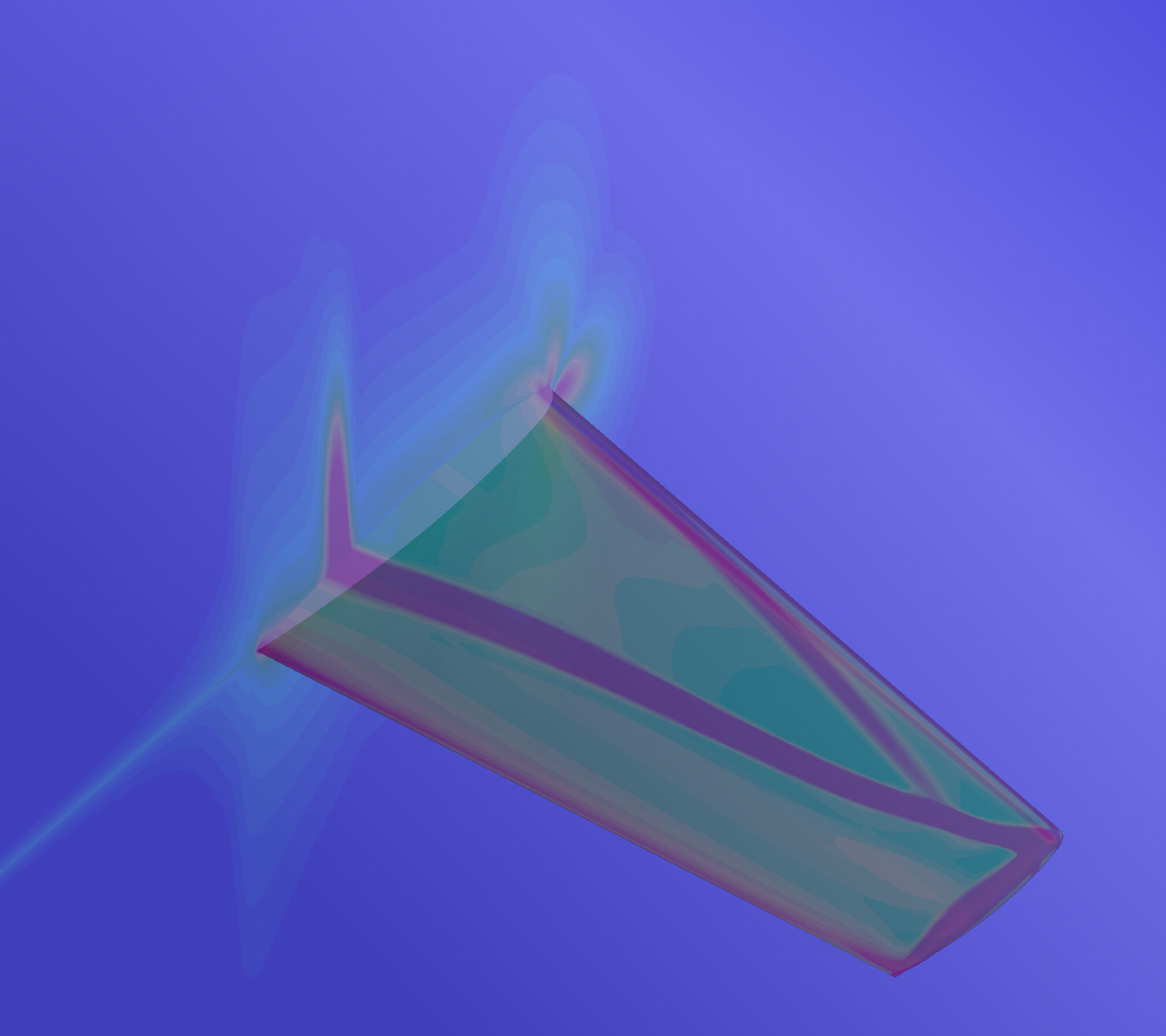}}
 \subfigure[Shock after 5 steps]{
 \includegraphics[width=4.8cm]{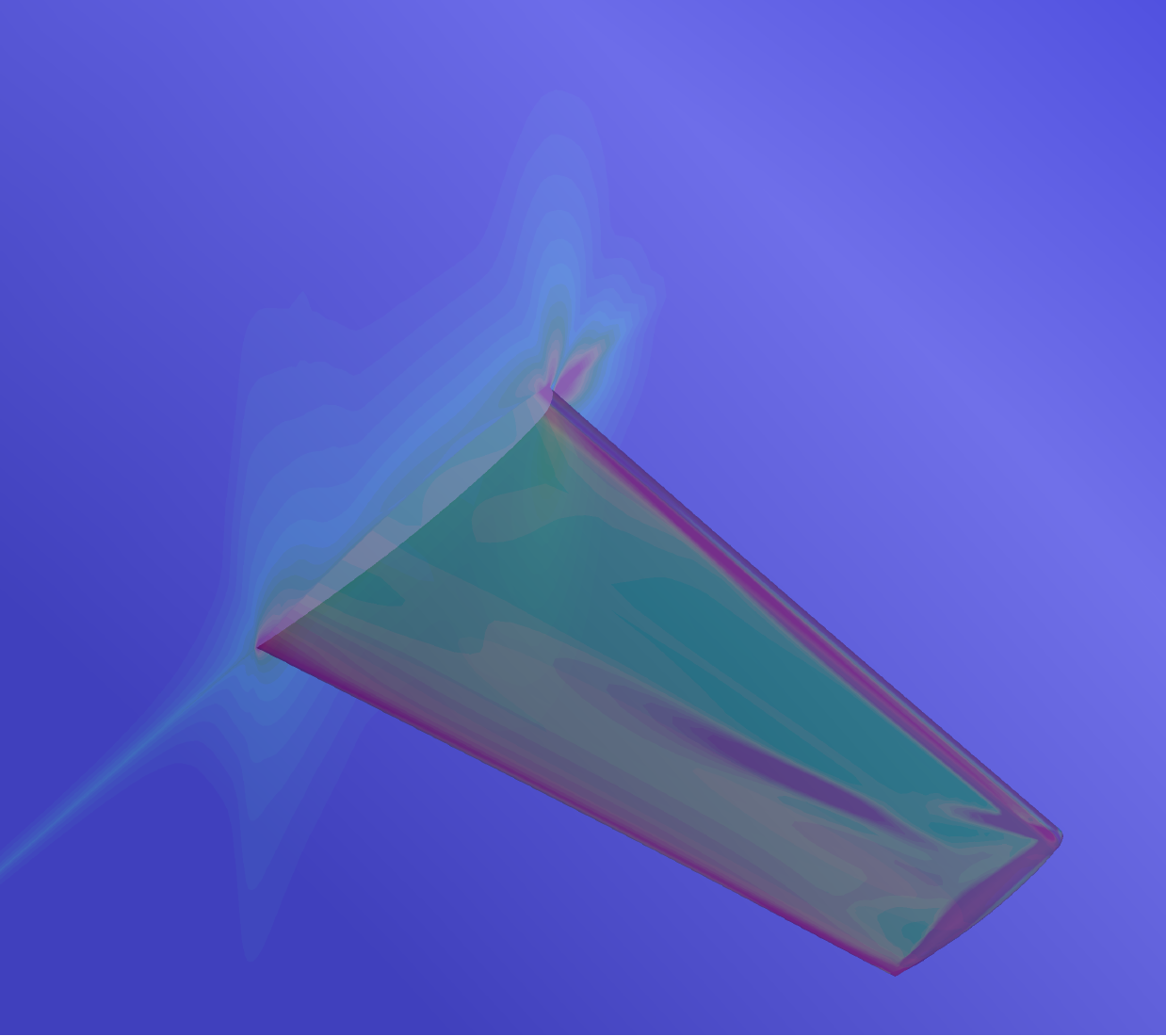}}
 \subfigure[Final shock (25 steps)]{
 \includegraphics[width=4.8cm]{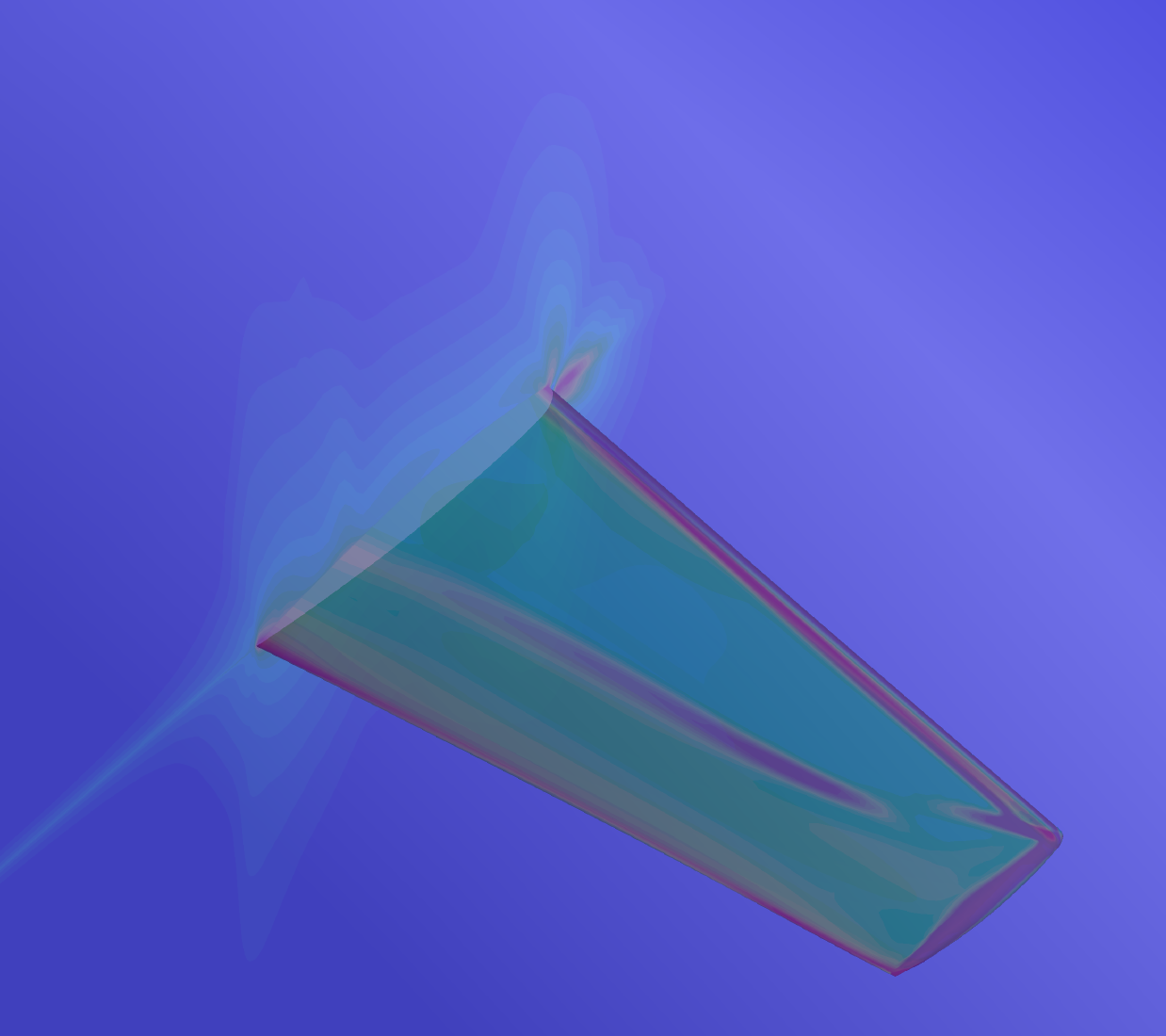}}
 \caption{Shock visualization of DGp1 results during the optimization}
 \label{fig:m6-shock1}
\end{figure}

\begin{figure}[htbp]
 \centering
 \subfigure[Original shock]{
 \includegraphics[width=4.8cm]{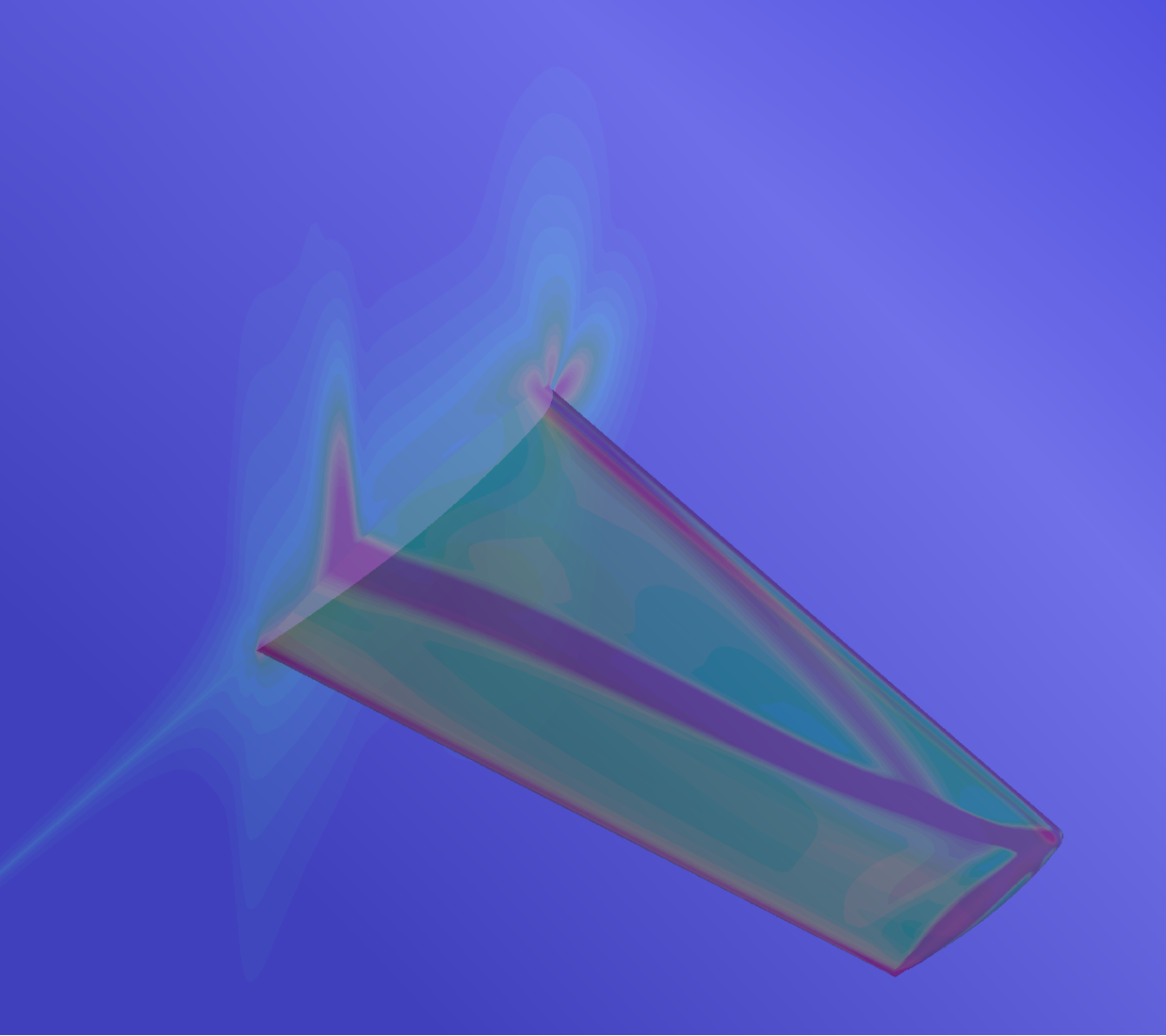}}
 \subfigure[Shock after 5 steps]{
 \includegraphics[width=4.8cm]{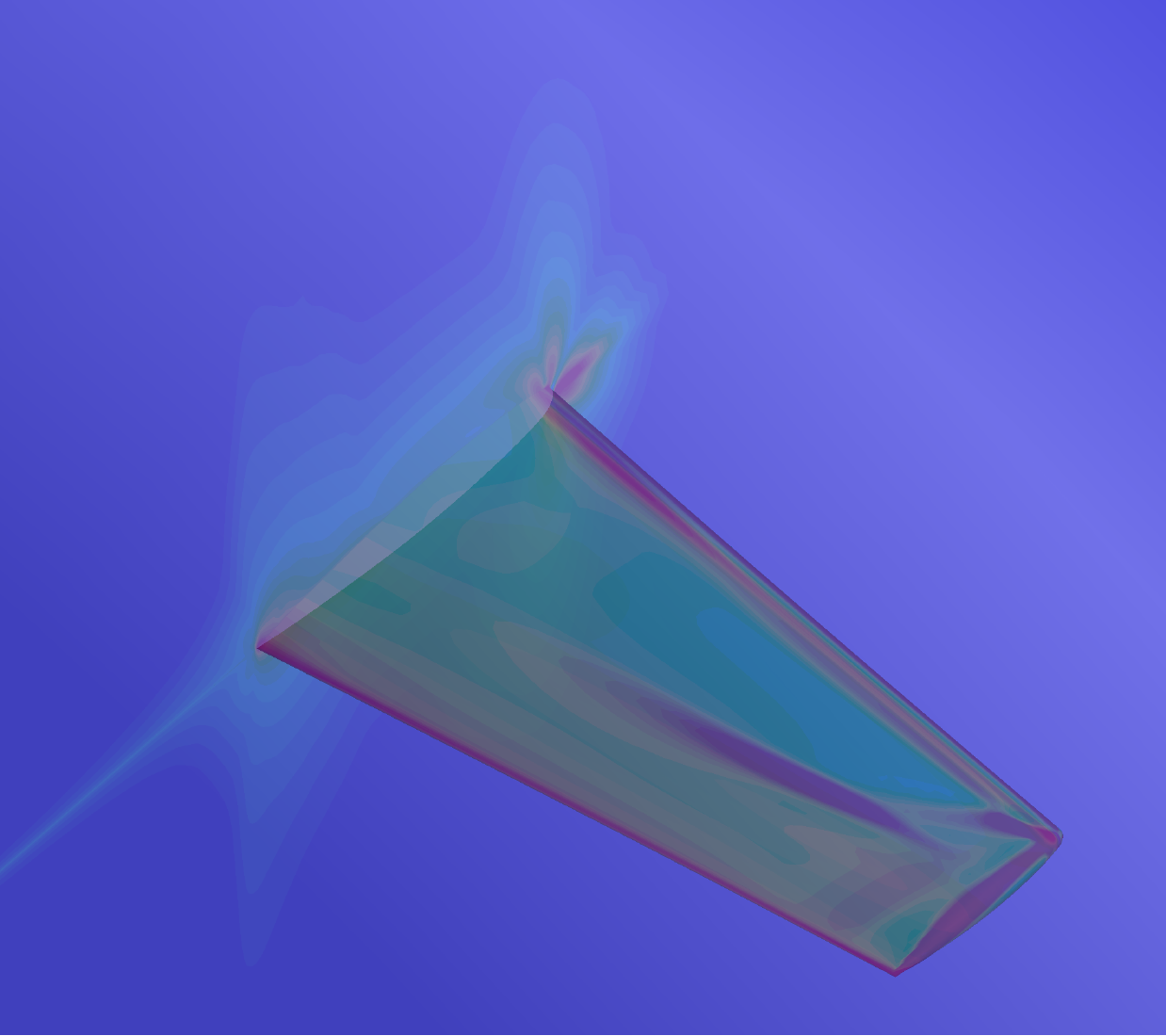}}
 \subfigure[Final shock (20 steps)]{
 \includegraphics[width=4.8cm]{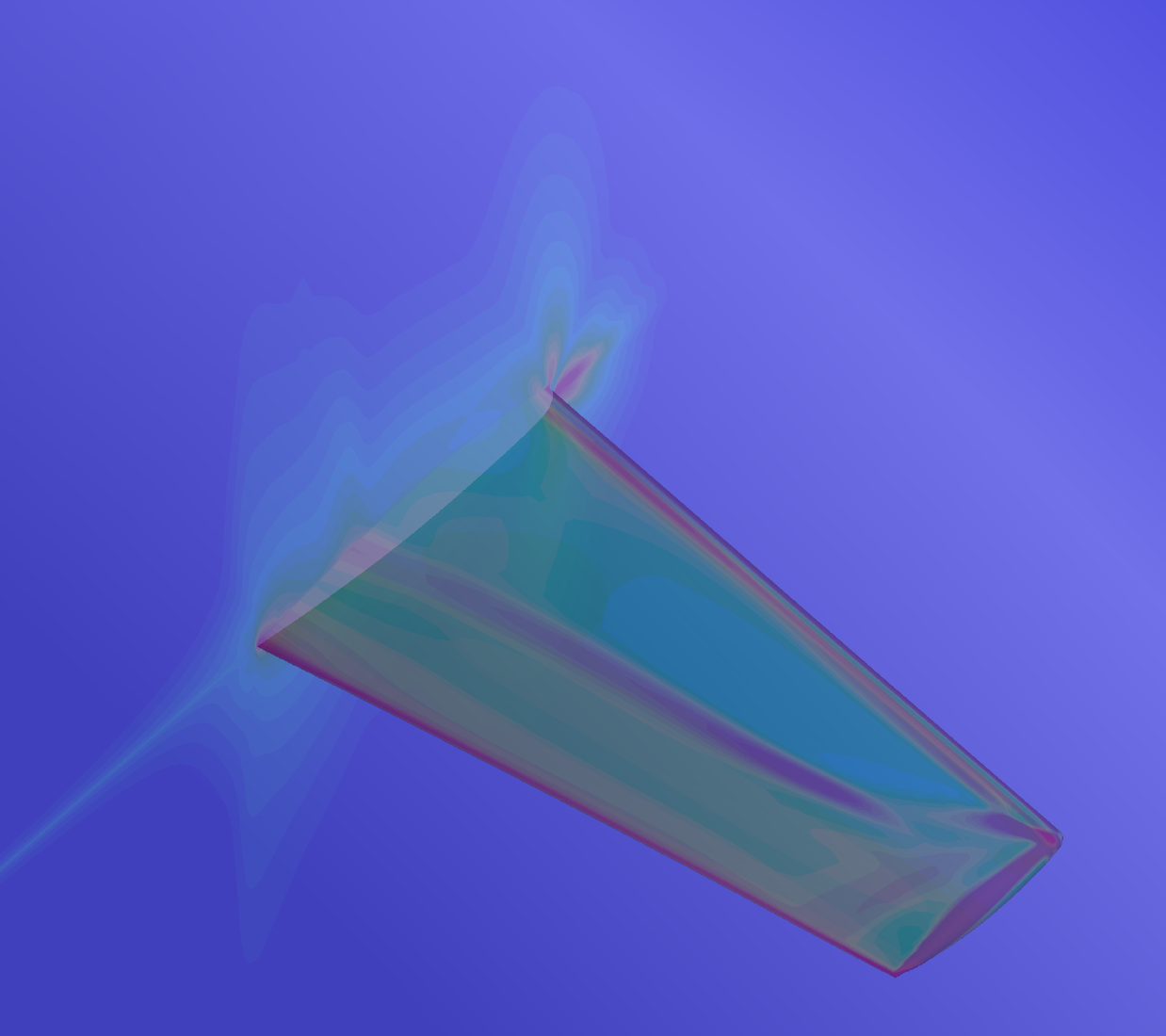}}
 \caption{Shock visualization of DGp2 results during the optimization}
 \label{fig:m6-shock2}
\end{figure}

\begin{figure}[htbp]
 \centering
 \subfigure[$y=20\%\text{ span}$]{
 \includegraphics[width=8.0cm]{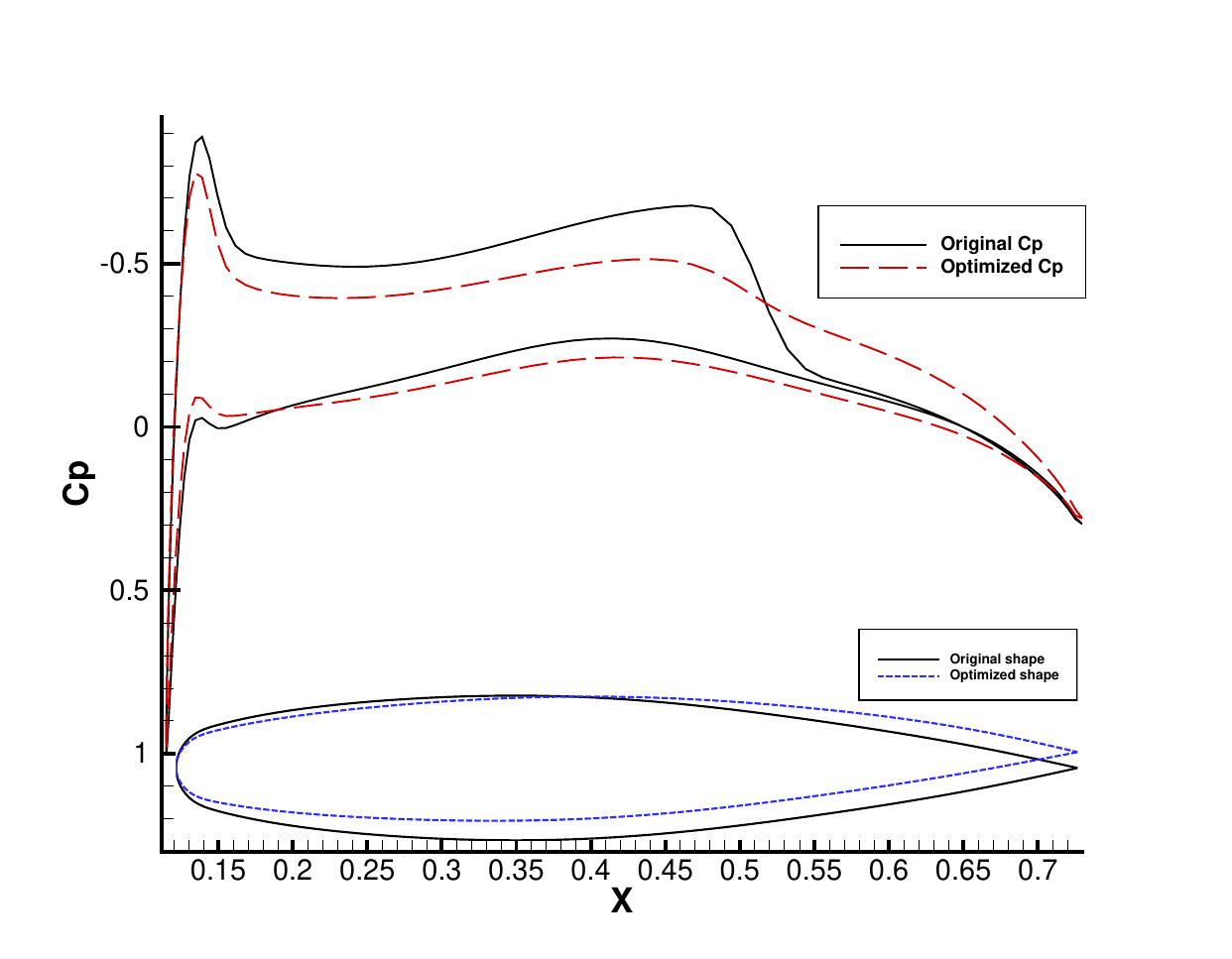}}
 \subfigure[$y=44\%\text{ span}$]{
 \includegraphics[width=8.0cm]{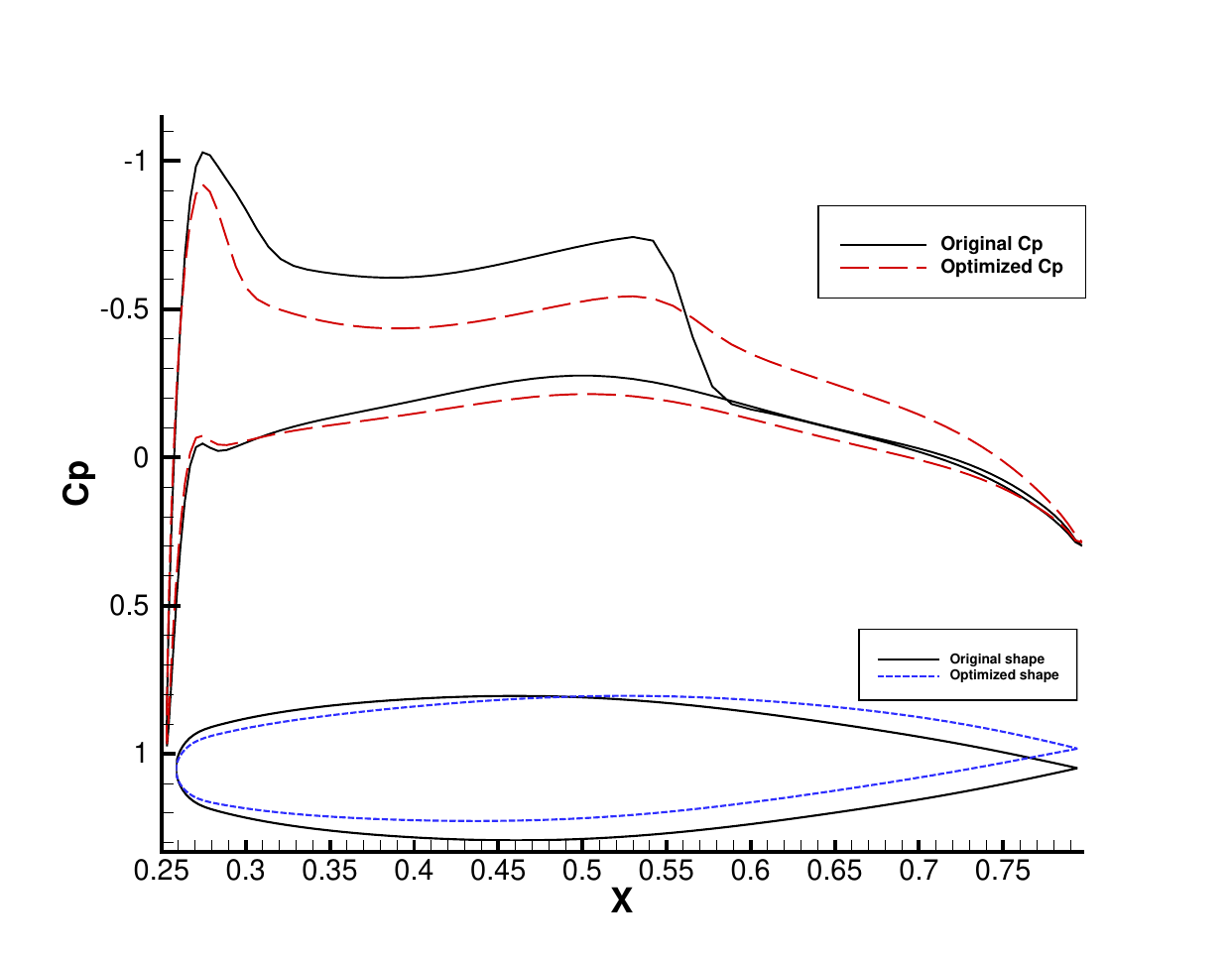}}
 \subfigure[$y=80\%\text{ span}$]{
 \includegraphics[width=8.0cm]{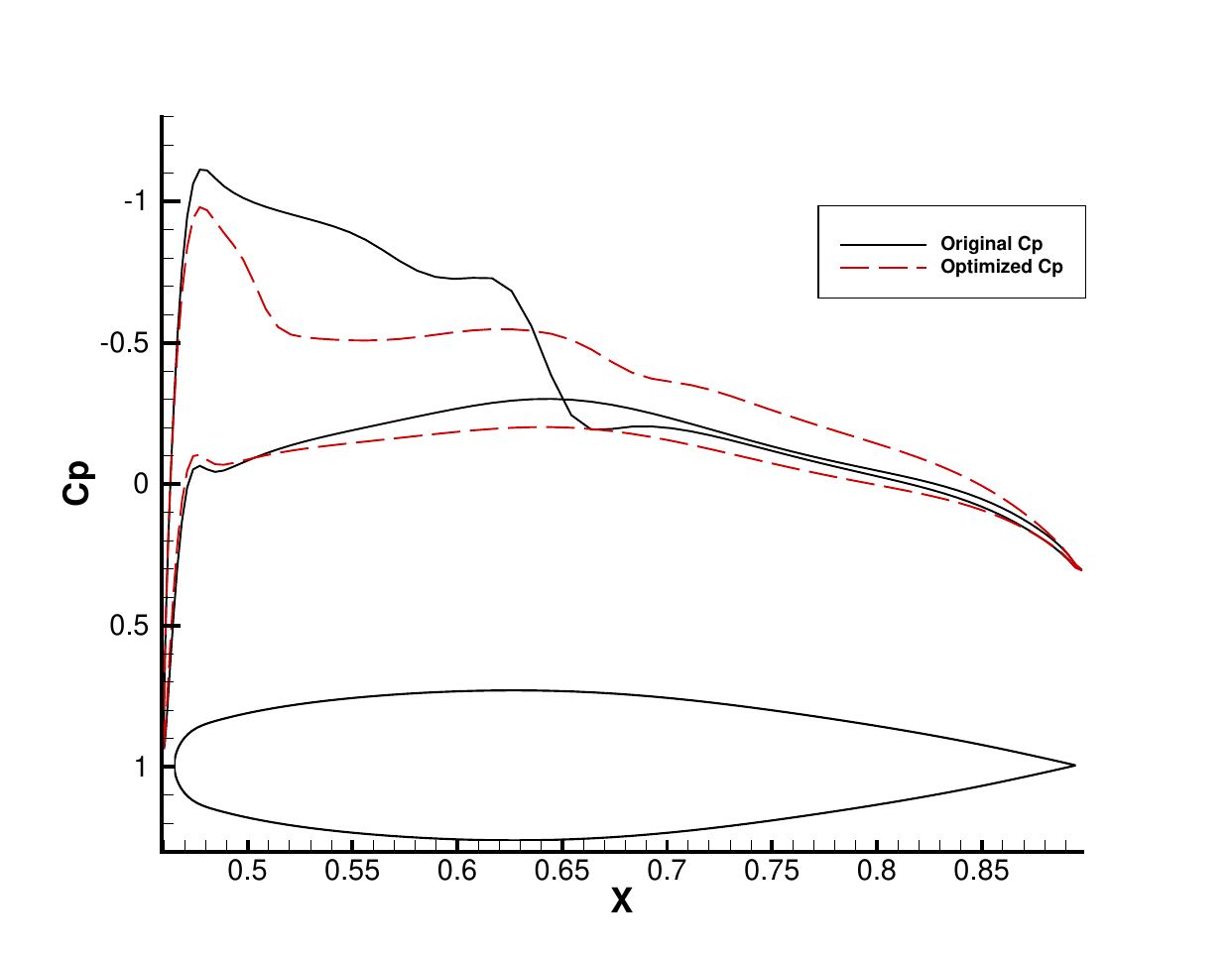}}
 \subfigure[$y=95\%\text{ span}$]{
 \includegraphics[width=8.0cm]{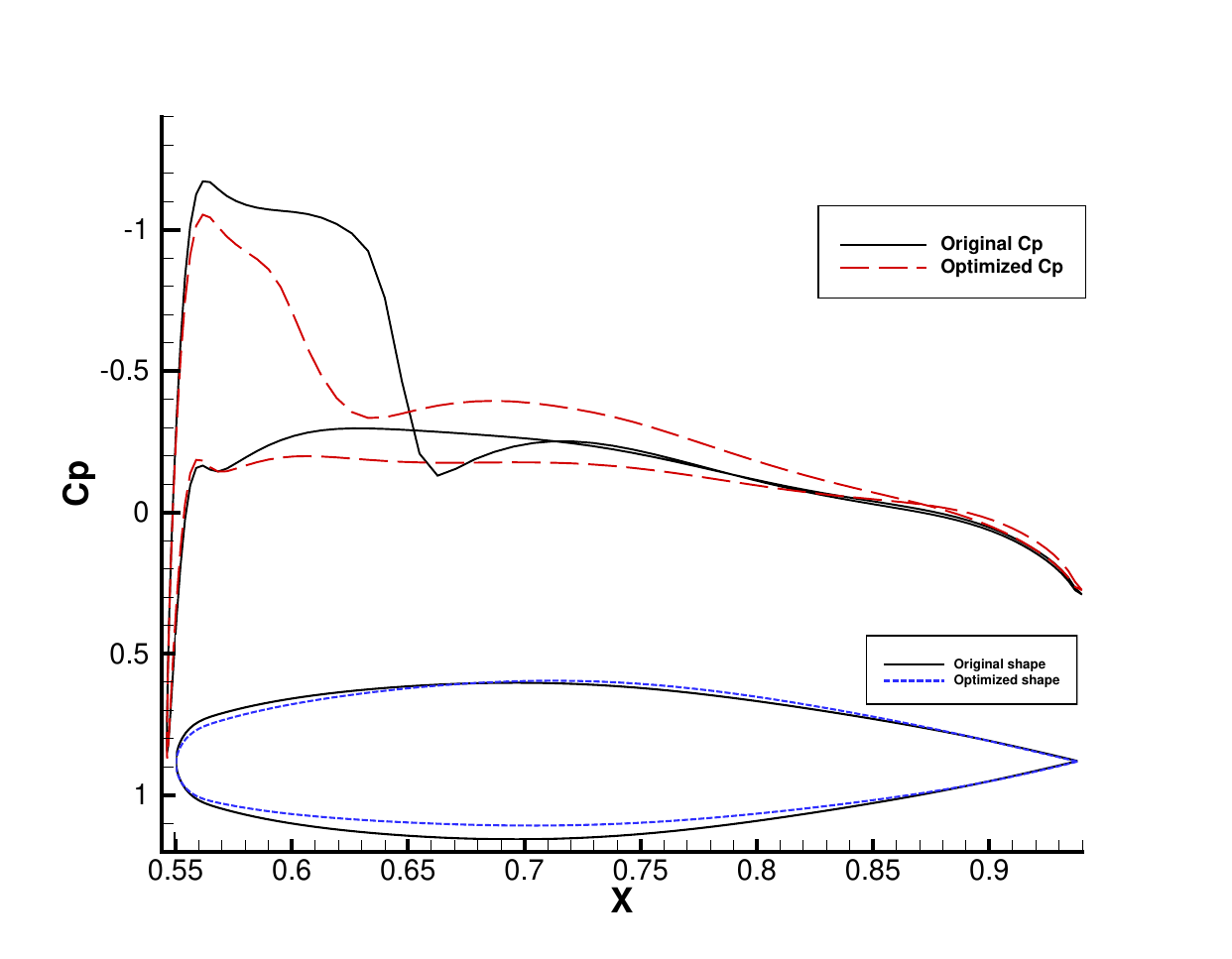}}
 \caption{Optimized shape and Cp distribution of the Onera M6 wing test case by DGp1 on mesh-2}
 \label{fig:m6-Cp1}
\end{figure}

\begin{figure}[htbp]
 \centering
 \subfigure[$y=20\%\text{ span}$]{
 \includegraphics[width=8.0cm]{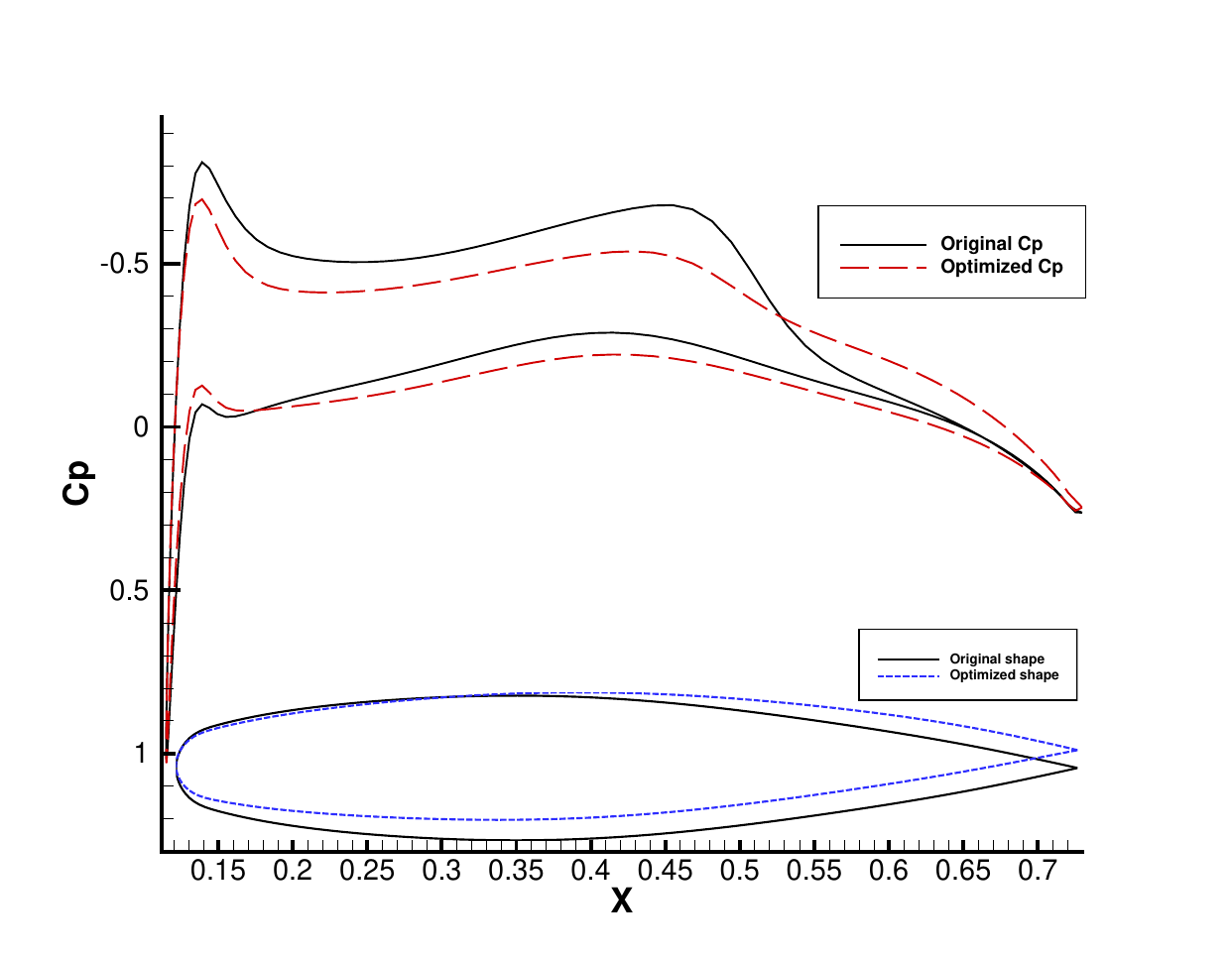}}
 \subfigure[$y=44\%\text{ span}$]{
 \includegraphics[width=8.0cm]{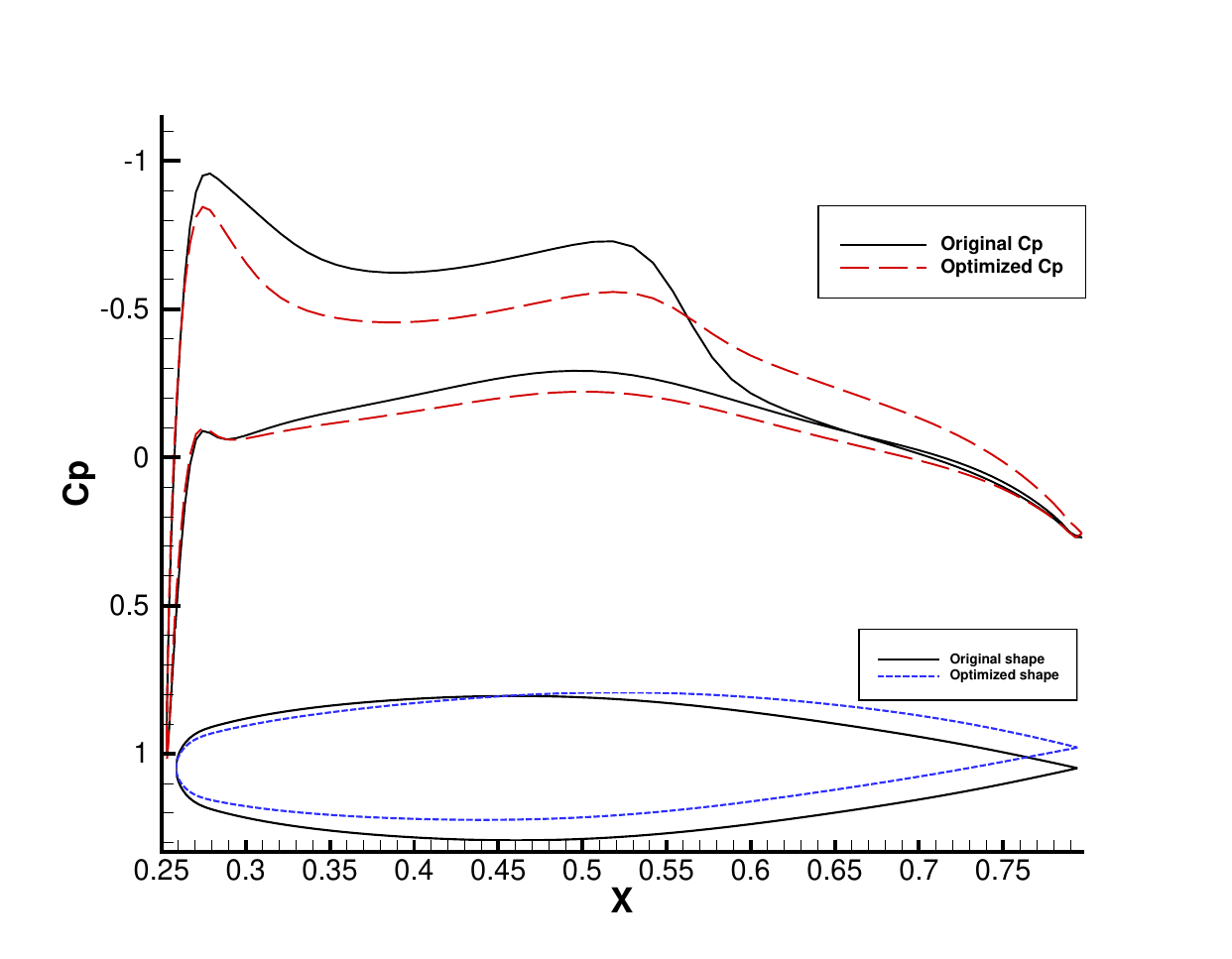}}
 \subfigure[$y=80\%\text{ span}$]{
 \includegraphics[width=8.0cm]{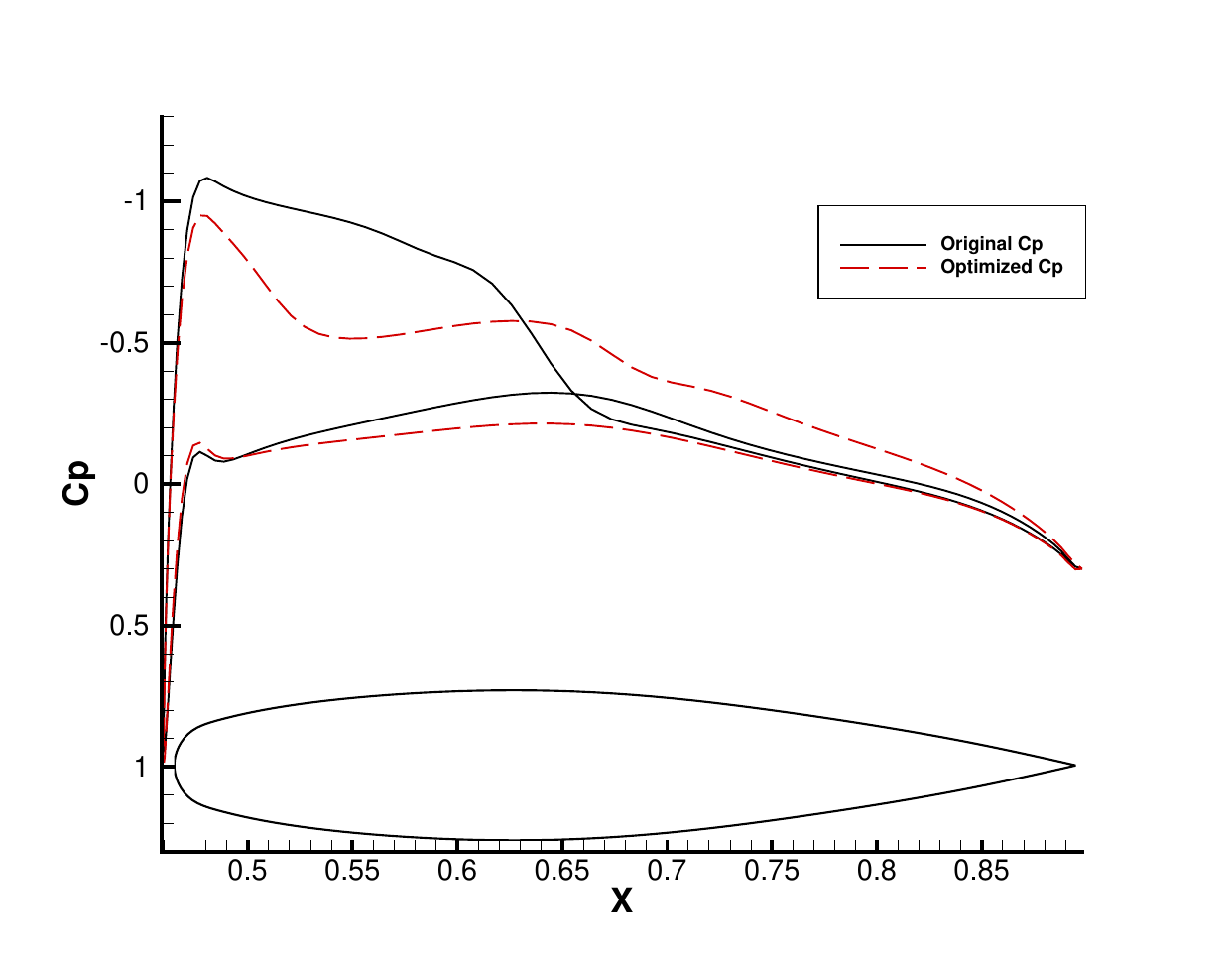}}
 \subfigure[$y=95\%\text{ span}$]{
 \includegraphics[width=8.0cm]{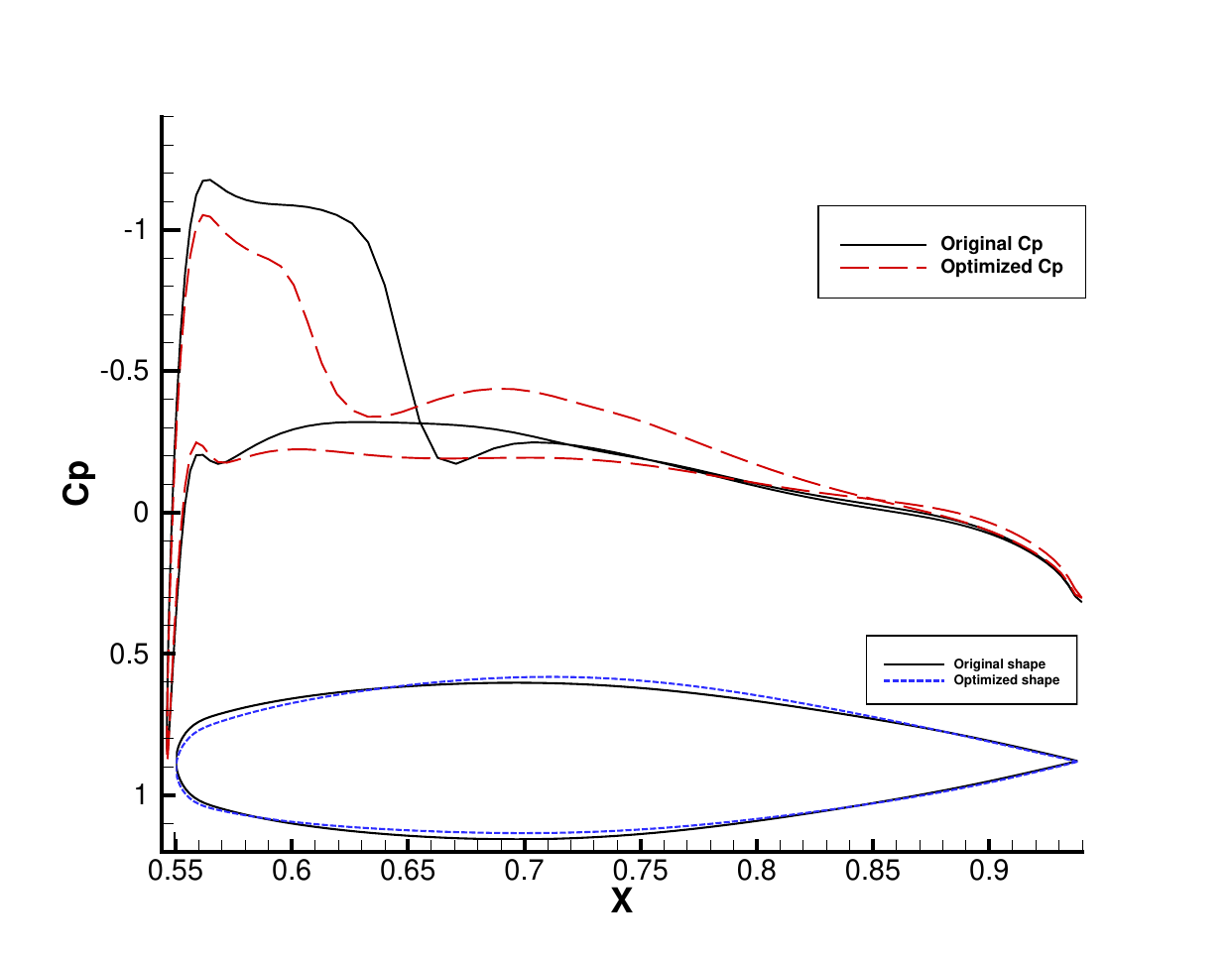}}
 \caption{Optimized shape and Cp distribution of the Onera M6 wing test case by DGp2 on mesh-1}
 \label{fig:m6-Cp2}
\end{figure}

Two different initial grids (mesh-1, mesh-2) as shown in Fig. \ref{fig:m6mesh} are used to test the performance of different CFD solvers (DGp1 and DGp2), the allocation of ierative DoFs of various CFD solvers on each grid is presented in Table \ref{tab:m6-DoFs}, and the initial $C_p$ distributions obtained by DG solvers are compared with experimental data \cite{schmitt1979pressure} in Fig. \ref{fig:m6Cp-exp}. Only the results produced by DGp1 on mesh-2 and DGp2 on mesh-1 are provided in Fig. \ref{fig:m6-flow1}-\ref{fig:m6-Cp2}. The distribution of the shock strength in Fig. \ref{fig:m6-shock1} and \ref{fig:m6-shock2} is visualized by the shock wave indicator in \cite{feng2020characteristic, feng2021characteristic2}. In addition, the iteration steps of optimization and the serial CPU run time are presented in Table \ref{tab:m6-time}, and the aerodynamic performance improvement driven by different CFD solvers are provided in Table \ref{tab:m6-opt}.

% mesh-1: 24*24*36*2+36*36*24*2; mesh-2: 32*32*64*2+64*64*32*2
\begin{table}[htbp]
\centering \caption{DoFs of various CFD solvers on each grid in the ONERA M6 test case}
\setlength{\tabcolsep}{12mm}
\begin{center}
\begin{tabular}{c|c|c}
\toprule
  & mesh-1 & mesh-2  \\ \hline
  DGp1 &2,073,600 &7,864,320 \\ \hline
  DGp2 &5,184,000 &/    \\
\bottomrule
\end{tabular}
\end{center}
\label{tab:m6-DoFs}
\end{table}

\begin{table}[htbp]
\centering \caption{Iteration step and serial CPU time (h) with/without DoFs remapping in the Onera M6 test case}
\begin{center}
\begin{tabular}{c|c|c|c}
\toprule
  & \hspace{0.6em} Iteration step \hspace{0.6em} & \small{CPU time (no remapping)} &\small{CPU time (remapping)} \\ \hline
  DGp1 on mesh-2 &20 & 241.58 &133.37  \\ \hline
  DGp2 on mesh-1 &25 & 177.47 &106.78  \\
\bottomrule
\end{tabular}
\end{center}
\label{tab:m6-time}
\end{table}

\begin{table}[htbp]
\centering \caption{Performance improvement of various CFD solvers in the Onera M6 wing test case}
\begin{center}
\begin{tabular}{c c c c c c c c c c}
\toprule
     &\footnotesize{$C_{d0}$} &\footnotesize{$C_{d1}$} &\footnotesize{$\Delta C_{d}$} &\footnotesize{$C_{l0}$} &\footnotesize{$C_{l1}$} &\footnotesize{$\Delta C_{l}$} &\footnotesize{$V_{0}$} &\footnotesize{$V_{1}$} &\footnotesize{$\Delta V$} \\
\midrule
  \footnotesize{Reference\cite{chakrabartyy1996computation,van1998discontinuous}} &\scriptsize{1.28$\sim$1.36e-2} &/ &/ &\scriptsize{2.90$\sim$2.96e-1} &/ &/ &/ &/ &/ \\
  \footnotesize{DGp1 on mesh-2} &\footnotesize{1.29e-2} &\footnotesize{9.74e-3} &\footnotesize{-24.5$\%$} &\footnotesize{2.86e-1} &\footnotesize{2.86e-1} &\footnotesize{+0$\%$} &\footnotesize{1.99e-2} &\footnotesize{1.98e-2} &\footnotesize{-0.5$\%$} \\
  \footnotesize{DGp2 on mesh-1} &\footnotesize{1.26e-2} &\footnotesize{9.60e-3} &\footnotesize{-23.8$\%$} &\footnotesize{2.84e-1} &\footnotesize{2.84e-1} &\footnotesize{+0$\%$} &\footnotesize{1.99e-2} &\footnotesize{1.92e-2} &\footnotesize{-3.5$\%$} \\
\bottomrule
\end{tabular}
\end{center}
\label{tab:m6-opt}
\end{table}

%\begin{figure}[htbp]
%  \centering
%  \includegraphics[width=6.2cm]{figure/m6-obj}
%  \caption{$C_d$ and $C_l$ variation of the Onera M6 wing test case}
%  \label{fig:obj}
%\end{figure}

It can be summarized from \ref{tab:m6-time} and Fig. \ref{fig:m6-flow1}-\ref{fig:m6-Cp2} that (a) under similar DoFs, DGp2 on mesh-1 costs less CPU time than DGp1 on mesh-2, even though it has more iteration steps; (b) DoFs remapping techniques can reduce the CPU time by around $40\%\sim50\%$; (c) even on a coarse grid, DGp2 provides similar trends of 3D shape morphing for drag minimization as compared with the DGp1 on a relatively fine gird; (d) under a small range of variation of the shape, the DGM-based CFD solvers can reduce the drag by around $24\%$ with lift coefficient and wing volume almost unchanged.

\section{Conclusion and perspectives}
\label{sec:7}
\qquad In this work, based on our recently open-sourced HODG platform, we developed a robust and efficient framework of the adjoint gradient-based aerodynamic shape optimization using DGMs as the CFD solver. It is found through theoretical analysis that DG representations show high-order strengths in the evaluation of the adjoint gradient. Compared with FVMs, DGMs can not only acquire the more accurate objective, but also resolve the modification term containing information of high-order moments of the numerical solution. The results of numerical experiments are almost consistent with the theoretical analysis, under the premise of similar computational errors and computational costs, the DGM-based flow respresentations or CFD solvers are capable of providing more precise adjoint-enabled gradient for optimization phase as compared with the FVM-based representations, and further able to explore the design space and produce more superior optimal aerodynamic shapes. Future works will focus on the extension of the governing equations (N-S, RANS equations), more aerodynamic objectives and requirements, and the high-performance parallel computation for the entire optimization chain.

%%%% Acknowledgments %%%%%%%%
\section*{Acknowledgments}
This work is supported by National Nature Science Foundation of China (No.12302380), and the National Numerical Wind Tunnel Project.

\bibliography{sample}

\end{document}